\documentclass[a4paper,12pt]{amsart}

\usepackage{amsmath,amsfonts,amssymb}
\usepackage{amsthm}
\usepackage{mathrsfs}

\usepackage{ulem}
\usepackage{graphics,graphicx,epsf}

\usepackage[usenames]{color}

\usepackage{float}
\usepackage{enumerate}
\usepackage[shortlabels]{enumitem}

\usepackage{tikz}
\usetikzlibrary{arrows}

\usepackage{csquotes}

\usepackage[
colorlinks=true,
linkcolor=blue!50!black,
urlcolor=teal!60!black,
citecolor=blue!70!black,
anchorcolor=blue!50!black,
]{hyperref}

\usepackage[left=2.22cm,right=2.22cm,top=2.9cm,bottom=2.9cm]{geometry} 

\newcommand{\cH}{\mathscr{H}}

\newtheorem{Assumption}{Assumption}[section]
\newtheorem{theorem}{Theorem}[section]
\newtheorem{lemma}[theorem]{Lemma}
\newtheorem{proposition}[theorem]{Proposition}
\newtheorem{Corollary}[theorem]{Corollary}
\newtheorem{remark}[theorem]{Remark}
\newtheorem{definition}{Definition}[section]

\numberwithin{equation}{section}

\numberwithin{equation}{section}

\title[Subquintic energy-damped wave equation]{Global Well-Posedness and Asymptotic Behavior of Energy-Damped Subquintic Wave Equations}

\author[Irena Lasiecka]{Irena Lasiecka$^{\star}$} \thanks{$^{\star}$ Research partially supported by the National Science Foundation, NSF-DMS-2205508.} \address[Irena Lasiecka]{Department of Mathematical Sciences, The University of Memphis, 3725 Norriswood Ave, Memphis, TN 38152, USA; and Systems Research Institute, Polish Academy of Sciences, Warsaw, Poland.} \email{lasiecka@memphis.edu}

\author[Yanan LI]{Yanan Li$^{\dagger}$}\thanks{$^{\dagger}$ Research partially supported by the Key Support Program of National Foreign Expert Affairs Administration: D20250158 and the Fundamental Research Funds for the Central Universities: 2026GH2401.}
\address[Yanan Li]{School of Mathematical Sciences, Harbin Engineering University, 145 Nantong Street,  Harbin, 150001 China}
\email{liyn@hrbeu.edu.cn}

\author[Vando Narciso]{Vando Narciso$^{\ddagger}$}\thanks{$^{\ddagger}$ Research partially supported by Fundect/CNPq Grant 15/2024.}
\address[V. Narciso]{Universidade Estadual de Mato Grosso do Sul, 79804-970, Dourados, MS, Brazil}
\email{vnarciso@uems.br}

\date{\today}

\begin{document}
	
	\begin{abstract}
		In this article, we consider a wave equation on a bounded three-dimensional domain, featuring a degenerate nonlocal damping mechanism that depends on the system's energy and a source term with subquintic growth. This class of models is motivated by applications in optical physics. We establish global existence in the Shatah--Strichartz sense by combining Galerkin approximations with space-time Strichartz estimates on bounded domains. A central difficulty in the long-time analysis is the interplay between three-dimensional wave dynamics, nonlocal damping, and forcing terms that are not locally Lipschitz continuous. We prove the existence of a compact global attractor and derive an upper bound for its Kolmogorov $\varepsilon$-entropy. In the non-degenerate damping case, assuming that the domain has a smooth, strictly convex boundary, we exploit the refined Strichartz estimates of Ivanovici {\it et al.}~\cite{Ivanovici-Lascar-Lebeau-Planchon-2023}, which yield, in particular, $u\in L^3_tL^{18}_x$. Under the further assumption $f\in C^2(\mathbb R)$, this improved integrability allows us 
		to establish quasi-stability of the associated dynamical system. Consequently, the global attractor has finite fractal dimension and higher regularity, and a generalized exponential attractor exists.

		
	\end{abstract}	
	\maketitle	
	
	\begin{quote}
		\tableofcontents
	\end{quote}

	\section{Introduction}
	\subsection{The model}
	This article is devoted to the analysis of the well-posedness and asymptotic dynamics of solutions $\mathbf{u}=(u,u_t)$ to the following class of wave equations:
	\begin{equation}\label{P}
		\left\{
		\begin{aligned}
			&u_{tt} - \Delta u + g\bigl(||\mathbf{u}||^2_{\mathscr{H}}\bigr)u_t + f(u) = h,
			&& \text{in } \Omega \times (0,+\infty),\\
			&u = 0,
			&& \text{on } \Gamma \times (0,+\infty),\\
			&u(0) = u_0,\quad u_t(0) = u_1,
			&& \text{in } \Omega,
		\end{aligned}
		\right.
	\end{equation}
	where $\Omega\subset\mathbb{R}^3$ is a bounded domain with smooth boundary $\Gamma=\partial\Omega$. The term
	$
	g\bigl(\|\mathbf{u}\|_{\mathscr{H}}^{2}\bigr)u_t
	$
	represents a nonlinear, nonlocal damping mechanism whose coefficient depends on the linear energy of the system through an increasing function of the form
	$
	g(s)\approx \kappa+\gamma s^{q},$ $\kappa\ge0,$ $\gamma>0,$ $q\ge 1/2.
	$
	The function $f \in C^1(\mathbb{R})$ denotes a source term satisfying the following growth condition:
	$
	|f'(u)| \le C\bigl(1 + |u|^{\rho-1}\bigr),$
	$\rho \in [1,5),$ and $h$ is an external force. The quantity
	$
	\|\mathbf{u}\|^2_{\mathscr{H}} := \|\nabla u\|^2 + \|u_t\|^2,$ $\mathbf{u}=(u,u_t),
	$
	defines the natural norm in the energy space
	$
	\mathscr{H} = H_0^1(\Omega) \times L^2(\Omega),
	$
	where $\|\cdot\|$ denotes the $L^2(\Omega)$-norm.
	
	For three-dimensional nonlinear wave equations, a source term with quintic growth corresponds to the energy-critical exponent, since the power five is the unique exponent for which the natural scaling of the wave equation leaves the energy functional invariant.
	This critical threshold separates the subcritical and supercritical regimes and plays a central role in the qualitative analysis of solutions.
	From the modeling viewpoint, fifth-order nonlinearities arise when higher-order interaction effects are retained in constitutive expansions, as in nonlinear optics, plasma dynamics, and continuum mechanics, where large-amplitude phenomena cannot be adequately described by purely cubic approximations. Quintic wave equations, which include nonlinear terms to the fifth power, are used to model complex physical phenomena such as solitons in optics and plasma physics, vibrations in materials, fluid dynamics, and the dynamics of wave propagation in various media. They are crucial for understanding and predicting the behavior of nonlinear systems and help in designing technologies in fiber optics, plasma confinement, solid state physics, and materials science [molecular vibrations and wave propagation in harmonic crystals], and nonlinear circuits. Of particular importance  are cubic-quintic media, which describe physical environments where  wave propagation is governed by two competing types of nonlinear responses. 
	The low intensity cubic reaction forces the wave to self-focus, while high intensity quintic reaction activates at intense thresholds and acts as a physical ceiling that prevents wave collapse, stabilizing localized packets into robust structures referred to as droplets. 
	Thus, the subquintic assumption $\rho \in [1,5)$ places the problem~\eqref{P} in the energy-subcritical regime while remaining close to the critical scaling that governs the most delicate analytical features of the three-dimensional wave equation.
	
	In order to ``stabilize''  the forced dynamics, an effective damping mechanism is necessary.  This results in  the second challenging feature of the problem under consideration, due to 
	the fact that the damping intensity is governed by a nonlocal coefficient depending on the system's energy. This leads to {\it nonmonotone} dissipation-in contrast with a  vast majority  of the existing literature dealing with  long time behavior, particularly in the case of  dynamics with critical forcing. 
	Thus, high power  forcing [$\rho<5$]  along with nonmonotone, energy modulated  dissipation constitutes the main  challenge of the problem under consideration. General comments presented below provide a historical perspective on these two  facets of the problem. 
	
	\subsection{General comments}
	We shall first discuss energy modulated dissipation and then we turn to quintic and subquintic forces acting upon  3 dimensional wave dynamics. 
	\subsubsection{\bf Energy modulated dissipation}
	This was  first introduced within  the framework of ordinary differential equations by Krasovskii~\cite{Krasovskii} (1963). Such models can be represented by the following class of ODEs:
	\begin{equation}\label{Krasovskii}
		\ddot{x}(t) + x(t) + \kappa(E(t)) \dot{x}(t) = 0,
	\end{equation}
	where $\kappa$ is a function with real values and
	$
	E(t) = \tfrac{1}{2}\bigl(x^2(t) + \dot{x}^2(t)\bigr)
	$
	denotes the total energy of the system.
	
	However, a physical application of Krasovskii's model is a PDE model used in aeronautics within the context of flying ``smart'' structures:  energy-dependent damping mechanisms  freely responding to vibrations, as studied by Balakrishnan~\cite{Balakrishnan} and Balakrishnan--Taylor~\cite{BalakrishnanTaylor}. In these models, the damping term is expressed directly in terms of a power law applied to the energy variable, namely
	$$
	\mathcal{D}(x(t),\dot{x}(t)) = E(t)^q \dot{x}(t), \qquad q>0.
	$$
	
	Consequently, the two formulations are consistent in the particular case
	$
	\kappa(s) = \gamma s^q
	$
	with the practical goal of bringing  given dynamics to  coherent structures-including compact attractors. From  a heuristic point of view one expects  that such a combination of forcing and dissipation  should result in  a well organized asymptotic dynamics described by an attractor defined on a  finite energy space $\cH$. 
	The issue, however, in infinite dimension is not a trivial  one, as pointed out   by the work of Chueshov in \cite[Section 5.3.3]{chueshov}, 
	where the abstract infinite-dimensional version of the Krasovskii system \eqref{Krasovskii} was considered via the following model:
	\begin{align}\label{model-chueshov}u_{tt}+Au+\kappa\left(\|A^{1/2}u(t)\|^2+\|u_t(t)\|^2\right)u_t=0,\quad u|_{t=0}=u_0,\;u_t|_{t=0}=u_1,
	\end{align}
	defined on a separable Hilbert space $H$. Here, $A$ is a positive operator with a discrete spectrum, and $\kappa(s)$ is a bounded Lipschitz scalar function on $\mathbb{R}^+$ such that
	$\kappa(1)=0$ and $\kappa(s)$ is strictly increasing for $s>1$. For the dynamical system in \eqref{model-chueshov}, the existence of a \textbf{noncompact} global attractor has been established. The main technical difficulty in proving the existence of global attractor is related to the {\bf  lack of monotonicity}  required to derive suitable estimates for the difference between two solutions, which is a crucial aspect of establishing compactness in infinite-dimensional spaces. A second difficulty is associated with {\bf  the lack of local Lipschitz property}  of the forcing terms.
	
	\subsubsection{\bf Quintic and subquintic forcings}
	In the case of  nonlinear  wave dynamics  in 3 dimensional domains,  one is interested in considering a nonlinear  feedback restoring  force $f(u)$ along with a  dissipative mechanism  $g(||u||_{\cH}) u_t$,  whose goal is to bring the forced dynamics to a coherent structure-such as an attractor. While typical  forces are of cubic structure, modern applications in nonlinear optics or plasma physics  resort to modeling solitary  pulses by quintic terms. The latter act as a balancing force [against cubic forces]  simulating artificial saturable absorbers. 
	\ifdefined\xxxxxxx
	In order to bring  asymptotically the resulting dynamics to a coherent structure, damping mechanism is necessary. And typical damping will be modulated by the strength of the total energy. This leads to a  nonlocal  and non-monotone effects described by a  nonlinear [power like] function  $g(||u(t)||_{\cH} )$m  a heuristic point of view one expects  that such combination of forcing and dissipation  should result in  a well organized asymptotic dynamics described by an attractor defined on a  finite energy space $\cH$. 
	Thus, no wonder the above  general  problem (\ref{P})  has attracted considerable attention in mathematical literature. Its physical relevance combined with mathematical challenge provides a rich .  Fro
	landscape  for investigations and development of new mathematical tools.  The issue was sparked  by work of Chueshov in \cite[Section 5.3.3]{chueshov} where
	the abstract infinite-dimensional version of the Krasovskii system \eqref{Krasovskii} was considered in \cite[Section 5.3.3]{chueshov} via  the following  model:
	\begin{align}\label{model-chueshov}u_{tt}+Au+\kappa\left(\|A^{1/2}u(t)\|^2+\|u_t(t)\|^2\right)u_t=0,\quad u|_{t=0}=u_0,\;u_t|_{t=0}=u_1,
	\end{align}
	defined on  a separable Hilbert space $H.$ Here $A$ is a positive operator with a discrete spectrum, $\kappa(s)$ is a scalar bounded Lipschitz function on $\mathbb{R}^+$ such that
	$\kappa(1)=0$ and $\kappa(s)$ is strictly increasing for $s>1$. The dynamics of problem \eqref{model-chueshov} were considered and the existence of a {\bf noncompact} global attractor has been established. The existence of a compact global attractor remained an open question, due to the difficulty in treating nonlocal dissipation that is  {\it not monotone}. In other words, analyzing the trajectories  difference in this term generates a monotone dissipative term plus a term that needs to be controlled, but lacks compactness in the topology of the problem's weak phase space. Thus  an ultimate  compactness of trajectories in the presence of non-monotone energy modulated damping can not be taken for granted. 
	The above example  along with   the fundamental  applicability  to  aeroelasticity[  large space structures -Balakrishnan- Taylor damping]  sparked great  interest in studying  dynamical properties of  a wave/plate
	\fi
	
	From the  mathematical point of view, the main obstacles are twofold.  The first lies at the level of  the restoring force $f(u)$, and the second lies at the level of  dissipation. The restoring force above cubic  {\bf loses its local Lipschitz character. The dissipation is no longer monotone} -the latter being a main driver in  long time asymptotic analysis. Of course, these difficulties can be addressed by incorporating various regularizing mechanisms into the dynamics, such as strong damping \cite{Carvalho-Cholewa,Pata-Zelik,Yang-Sun}. However, when dealing with a ``pure''  unregularized model of  wave dynamics  above the local Lipschitz level, one needs to resort to more refined  tools, such as Strichartz-Struve estimates [SS]. These estimates originated in the work of Strichartz
	\cite{Strichartz-70,Strichartz-77} and were subsequently developed
	and applied to critical nonlinear wave equations
	\cite{Struwe,Shatah-Struwe-93,Struwe-94}. These were first established for the full space \cite{Shatah-Struwe-93} and later extended to bounded domains
	\cite{Blair,Burq-Lebau-Planchon,Burq-Planchon}.
	
	SS estimates allow one to trade spatial integrability for temporal integrability. Finite energy lies at the endpoint of the scale. While Strichartz estimates provide the tool for establishing local well-posedness of solutions, the word ``local'' is essential. At the critical quintic level, these estimates do not provide any control of long-time behavior-as the estimates themselves depend on the solution rather than the norm \cite[Section 3]{KSZ}. Thus, the question of whether Strichartz estimates survive asymptotically in time is still an open problem [in contrast to energy estimates]. This has been a main predicament in dealing with long-time dynamics, starting with the global existence of solutions. In the case of linear [linearly bounded] dissipation, the arguments based on non-concentration of Strichartz norms proved successful in resolving the global existence problem of Shatah-Struwe solutions \cite[Section A]{Savostianov-thesis}. Since the argument is based on a contradiction principle, no estimates are available to extend the validity of the SS estimates obtained for local solutions. Instead, in the subquintic case, control of SS estimates over  long-time horizons can be achieved \cite[Section 4]{KSZ},  as detailed in the present work. However, {\it monotonicity of the dissipation} plays a critical role in the arguments. 
	This brings us to the core of the problem.

	
	The above considerations, together with  fundamental applications to
	aeroelasticity [in particular, large space structures with Balakrishnan--Taylor damping],
	have sparked great interest in studying the dynamical properties of 
	wave/plate systems equipped with this type of non-monotone dissipation
	and nonlinear restoring  force:
	\begin{align}\label{model}
		u_{tt}+Au+g\left(\|A^{1/2}u(t)\|^2+\|u_t(t)\|^2\right)u_t
		+f(u)=0,\quad u|_{t=0}=u_0,\;u_t|_{t=0}=u_1,
	\end{align}
	where the function $g(s)$ is of rational/polynomial type,
	$g(s)=\kappa+\gamma s^q$, determined by values of $q$. Our problem  (\ref{P}) falls into this category with $ A = -\Delta $ under Dirichlet boundary conditions and $ q \in [1/2,2]$. Moreover, the restoring force $f(u)$ is {\it subquintic} in contrast to  the majority of the literature where locally Lipschitz forces  have been considered. For this model we shall establish the existence of  a well-posed dynamical system  which admits a  {\bf global attractor}  in the strong topology of the phase space $\cH$. Under  an additional assumption $\kappa >0 $ the said  dynamical system is quasi-stable, hence  it  is {\bf  smooth, finite-dimensional, and admits an exponential attractor}.
	
	\subsubsection{\bf What is new?}  
	The literature regarding the general model (\ref{model}) is large and relatively recent.  Thus, it is appropriate to elucidate  novel and original aspects of the present work. The vast majority of papers dealt with either {\it regularized} models [e.g., by introducing strong damping] \cite{Carvalho-Cholewa,Pata-Zelik,Yang-Sun} 
	or with {\it locally Lipschitz internal forces} \cite{BJN-AA,Freitas-Narciso,GJNV,GJLN,Narciso-Nodea,Zhou-Yang-2026}.
	In order to focus our exposition, we will not  discuss these works  extensively  [see \cite{GJLN} for more information], but rather provide a perspective on {\it higher order forcing subject to energy modulated damper}, which commands a new set of tools and mathematical techniques.

	{\it Results  with forcings which are not locally Lipschitz} include  \cite{SavostianovZelik2014a,SavostianovZelik2014b,Savostianov-thesis,YZZT,YTZ,ZLZ,LasieckaNarciso2025}. However, these works consider only  monotone damping modulated by kinetic energy. This leads to a different type of analysis. While  SS estimates are used for both well-posedness and long time behavior, the monotonicity of the dissipation permits a   direct  treatment of the difference of two solutions. 
	Indeed, the monotonicity of the dissipation, together with  the SS estimates  for the nonlinear forced term, reduces the problem  at the technical level to  the  analysis  developed   for dissipative systems in the quintic case \cite{KSZ,Savostianov-thesis}, which relies on classical SS estimates in $L^4(L^{12})$ or $L^5(L^{10})$.
	This is not the case in the non-monotone case of energy modulated damping, where more refined tools  such as specially designed multipliers and nonstandard SS estimates are employed.
	New classes of SS estimates are called for and new constructions of multipliers are needed. 
	\ifdefined\xxxxxx


	\begin{align}\label{ADE-model}u_{tt}+\Delta^2 u-\kappa \Delta u+\gamma\left[\,\|(-\Delta)^{\alpha}u(t)\|^2+\|u_t(t)\|^2\,\right]^qu_t+f(u)=h,\end{align}
	where $\Omega$ is a bounded domain of $\mathbb{R}^n$, $\gamma>0$, $\kappa\ge 0$, $q\ge1/2$, $\alpha \in [0,1]$, $f$ a nonlinearity with growth exponent $\rho\le \frac{4}{n-4}$ if $n\ge 5$, and $\|\cdot\|$ represents the $L^2$-norm. By restricting  $\alpha<1$, along with the range of parameter $ q \in [1,1 +\epsilon] $  the authors prove that the dynamical system generated by the weak solutions of the problem is asymptotically compact, and therefore, the existence of a compact global attractor for the problem \eqref{ADE-model} was then established.
	Subsequently the restrictions such as $\alpha < 1 $   have been removed   \cite{ccc} , the range of parameters $q$  enlarged  and the Krasowski model with locally Lipschits restoring force $f(u)$ became treatable. This brings us to the core of the problem studied in this work - which is to account for not locally Lipschitz  restoring force with non-monotone energy dependent damping.

	Recently, Gomes {\it et al.} \cite{GJLN}, using new nonlinear multipliers combined with the multiplier technique used by Chueshov and Lasiecka in \cite{chueshov-lasiecka-2008}, managed to overcome the difficulty of the case $\alpha=1$ and proved that the dynamical system is also asymptotically compact for $\alpha=1$, and therefore, the existence of a compact global attractor for the problem \eqref{ADE-model} in the critical case $\alpha=1$ has been achieved. Also recently, Narciso \cite{Narciso-Nodea} considered the Krasovskii model \eqref{model-chueshov} with $\kappa(s)\sim s^q$ and an additional source term $f$. Using nonlinear multipliers as in \cite{GJLN}, the existence of a {\bf compact global attractor} was established in the regime $q\in[1,1+\varepsilon)$, with $\varepsilon>0$ sufficiently small. Moreover, an upper bound for the Kolmogorov $\varepsilon$-entropy of the global attractor was derived. Under the additional assumption that $\kappa'(0)>0$, the finite-dimensionality and additional regularity of the attractor were also established. More recently, Zhou and Yang \cite{Zhou-Yang-2026} extended the results of Narciso \cite{Narciso-Nodea} to the range $q\in(0,2]$.
	
	An application of Krasovskii's model in the context of wave equations was considered by Tang \textit{et al.} \cite{TYXZ-DCDS}. More specifically, in a bounded domain $\Omega\subset \mathbb{R}^3$, the following wave equation endowed with Dirichlet boundary conditions was proposed
	\begin{align}\label{model-TYXZ}
		u_{tt}-\Delta u+\left[\,\|\nabla u(t)\|^p+\|u_t(t)\|^p\,\right]u_t+f(u)=0.
	\end{align}
	In \eqref{model-TYXZ}, the exponent $p\in (1,2]$ and the nonlinearity $f\in C^2(\mathbb{R})$ satisfies the growth condition $|f''(s)|\le C(1+|s|)$, together with the dissipativity assumption $\liminf_{|s|\to +\infty} f'(s)>-\lambda_1$, where $\lambda_1$ denotes the first eigenvalue of the Laplacian $-\Delta$. Under these assumptions, the authors establish the existence of a global attractor and, moreover, develop a refined analysis near the degenerate regime, which yields the finite dimensionality of the attractor. We emphasize that the Krasovskii damping in this model is sublinear with respect to the energy, since it behaves as
	\[
	\kappa(E(t)) \approx c\,E(t)^{p/2}, \qquad p\in(1,2],
	\]
	which plays a crucial role in the treatment of the degeneracy of the damping mechanism. Building on this framework, the non-autonomous version was recently studied by Bezerra \textit{et al.} \cite{BJN-AA}, where the following model was considered:
	\begin{align}\label{non-autonomous}
		u_{tt}-\Delta u+\gamma(t)\left[\,\|\nabla u(t)\|^2+\|u_t(t)\|^2\,\right]^q u_t+f(u)=0.
	\end{align}
	Here $f(u)\sim u^3$ in the case $n=3$, and $\gamma:\mathbb{R}\to \mathbb{R}$ is bounded, globally Lipschitz, and satisfies $0<\gamma_0\le \gamma(t)\le \gamma_1$ for some constants $\gamma_1\ge \gamma_0>0$. For $q\in [1,1+\epsilon)$, with $\epsilon>0$ sufficiently small, the authors employ nonlinear multipliers as in \cite{GJLN} to prove the existence of a pullback attractor for problem \eqref{non-autonomous}.
	
	\medskip
	Wave equations related to the Krasovskii and Balakrishnan-Taylor models, {\bf without the potential energy term} in the dissipative coefficient, were studied by Aloui {\it et al.}~\cite[Section~8, Eqs.~(8.1) and~(8.8)]{ABH}. In this context, the corresponding dissipative mechanism is referred to as {\bf averaged damping}.
	More recently, Yan {\it et al.}~\cite{YTZ} investigated the following wave equation with averaged damping and Dirichlet boundary conditions:
	\begin{align}\label{model_averaged_damping-1}
		u_{tt} - \Delta u + k\|u_t(t)\|^p u_t + f(u) = g(x),
		\quad \text{in } \Omega \times \mathbb{R}^+.
	\end{align}
	Here, $\Omega \subset \mathbb{R}^3$ is a bounded domain with smooth boundary,
	$k, p > 0$ are fixed constants,
	$g \in L^2(\Omega)$ is an external force,
	and the source term $f$ satisfies the following subquintic growth condition
	$|f'(s)| \le C(1 + |s|^{4 - \kappa})$,
	for $s \in \mathbb{R}$ and $\kappa \in (0,4]$.
	Using Strichartz estimates adapted to bounded domains,
	the authors established the existence of global Shatah--Struwe solutions
	and proved that the dynamical system generated by these solutions
	possesses a global attractor. The non-autonomous version of problem~\eqref{model_averaged_damping-1},
	with $g = g(x,t)$, was considered by Zhao {\it et al.}~\cite{YZZT}.
	The authors introduced the concept of a uniform $\varphi$-attractor
	and showed that the Shatah--Struwe family of solution processes
	admits a uniform polynomial attractor.
	
	To complete the contextualization of this class of nonlocal problems, we note that the most challenging models related to the energy-damped framework are those in which the {\bf kinetic energy term is absent} from the nonlocal damping coefficient. In such cases, standard multiplier techniques are insufficient to recover the kinetic energy term and to derive the inequalities that govern the decay of solutions. In this setting, Cavalcanti {\it et al.}~\cite{Cavalcanti} studied the following equation
	in a bounded domain~$\Omega \subset \mathbb{R}^n$, subject to Dirichlet boundary conditions:
	\begin{equation}\label{model-cavalcanti}
		u_{tt} - \Delta u + \|\nabla u(t)\|^2 u_t + f(u) = 0,
	\end{equation}
	where the nonlinearity $f$ has growth exponent $2/(n-2)$ for $n \ge 3$.
	To overcome the difficulty of applying direct multiplier methods, the authors employed contradiction arguments to prove that the energy of the model \eqref{model-cavalcanti} decays to zero for initial data in bounded sets of
	$(H^2(\Omega) \cap H^1_0(\Omega)) \times H^1_0(\Omega)$. By contrast, results concerning the asymptotic behavior of the nonlocal degenerate model \eqref{model-cavalcanti} in the weak topological space $H^1_0(\Omega) \times L^2(\Omega)$, under weak initial data, remain an open question.

	\subsection{A brief overview of the literature on quintic waves}
	It is well known that the global well-posedness of the quintic wave equation in the whole space \(\mathbb{R}^3\) was first established by Struwe~\cite{Struwe} for radially symmetric solutions. Later, Grillakis extended this result to the nonradial setting for smooth initial data. These works rely essentially on explicit representation formulas for solutions of the wave equation in \(\mathbb{R}^3\), together with the Morawetz--Pohozaev identity. Subsequently, Shatah and Struwe~\cite{Shatah-Struwe-93} rigorously established the existence and uniqueness of global solutions in the natural energy phase space, with further refinements developed in~\cite{Struwe-94}. The resulting class of Shatah--Struwe solutions possesses additional space--time regularity beyond the standard energy framework, providing a suitable setting for the analysis of the long-time behavior of quintic wave equations.
	
	To the best of our knowledge, the first works devoted to the asymptotic dynamics of dissipative wave equations with quintic critical nonlinearities in bounded three-dimensional domains are due to
	\cite{Carvalho-Cholewa,Pata-Zelik,Yang-Sun}, where strongly damped models of the form
	\begin{equation}\label{strongly-damped}
		u_{tt} - \Delta u - \Delta u_t + f(u) = h
	\end{equation}
	were investigated.
	
	In particular, the existence of a compact global attractor for bounded smooth domains \(\Omega \subset \mathbb{R}^n\), under nonlinearities with critical growth
	\[
	f(u)\sim |u|^{\frac{n+2}{n-2}}
	\]
	(and hence \(f(u)\sim |u|^5\) when \(n=3\)),
	was extensively studied by Carvalho and Cholewa in
	\cite[Section~4.3]{Carvalho-Cholewa}, where the analysis relies on Alekseev's nonlinear variation-of-constants formula. Later, Pata and Zelik~\cite{Pata-Zelik} established the existence of compact global attractors with optimal regularity properties.
	
	Further developments concerning global well-posedness, dissipativity, and the existence of attractors for wave equations with quintic nonlinearities and fractional structural damping were obtained in
	\cite{KSZ,SavostianovZelik2014a,SavostianovZelik2014b,Savostianov-thesis}. More precisely, the authors considered equations of the form
	\begin{equation}\label{linear-structural-damping}
		u_{tt} - \Delta u + \gamma(-\Delta)^{\theta}u_t + f(u) = h,
	\end{equation}
	posed in bounded domains \(\Omega\subset\mathbb{R}^3\), with
	$
	\theta =0$ $\text{or}$ $\theta=\frac12.
	$ For the case \(\theta=0\), the works
	\cite[Chapter~2]{SavostianovZelik2014a} and \cite{KSZ}
	combine Strichartz estimates in bounded domains with suitable extensions of the Morawetz--Pohozaev identity in order to establish the global existence of Shatah--Struwe solutions and the existence of a compact global attractor for the associated dynamical system.
	
	For the structural damping case \(\theta=\frac12\),
	the works
	\cite[Chapter~4]{Savostianov-thesis} and
	\cite{SavostianovZelik2014a}
	derive additional regularity of energy solutions through the construction of a suitable Lyapunov-type functional. This approach yields global well-posedness, dissipativity, and the existence of smooth global and exponential attractors with finite Hausdorff and fractal dimensions. In addition, the hidden regularizing effects of the structural damping term were analyzed in detail in~\cite{SavostianovZelik2014b}. Owing to this additional smoothing property, the existence of smooth attractors follows naturally in the case \(\theta=\frac12\).
	
	Building upon the works
	\cite{KSZ,SavostianovZelik2014a,SavostianovZelik2014b,Savostianov-thesis},
	and motivated by the successful application of Strichartz estimates
	to quintic wave equations, several authors have investigated the extension
	of these techniques to dissipative wave models compatible with dispersive
	estimates. In this direction, Cavalcanti {\it et al.}~\cite{Cavalcanti}
	established the exponential decay of the energy for the quintic wave
	equation with a localized damping term of the form \(a(x)u_t\).
	Combining Strichartz estimates with the Galerkin approximation method,
	the authors proved the global existence of Shatah--Struwe solutions,
	while the exponential decay of the energy was obtained through a contradiction argument.
	
	Subsequently, Zhou {\it et al.}~\cite{ZLZ} investigated the quintic wave equation with a nonlocal averaged damping mechanism, namely model~\eqref{model_averaged_damping-1}, involving a nonlinear dissipation of the form
	\[
	\|u_t(t)\|^{q}u_t(t).
	\]
	The authors established the existence and described the structure of the weak, strong, and exponential attractors associated with the generated dynamical system. Moreover, the monotone character of the dissipative coefficient allowed for an effective use of Strichartz estimates in the analysis.
	
	For classes of quintic wave equations whose dissipative mechanisms are not compatible with Strichartz estimates, different techniques are required. In this direction, Lasiecka and Narciso~\cite{LasieckaNarciso2025} studied a quintic wave equation in a bounded three-dimensional domain with nonlinear damping of the form
	\[
	g(u_t)\sim |u_t|^p u_t,
	\qquad p=4.
	\]
	The authors proved the existence of a compact global attractor and showed that, under the additional assumption \(g'(0)>0\), the associated dynamical system is quasi-stable. As a consequence, the attractor has finite fractal dimension, possesses additional regularity, and admits an exponential attractor.
	
	In this energy-critical setting, the analysis relies on refined energy-dissipation mechanisms, the derivation of suitable energy identities for weak solutions, appropriate modifications of Ball's method, and the general theory of quasi-stable dynamical systems. Related dissipative models have recently been investigated in
	\cite{LiuMengHanZhang2023,LiuMengZhang2025,Zhou2025}.
	\fi 
	
	\subsubsection{\bf How do we handle the problem?}
	
	It is important to emphasize that the analysis of global solutions and the long-time dynamics of problem~\eqref{P} cannot be carried out by means of developed  techniques such as those employed in
	\cite{Carvalho-Cholewa,LasieckaNarciso2025,LiuMengHanZhang2023,LiuMengZhang2025,Pata-Zelik,Yang-Sun}. The main difficulty stems from the nonlinear and nonlocal character of the damping term
	\[
	g\bigl(\|\mathbf u\|_{\mathscr{H}}^2\bigr)u_t,
	\]
	which does not provide sufficient regularization to control the higher-order terms generated by the source nonlinearity \(f\). Consequently, dispersive tools, and in particular Strichartz estimates, become essential in the analysis of problem~\eqref{P}.
	
	In this direction, the fundamental reference is the work of Kalantarov \textit{et al.}~\cite{KSZ}, concerning the three-dimensional wave equation with linear damping \(\gamma u_t\) and source nonlinearities of quintic and subquintic growth. In the critical quintic case, the authors established the global existence of Shatah--Struwe solutions by combining Strichartz estimates adapted to bounded domains with a non-concentration argument based on Morawetz-type identities.
	
	A major difficulty in the critical regime is that the local existence time
	$
	T=T(\mathbf u_0)
	$
	depends, \textit{a priori}, on the particular initial data
	$
	\mathbf u_0\in\mathscr{H}
	=H_0^1(\Omega)\times L^2(\Omega),
	$
	rather than solely on its energy norm. As a consequence, the corresponding Strichartz estimates are not uniform on bounded subsets of the phase space, and one cannot directly construct a strongly continuous dynamical semigroup associated with the equation.
	
	To overcome this obstruction, Kalantarov \textit{et al.}~\cite{KSZ} combined Strichartz theory with the trajectory dynamical systems framework. Local Shatah--Struwe solutions were first constructed; these solutions are unique, satisfy the exact energy identity, and can be extended globally through a non-concentration mechanism. In parallel, global Galerkin solutions were obtained as limits of finite-dimensional approximations. Although these solutions satisfy uniform dissipative estimates at the energy level, uniqueness is not available in this class. The long-time analysis is therefore performed in the space of complete Galerkin trajectories, which is weakly closed and possesses a gradient structure induced by the dissipative energy functional. This approach yields the existence of a weak trajectory attractor. A further backward regularity argument then allows one to recover uniqueness and the exact energy identity on complete bounded trajectories, which makes it possible to apply Ball's energy method and ultimately obtain a compact global attractor.
	
	However, the extension of the approach developed in~\cite{KSZ} to nonlinear nonlocal dissipative models of the form considered here remains largely open in the critical quintic regime. Indeed, the Morawetz-type arguments employed in~\cite{KSZ} strongly depend on the particular structure of the linear damping and do not directly apply to nonlinear degenerate dissipation mechanisms.
	
	In contrast, the subquintic regime exhibits a substantially more favorable structure. In this setting, Strichartz estimates are uniform on bounded subsets of the energy space, and the maximal existence time depends only on the energy norm of the initial data. This allows one to combine dispersive estimates with suitable dissipative a priori bounds in order to establish global well-posedness and generate a strongly continuous dynamical semigroup in the natural phase space, without resorting to the trajectory formalism.
	
	Building upon this well-posedness theory, we prove the existence of a compact global attractor for problem~\eqref{P}. Dissipativity is obtained through Nakao's method, while asymptotic smoothness follows from a compensation argument based on contractive functionals. In the non-degenerate case, assuming that $\Gamma$ is smooth and strictly convex, we exploit the refined Strichartz estimates of Ivanovici {\it et al.}~\cite{Ivanovici-Lascar-Lebeau-Planchon-2023}, which yield, in particular, $u\in L^3(0,T;L^{18}(\Omega))$. Under the further assumption that $f\in C^2(\mathbb R)$, this improved integrability allows us to close the estimates for differences of solutions and establish quasi-stability of the associated dynamical system. Consequently, the global attractor has finite fractal dimension and additional regularity, and a generalized exponential attractor exists.
	
	Therefore, although the subquintic regime preserves the essential nonlinear and nonlocal features of the model, it provides a  {\it robust}  framework for the study of long-time dynamics, including global well-posedness, asymptotic compactness, finite-dimensionality of the attractor, and the existence of both global and exponential attractors.

	\section{Well-posedness}
	\subsection{Phase space and assumptions}
	We begin by recalling some fundamental concepts and key properties of the Laplacian operator on
	$L^{2}(\Omega)$ that will be relevant for the subsequent analysis. Let $L^{2}(\Omega)$ denote the
	usual Lebesgue space, endowed with the scalar product $(\cdot,\cdot)$ and the associated norm
	$\|\cdot\|$. We consider the negative Laplace operator $-\Delta$ with domain
	$
	D(-\Delta)=H^{2}(\Omega)\cap H^{1}_{0}(\Omega).
	$
	It is well known that $-\Delta$ is a positive self-adjoint operator with a compact inverse.
	Consequently, there exists a complete orthonormal family $\{\omega_j\}_{j\in\mathbb{N}}\subset
	L^{2}(\Omega)$ of eigenfunctions of $-\Delta$ such that
	\[
	-\Delta \omega_j=\lambda_j \omega_j,\qquad j=1,2,\ldots,
	\]
	where the corresponding eigenvalues satisfy
	\[
	0<\lambda_1\le \lambda_2\le \cdots,\qquad \lambda_j\to +\infty \quad \text{as } j\to +\infty.
	\]
	
	In this paper, we investigate the dynamics of solutions to problem~\eqref{P} in the energy phase space
	\begin{align*}
		\mathscr{H} := H^{1}_{0}(\Omega)\times L^{2}(\Omega),
	\end{align*}
	equipped with the natural (graph) norm
	$$
	||\mathbf{u}||^2_{\mathscr{H}}
	:= \|\nabla u\|^{2}_{L^{2}(\Omega)} + \|v\|^{2}_{L^{2}(\Omega)},
	\qquad \mathbf{u}=(u,v).
	$$

	We next establish the well-posedness of the nonlinear problem \eqref{P} under the following assumptions on the nonlinearities $g$ and $f$ and the external force $h$.
	\begin{Assumption} \label{Assumption}
		\begin{itemize}
			\item The damping coefficient $g$ is an increasing function of the form
			\begin{align}
				\label{hyp-g}
				g(s)=\kappa+\gamma s^q,
				\qquad
				\kappa\ge0,\quad \gamma>0,\quad q\ge\frac12,
				\qquad s\in \mathbb{R}^+:=[0,\infty).
			\end{align}
			
			\item The nonlinearity $f\in C^1(\mathbb{R})$ satisfies $f(0)=0$ and the
			polynomial growth condition
			\begin{align}
				\label{hyp_f'}
				|f'(s)|\le C_{f'}(1+|s|^{\rho-1}),
				\qquad s\in\mathbb{R},
			\end{align}
			for some positive constant $C_{f'}$ and  some $\rho\in[1,5)$. Moreover, $f$
			satisfies the dissipativity condition
			\begin{align}
				\label{hyp-inf-f}
				\liminf_{|s|\to+\infty} f'(s)>-\lambda_1,
			\end{align}
			where $\lambda_1$ denotes the first eigenvalue of the Dirichlet
			Laplacian $-\Delta$.
			
			\item The external force satisfies
			\begin{align*}
				h\in L^2(\Omega).
			\end{align*}
		\end{itemize}
	\end{Assumption}
	\begin{remark}
		From Assumption \eqref{hyp_f'} and the Mean Value Theorem, it follows that there exists $C_f>0$ such that
		\begin{align}
			\label{hyp_f}|f(s)|\le C_f(1+|s|^{\rho}),\quad \forall s\in \mathbb{R}.
		\end{align}
		Moreover, Assumption \eqref{hyp-inf-f}  implies that
		\begin{align}
			\label{hyp-f-2}
			-c_f-\frac{\theta}{2}|s|^2\le F(s)\le f(s)s+\frac{\theta}{2}|s|^2+c_f,\quad \forall s\in \mathbb{R},
		\end{align}
		for some $c_f>0$ and $\theta\in [0,\lambda_1)$, where $F(s) := \int_{0}^{s} f(\tau)\,d\tau.$
	\end{remark}
	
	In the sequel, we will use
	\begin{equation}\label{def-w} \omega:=1-\frac{\theta}{\lambda_1}>0.
	\end{equation}
	
	\begin{remark}\label{lowerbounded}
		From the conditions \eqref{hyp_f'} and \eqref{hyp-inf-f}, it follows that
		\begin{equation}\label{lower1}
			f'(s)\geq -K_f, \qquad s\in\mathbb{R},
		\end{equation}
		for some positive constant $K_f$. Indeed, by the dissipativity condition \eqref{hyp-inf-f}, there exist $\mu \in [0,\lambda_1)$ and a sufficiently large $N>0$ such that
		\begin{equation*}
			f'(s) \ge -\mu, \quad \text{for } |s|>N.
		\end{equation*}
		On the other hand, from \eqref{hyp_f'}, we have
		\begin{align*}
			|f'(s)| \le C_{f'}(1+|N|^{\rho-1}), \quad \text{for } |s| \le N.
		\end{align*}
		Taking $K_f = \max\{\mu, C_{f'}(1+|N|^{\rho-1})\}$ yields \eqref{lower1}.
	\end{remark}


The energy functional $E(\mathbf{u}(t))$ associated with problem~\eqref{P} is defined by
\begin{align}\label{functional-energy}
	E(\mathbf{u}(t))
	:= \frac{1}{2}\|\mathbf{u}(t)\|_{\mathscr{H}}^{2}
	+ \int_{\Omega} F(u)\,dx-\int_{\Omega}hu\,dx.
\end{align}

\subsection{Strichartz estimates and the definition of Shatah--Struwe solutions}

Global well-posedness for problem~\eqref{P} in the subquintic regime relies on Strichartz estimates for wave equations on bounded domains; see \cite{Blair, Burq-Lebau-Planchon, Burq-Planchon, Ivanovici-Lebeau-Planchon-2014,Ivanovici-Lascar-Lebeau-Planchon-2023}. When the energy-dependent damping coefficient may degenerate, we use the estimates of Blair {\it et al.}~\cite{Blair} to obtain the space-time regularity required for the construction of Shatah--Strichartz solutions. In the non-degenerate case, under the additional assumption that $\Gamma$ is smooth and strictly convex, the refined estimates of Ivanovici {\it et al.}~\cite{Ivanovici-Lebeau-Planchon-2014,Ivanovici-Lascar-Lebeau-Planchon-2023} provide an enlarged range of admissible pairs, including $(r,s)=(3,18)$. This additional integrability will be used later in the quasi-stability analysis.

We begin with the associated linear problem
\begin{equation}\label{wave-linear}
	\left\{
	\begin{aligned}
		u_{tt}-\Delta u &= G(t),
		&& (x,t)\in\Omega\times\mathbb{R}^{+},\\
		u|_{\partial\Omega} &= 0,
		&& t\in\mathbb{R}^{+},\\
		\mathbf{u}(0) &= (u_0,u_1)=:\mathbf{u}_0,
	\end{aligned}
	\right.
\end{equation}
where $\Omega\subset\mathbb{R}^3$ is a bounded domain with smooth
boundary $\Gamma=\partial\Omega$ and
$G:\mathbb{R}^{+}\to L^2(\Omega)$ is given.

We shall use the following standard notion.

\begin{definition}
	A pair
	$
	(r,s)\in(2,\infty]\times[2,\infty)
	$
	is called an {\bf admissible pair} if
	\begin{equation}\label{admissible-pair}
		\frac{1}{r}+\frac{3}{s}=\frac{1}{2}.
	\end{equation}
\end{definition}

Considering $\overline{\Omega}$ as a compact Riemannian manifold with
boundary, endowed with the Euclidean metric, the Strichartz estimates
of Blair {\it et al.}~\cite[Corollary 1.2]{Blair} yield the following
result.

\begin{proposition}[Strichartz estimate]
	\label{Prop-Strichartz-estimate}
	Let $\mathbf{u}_0\in\mathscr{H}$ and
	$G\in L^1(0,T;L^2(\Omega))$, and let
	$\mathbf{u}(t)=(u(t),u_t(t))$ be a solution of~\eqref{wave-linear}
	such that
	\begin{equation*}
		\mathbf{u}\in C([0,T];\mathscr{H}).
	\end{equation*}
	If $(r,s)$ is an admissible pair satisfying
	\begin{equation}\label{BSS-scaling}
		\frac{3}{r}+\frac{2}{s}\leq1,
	\end{equation}
	then
	$
	u\in L^r(0,T;L^s(\Omega))
	$
	and
	\begin{equation}\label{Strichartz-estimate}
		\|u\|_{L^r(0,T;L^s(\Omega))}
		\leq
		C_T
		\left(
		\|\mathbf{u}_0\|_{\mathscr{H}}
		+\|G\|_{L^1(0,T;L^2(\Omega))}
		\right),
	\end{equation}
	where $C_T>0$ depends on $T$ and $\Omega$, but is independent of
	$\mathbf{u}_0$ and $G$.
\end{proposition}

The admissibility relation~\eqref{admissible-pair}, together with
the additional restriction~\eqref{BSS-scaling}, determines the range
of Strichartz exponents available for a general smooth boundary.
Indeed,
\begin{equation*}
	\frac{3}{r}+\frac{2}{s}
	=
	\frac{3}{2}-\frac{7}{s},
\end{equation*}
and hence~\eqref{BSS-scaling} is equivalent to \(s\leq14\).
Combining this restriction with the admissibility relation, we obtain
\begin{equation*}
	6\leq s\leq14,
	\qquad
	r=\frac{2s}{s-6},
\end{equation*}
where \(r=\infty\) when \(s=6\). Equivalently,
\begin{equation*}
	\frac72\leq r\leq\infty.
\end{equation*}
Thus, Proposition~\ref{Prop-Strichartz-estimate} includes, for
example, the pairs \((r,s)=(4,12)\) and \((r,s)=(5,10)\).
Notice, however, that the additional restriction~\eqref{BSS-scaling}
excludes admissible pairs with \(s>14\).

When \(\Omega\) has a smooth, strictly convex boundary, the Strichartz
estimates established by Ivanovici \textit{et al.}~
\cite{Ivanovici-Lebeau-Planchon-2014} and their subsequent refinement
by Ivanovici \textit{et al.}~
\cite{Ivanovici-Lascar-Lebeau-Planchon-2023} yield an enlarged range
of admissible exponents. More precisely, by
\cite[Theorem 1.5]{Ivanovici-Lascar-Lebeau-Planchon-2023}, the scaling
relation~\eqref{admissible-pair} remains unchanged, whereas the
additional restriction~\eqref{BSS-scaling} is replaced by
\begin{equation}\label{Ivanovici-admissibility}
	\frac{1}{r}
	\leq
	\left(\frac{5}{6}-\varepsilon\right)
	\left(\frac{1}{2}-\frac{1}{s}\right),
	\qquad \varepsilon>0.
\end{equation}
Here, \(\varepsilon\) is an arbitrarily small positive parameter in
the Strichartz estimate. Combining
\eqref{admissible-pair} and
\eqref{Ivanovici-admissibility}, we obtain
\begin{equation*}
	6\leq s\leq
	\frac{2(13+6\varepsilon)}{1+6\varepsilon},
	\quad\mbox{or, equivalently,}\quad
	\frac{13+6\varepsilon}{5-6\varepsilon}
	\leq r\leq\infty,
\end{equation*}
for \(0<\varepsilon<5/6\), where \(s=6\) corresponds to
\(r=\infty\). As \(\varepsilon\to0^+\), the upper endpoint for \(s\)
approaches \(26\), while the lower endpoint for \(r\) approaches
\(13/5\).

The improvement is particularly important for the pair
\((r,s)=(3,18)\), which is needed in our subsequent analysis.
Indeed, this pair satisfies the scaling relation
\eqref{admissible-pair}, since
\begin{equation*}
	\frac{1}{3}+\frac{3}{18}=\frac{1}{2}.
\end{equation*}
However, it does not satisfy the restriction in
Proposition~\ref{Prop-Strichartz-estimate}, because
\begin{equation*}
	\frac{3}{3}+\frac{2}{18}
	=1+\frac{1}{9}>1.
\end{equation*}
In contrast, the refined estimate permits this pair under the strict
convexity assumption. Indeed, for \((r,s)=(3,18)\), condition
\eqref{Ivanovici-admissibility} becomes
\begin{equation*}
	\frac{1}{3}
	\leq
	\left(\frac{5}{6}-\varepsilon\right)
	\left(\frac{1}{2}-\frac{1}{18}\right)
	=
	\left(\frac{5}{6}-\varepsilon\right)\frac{4}{9},
\end{equation*}
which holds whenever \(0<\varepsilon\leq 1/24\).
Consequently, provided the corresponding inhomogeneous Strichartz
estimate applies to the solution under consideration, we obtain
\begin{equation}\label{L3L18}
	u\in L^3(0,T;L^{18}(\Omega)).
\end{equation}
This additional space-time regularity plays a key role in the
analysis of the long-time dynamics in the non-degenerate damping
regime.

For completeness, we recall that the refined Strichartz estimates of
Ivanovici \textit{et al.}~
\cite[Theorem 1.5]{Ivanovici-Lascar-Lebeau-Planchon-2023}
are established for the homogeneous wave equation on manifolds with
smooth, strictly convex boundary. Their extension to the inhomogeneous
problem \eqref{wave-linear} follows the same reasoning used by Blair
\textit{et al.}~\cite{Blair} for Riemannian manifolds with boundary.
More specifically, the passage from homogeneous to inhomogeneous
Strichartz estimates is standard; see, for instance,
\cite[Section 3, Theorem 3.2]{Smith-Sogge}, where the
Christ--Kiselev lemma~\cite{Christ-Kiselev} is employed.

\begin{proposition}[Strichartz estimate for strictly convex domains]
	\label{Thm-Ivanovici-Strichartz}
	Assume that $\Omega\subset\mathbb{R}^3$ has a smooth, strictly
	convex boundary. Let $\mathbf{u}_0\in\mathscr{H}$ and
	$G\in L^1(0,T;L^2(\Omega))$, and let
	$\mathbf{u}(t)=(u(t),u_t(t))\in C([0,T];\mathscr{H})$
	be the energy solution of \eqref{wave-linear}.
	Then there exists $T>0$ such that, for every $\varepsilon>0$
	and every admissible pair $(r,s)$ satisfying
	\eqref{Ivanovici-admissibility} and the scaling relation
	\eqref{admissible-pair}, there exists a constant
	$C_{T,\varepsilon}>0$ such that
	\begin{align*}
		\|u\|_{L^r(0,T;L^s(\Omega))}
		\leq C_{T,\varepsilon}
		\left(
		\|\mathbf{u}_0\|_{\mathscr{H}}
		+\|G\|_{L^1(0,T;L^2(\Omega))}
		\right).
	\end{align*}
	The constant $C_{T,\varepsilon}$ depends on $T$, $\varepsilon$,
	and $\Omega$, but is independent of $\mathbf{u}_0$ and $G$.
\end{proposition}

The distinction between the two geometric settings is reflected in the corresponding admissible ranges of Strichartz exponents, summarized in Table~\ref{tab:admissible-pairs}.

\begin{table}[htbp]
	\centering
	\caption{Admissible ranges of the Strichartz pairs $(r,s)$.}
	\label{tab:admissible-pairs}
	\renewcommand{\arraystretch}{1.4}
	\begin{tabular}{lcc}
		\hline
		\textbf{Boundary geometry}
		&
		\textbf{Range of $s$}
		&
		\textbf{Corresponding range of $r$}
		\\
		\hline
		Smooth
		&
		$\displaystyle 6\leq s\leq14$
		&
		$\displaystyle \frac72\leq r\leq\infty$
		\\[8pt]
		Smooth and strictly convex
		&
		$\displaystyle
		6\leq s\leq
		\frac{2(13+6\varepsilon)}{1+6\varepsilon}$
		&
		$\displaystyle
		\frac{13+6\varepsilon}{5-6\varepsilon}
		\leq r\leq\infty$
		\\
		\hline
	\end{tabular}
\end{table}

Strict convexity leaves the scaling relation~\eqref{admissible-pair} unchanged while enlarging the range of admissible Strichartz pairs. In particular, choosing \(\varepsilon>0\) sufficiently small, the pair \((r,s)=(3,18)\) satisfies the refined admissibility condition and is therefore available in the strictly convex case, whereas it lies outside the admissible range for a general smooth boundary. This yields estimate~\eqref{L3L18}, which plays a key role in the subsequent analysis of the long-time dynamics in the non-degenerate damping regime.\\

We also highlight the importance of the following lemma, which provides a simple yet powerful mechanism to control quantities satisfying superlinear inequalities. In particular, this result guarantees that small perturbations remain uniformly bounded, preventing blow-up in iterative or approximation schemes.
\begin{lemma}\cite[Lemma 3.1]{Savostianov}\label{lemma-Strichartz}
	Let \( 0 < C_0 < \infty \) and suppose that the nonnegative function \( y \in C([a,b)) \) satisfies
	\( y(a) = 0 \) and
	\begin{align*}
		y(s) \le C_0 y(s)^{\sigma} + \varepsilon,
	\end{align*}
	for some \( \sigma > 1 \) and
	$
	0 < \varepsilon < \frac{1}{2} \left( \frac{1}{2C_0} \right)^{1/(\sigma - 1)}.
	$
	Then
	\begin{align*}
		y(s) \le 2\varepsilon, \quad s \in [a,b).
	\end{align*}
\end{lemma}

Before stating the precise notion of solution, we recall that the nonlinear terms and the damping require a weak formulation that preserves the energy structure of the problem. The following definition introduces weak (energy) solutions and, in particular, the class of Shatah--Struwe solutions, which enjoy additional space-time regularity suitable for controlling the subquintic nonlinearity.

\begin{definition}[Weak (energy) solution]
	For any $T>0$, a function $\mathbf{u}(t)=(u(t),u_t(t))\in C_w([0,T];\mathscr{H})$ satisfying $\mathbf{u}(0)=(u_0,u_1)=\mathbf{u}_0$ is said to be a {\bf weak (energy) solution} to problem \eqref{P} on the interval $[0,T]$ if $\mathbf{u}\in L^{\infty}(0,T;\mathscr{H})$ and equation~$\eqref{P}$ holds in the sense of distributions, namely,
	\begin{align*}
		\int_0^T\int_{\Omega}\Big[\,-u_t\psi_t+\nabla u\cdot\nabla \psi
		+g\left(\|\mathbf{u}\|_{\mathscr{H}}^{2}\right)u_t\psi+f(u)\psi \Big]dx\,dt
		=-\left.\int_{\Omega}u_t\psi dx\right|_0^T+\int_0^T\int_{\Omega}h\psi dx\,dt,
	\end{align*}
	for all $\psi \in C([0,T];H^1_0(\Omega))\cap C^1([0,T];L^2(\Omega))$.
\end{definition}
Here, $C_w([0,T];Y)$ denotes the space of weakly continuous functions with values in a Banach space $Y$.

\begin{definition}[Shatah--Struwe solution]
	A weak solution is called a {\bf Shatah--Struwe solution} of problem \eqref{P} if the following additional regularity holds:
	\[
	u\in L^4(0,T;L^{12}(\Omega)).
	\]
\end{definition}
\begin{remark}\label{5-10} If $\mathbf{u}(t)=(u(t),u_t(t))$ is a Shatah--Struwe solution of problem \eqref{P} on the interval $[0,T]$, then $u\in L^{5}(0, T; L^{10}(\Omega))$. Indeed, by interpolation, we have
	$$\|u\|_{L^5(0, T; L^{10}(\Omega))}\leq C \|u\|^{\frac45}_{L^4(0, T; L^{12}(\Omega))}\|u\|^{\frac15}_{L^\infty(0, T; H^1_0(\Omega))}.$$
\end{remark}

\subsection{Local solutions}
We are now in a position to establish the existence of local solutions
to problem~\eqref{P}.
\begin{theorem}[Local existence]
	\label{theo-existence}
	Assume that Assumption~\ref{Assumption} holds. If the initial data
	$\mathbf{u}_0=(u_0,u_1)\in\mathscr{H}$, then there exists
	$T=T(\|\mathbf{u}_0\|_{\mathscr{H}})>0$ such that problem~\eqref{P}
	possesses a Shatah--Strichartz solution $\mathbf{u}(t)$ on $[0,T]$.
	Moreover,
	\[
	\mathbf{u}\in C([0,T];\mathscr{H}),\qquad
	u\in L^{r}(0,T;L^{s}(\Omega)),\qquad
	u_t\in L^{2(q+1)}(0,T;L^2(\Omega)),
	\]
	and
	\begin{equation}
		\label{r-s}
		\sup_{t\in[0,T]}\|\mathbf{u}(t)\|_{\mathscr{H}}
		+\|u\|_{L^{r}(0,T;L^{s}(\Omega))}
		\leq
		\mathcal{Q}(\|\mathbf{u}_0\|_{\mathscr{H}}),
	\end{equation}
	for every admissible pair $(r,s)$ for which the corresponding
	Strichartz estimate is available. Here, $\mathcal{Q}$ is a positive,
	monotone function depending only on $\|\mathbf{u}_0\|_{\mathscr{H}}$.
\end{theorem}
\begin{remark}
	The admissible pairs in Theorem~\ref{theo-existence} depend on the
	geometry of $\Omega$. For a general smooth boundary, they are determined
	by
	\[
	\frac{3}{r}+\frac{2}{s}\leq 1,
	\]
	whereas, for a smooth strictly convex boundary, the refined estimates
	of Ivanovici {\it et al.}~\cite{Ivanovici-Lascar-Lebeau-Planchon-2023}
	yield the less restrictive condition
	\[
	\frac{1}{r}
	\leq
	\left(1-\chi(3)\right)
	\left(\frac{1}{2}-\frac{1}{s}\right),
	\qquad
	\chi(3)=\frac{1}{6}+\varepsilon.
	\]
	In particular, the pair $(r,s)=(3,18)$ is admissible in the strictly
	convex case for $0<\varepsilon\leq 1/24$; see the discussion above for
	the precise ranges of admissible pairs.
\end{remark}

\paragraph{\bf Proof of Theorem \ref{theo-existence}}
The proof is based on the Galerkin method combined with Strichartz estimates. Throughout this proof, we will repeatedly use, without explicit mention, the Young, H\"older and Poincar\'e inequalities, as well as  standard Sobolev embeddings, e.g., $H^1_0(\Omega) \hookrightarrow L^6(\Omega)$. The symbol $C$ denotes generic positive constants.
\subsubsection{\bf The approximate problem}

Assume that the boundary $\Gamma$ of $\Omega$ is sufficiently regular to ensure that the eigenfunctions
$\{\omega_j\}_{j\in \mathbb{N}} \subset H^2(\Omega)\cap H^1_0(\Omega)$.
For every $N\in \mathbb{N}$, we define the finite-dimensional subspace
$$
V_N:=\mathrm{span}\{\omega_1,\cdots,\omega_N\}\subset L^2(\Omega),
$$
and denote by $\mathcal{P}_N:L^2(\Omega)\rightarrow V_N$ the orthogonal projection onto $V_N$. We approximate the initial datum $\mathbf{u}_0=(u_0,u_1)$ by a sequence $(u_{0N},u_{1N})$, where
\begin{align}
	\label{convergence-1}
	\mathbf{u}_{0N}=(u_{0N},u_{1N})=(\mathcal{P}_N u_0, \mathcal{P}_N u_1) \rightarrow (u_0,u_1)=\mathbf{u}_0
	\quad \text{strongly in} \quad \mathscr{H}.
\end{align}
For every $N\in \mathbb{N}$, we look for $t_N>0$ and
$$
u^N: [0, t_N] \to V_N
$$
satisfying the Galerkin approximate problem
\begin{eqnarray}\left\{
	\begin{array}{l}
		\label{approximate}
		\displaystyle{
			u^N_{tt} - \Delta u^N + g\left(\|\mathbf{u}^N\|^{2}_{\mathscr{H}}\right)u^N_t
			+\mathcal{P}_N f(u^N)=\mathcal{P}_Nh,}\\
		\displaystyle{
			\mathbf{u}^N(0)= (u^N(0),u^N_t(0))=\mathbf{u}_{0N}.}
	\end{array}
	\right.
\end{eqnarray}

\begin{lemma}
	For every $N\in \mathbb{N}$, there exist $t_N>0$ and a function $u^N$ solving problem \eqref{approximate}, where $u^N$ is of the form
	$$
	u^N(t)=\sum_{j=1}^N d_{jN}(t)\omega_j\in V_N, \qquad d_{jN}\in C^2([0, t_N]).
	$$
\end{lemma}

The proof follows directly from classical ODE theory and is therefore omitted.

\subsubsection{\bf A priori estimates}
\smallskip
\paragraph{\it{The first estimate}}
Multiplying \eqref{approximate} by $u_t^N(t)$ and integrating over $\Omega$, we obtain
\begin{align}\label{approx_1}
	\frac{d}{dt}E(\mathbf{u}^N(t))
	+ g\!\left(\|\mathbf{u}^N(t)\|_{\mathscr{H}}^{2}\right)\|u_t^N(t)\|^{2}
	= 0,
\end{align}
where  the energy functional $E(\cdot)$ is given by \eqref{functional-energy}. Integrating \eqref{approx_1} over $(0,t)$ for $t\le t_N$, and using \eqref{hyp-g}, we deduce
\begin{align}\label{aprox_2}
	E(\mathbf{u}^N(t))
	+\gamma\int_0^t\|u^N_t(s)\|^{2(q+1)}\,ds
	\le E(\mathbf{u}^N(0)).
\end{align}
From \eqref{hyp-f-2} and Poincaré's inequality, we have
\begin{align}\label{existA}
	\int_{\Omega}F(u^N(t))\,dx
	\ge -\frac{\theta}{2}\|u^N(t)\|^{2}-c_f|\Omega|
	\ge -\frac{\theta}{2\lambda_1}\|\nabla u^N(t)\|^{2}-c_f|\Omega|.
\end{align}
A fundamental calculation yields
\begin{align}\label{existAA}
	-\int_{\Omega}hu^N\,dx
	\ge
	-\frac{1}{\lambda_1\omega}\|h\|^2
	-\frac{\omega}{4}\|\nabla u^N(t)\|^2.
\end{align}
Combining \eqref{existA}, \eqref{existAA}, and \eqref{functional-energy}, we obtain the lower bound
\begin{align}\label{equiv_E}
	E(\mathbf{u}^N(t))
	\ge
	\frac{1}{2}\|u_t^N(t)\|^2
	+\frac{\omega}{4}\|\nabla u^N(t)\|^2
	-c_f|\Omega|-\frac{1}{\omega\lambda_1}\|h\|^2.
\end{align}
Substituting \eqref{equiv_E} into \eqref{aprox_2} and making use of \eqref{hyp-g}, we infer
\begin{align}\label{Estimate_I}
	\frac{\omega}{4}\|\mathbf{u}^N(t)\|^2_{\mathscr{H}}
	+\gamma\int_0^t\|u^N_t(s)\|^{2(q+1)}\,ds
	\le E(\mathbf{u}^N(0))
	+c_f|\Omega|+\frac{1}{\omega\lambda_1}\|h\|^2,
\end{align}
for every $t\in [0,t_N]$ and every $N\in \mathbb{N}$.
\smallskip

Next, we estimate the term $E(\mathbf{u}^N(0))$.
We first observe that, owing to \eqref{convergence-1},
\begin{equation}
	\label{initialbound}
	\sup_{N\in \mathbb{N}}\|(u_{0N},u_{1N})\|_{\mathscr{H}}
	\le \|(u_0,u_1)\|_{\mathscr{H}}
	= \|\mathbf{u}_0\|_{\mathscr{H}}.
\end{equation}
Using \eqref{hyp_f} and the uniform boundedness \eqref{initialbound}, we obtain
\begin{align*}
	\int_{\Omega}F(u_{0N})\,dx
	\le C_f\left(1+\|u_{0N}\|^{\rho+1}_{L^{\rho+1}(\Omega)}\right)
	\le C\left(1+\|\nabla u_{0N}\|^{\rho+1}\right)
	\le C\left(1+\|\mathbf{u}_0\|_{\mathscr{H}}^{\rho+1}\right).
\end{align*}
Using again \eqref{initialbound}, we deduce from \eqref{functional-energy} that
\begin{equation}
	\begin{split}
		\label{inienergy}
		E(\mathbf{u}^N(0))
		&\le \frac{1}{2}\|\mathbf{u}_{0N}\|^2_{\mathscr{H}}
		+ C\left(1+\|\mathbf{u}_0\|_{\mathscr{H}}^{\rho+1}\right)
		+\|h\|^2
		+\frac{1}{4\lambda_1}\|\nabla u_{0N}\|^2\\
		& \le C\left(1+\|\mathbf{u}_0\|^2_{\mathscr{H}}
		+\|\mathbf{u}_0\|_{\mathscr{H}}^{\rho+1}\right).
	\end{split}
\end{equation}
Hence, returning to \eqref{Estimate_I}, we obtain
\begin{align}\label{Estimate_I0}
	\|\mathbf{u}^N(t)\|^2_{\mathscr{H}}
	+\int_0^t\|u^N_t(s)\|^{2(q+1)}\,ds
	\le \mathcal{Q}(\|\mathbf{u}_0\|_{\mathscr{H}}),\qquad t\in [0, t_N],
\end{align}
where the monotone function $\mathcal{Q}(\|\mathbf{u}_0\|_{\mathscr{H}})$ is independent of $N$. Estimate \eqref{Estimate_I0} allows us to extend the local solution $u^N$ of \eqref{approximate} to the entire interval $[0,T]$ for any given $T>0$.
\medskip

\paragraph{\it{The second estimate}}
In what follows, we employ Strichartz estimates to show that the approximate solutions $u^{N}$ are bounded in $L^{r}(0,T;L^{s}(\Omega))$ for any admissible pair $(r, s)$.  To this end, we decompose the solution $u^N=v^N+w^N$, where $v^N$ solves
\begin{eqnarray}\left\{
	\begin{aligned}
		&	v^N_{tt}-\Delta v^N=G_1(t),\label{EST01}\\
		&\mathbf{v}^N(0)=(v^N(0),v^N_t(0))=\mathbf{u}^N(0),
	\end{aligned}
	\right.\end{eqnarray}
with $G_1(t):=-g\left(||\mathbf{u}^N||_{\mathscr{H}}^{2}\right)v^N_t+\mathcal{P}_N h$,
and $w^N$ solves
\begin{eqnarray}
	\left\{
	\begin{aligned}
		&w^N_{tt}-\Delta w^N=G_2(t),\label{EST02}\\
		&\mathbf{w}^N(0)=(w^N(0),w^N_t(0))=(0,0),
	\end{aligned}\right.
\end{eqnarray}
with $G_2(t):=-g\left(||\mathbf{u}^N||_{\mathscr{H}}^{2}\right)w^N_t-\mathcal{P}_N f(v^N+w^N)$.
Arguing as above, we deduce that
$$\mathbf{v}^N,  \mathbf{w}^N\in C([0, T]; \mathscr{H})\qquad \text{for all}\quad N\in \mathbb{N},\, T>0.$$
\paragraph{\underline{Strichartz estimate for Eq.~\eqref{EST01}}} Using the multiplier $v^N_t(t)$ in $\eqref{EST01}_1$ yields
\begin{eqnarray*}
	\frac{d}{dt}\left\{\,\frac{1}{2}||\mathbf{v}^N(t)||^2_{\mathscr{H}}-\int_{\Omega}(\mathcal{P}_N h)v^Ndx\,\right\}
	=-g\left(||\mathbf{u}^N(t)||_{\mathscr{H}}^{2}\right)\|v^N_t(t)\|^2\le 0.
\end{eqnarray*}
Using the estimate
$$
\left|-\int_{\Omega}(\mathcal{P}_N h)v^Ndx\right|
\le
\frac{1}{\lambda_1}\|h\|^2+\frac{1}{4}\|\nabla v^N(t)\|^2,
$$
and combining it with $\eqref{EST01}_2$ and \eqref{initialbound}, we obtain
\begin{align}
	\label{EST3'}
	\frac{1}{4}\|\mathbf{v}^N(t)\|^2_{\mathscr{H}}
	\le
	\frac{1}{2}\|\mathbf{v}^N(0)\|^2_{\mathscr{H}}
	+\frac{1}{4}\|\nabla v^N(0)\|^2
	+\frac{2}{\lambda_1}\|h\|^2
	\le
	\mathcal{Q}\bigl(\|\mathbf{u}_0\|_{\mathscr{H}}\bigr),
\end{align}
for all $t\in [0,T]$. On the other hand, thanks to \eqref{Estimate_I0} and the continuity of the function $g$, we observe that
\begin{equation}
	\label{est-g0}
	g\left(\|\mathbf{u}^N(t)\|^2_{\mathscr{H}}\right)\le \mathcal{Q}(||\mathbf{u}_0||_{\mathscr{H}}),
\end{equation}
for all $t\in [0, T]$ and $N\in \mathbb{N}$. Thus, using \eqref{initialbound}, \eqref{EST3'} and \eqref{est-g0}, we obtain
\begin{align*}
	\int_0^1 \|G_1(s)\|d s\leq \int_0^1g\left(||\mathbf{u}^N(s)||^{2}_{\mathscr{H}}\right)\|v_t^N(s)\|ds \leq \mathcal{Q}(||\mathbf{u}_0||_{\mathscr{H}}),
\end{align*}
that is, $G_1\in L^1(0, 1; L^2(\Omega))$. Recalling that $\mathbf{v}^N\in C([0,T];\mathscr{H})$ for all $T>0$, and using \eqref{initialbound} together with Proposition~\ref{Prop-Strichartz-estimate} when $\Omega$ has a smooth boundary, or Proposition~\ref{Thm-Ivanovici-Strichartz} when $\Omega$ has a smooth, strictly convex boundary, we obtain
$
v^N\in L^r(0,1;L^s(\Omega))
$
and
\begin{equation}
	\label{est-vn}
	\|v^N\|_{L^r(0,1;L^s(\Omega))}
	\leq C\left[
	\|\mathbf{u}_0^N\|_{\mathscr{H}}
	+\|G_1\|_{L^1(0,1;L^2(\Omega))}
	\right]
	\leq \mathcal{Q}(\|\mathbf{u}_0\|_{\mathscr{H}}),
\end{equation}
for all $N\in\mathbb{N}$ and for every admissible pair $(r,s)$ allowed by the Strichartz estimates corresponding to the geometry of $\Omega$.

\medskip
\paragraph{\underline{Strichartz estimate for Eq.~\eqref{EST02}}}
Next, from \eqref{Estimate_I0} and \eqref{EST3'}, it follows that
\begin{align}\label{EST55}
	||{\mathbf{w}^N}(t)||^2_{\mathscr{H}}\le 2\left\{||\mathbf{u}^N(t)||^2_{\mathscr{H}}
	+||\mathbf{v}^N(t)||^2_{\mathscr{H}}\right\}
	\le \mathcal{Q}(||\mathbf{u}_0||_{\mathscr{H}}),\quad \forall t\in [0,T].
\end{align}
Moreover, since $\rho\in[1, 5)$, estimate~\eqref{hyp_f} and the embedding $L^{10}(\Omega)\hookrightarrow  L^{2\rho}(\Omega)$ give
$$\|\mathcal{P}_Nf(v^N+w^N)\|_{L^2(\Omega)}\leq C\left[1+\|v^N\|^{\rho}_{L^{10}(\Omega)}
+\|w^N\|^{\rho}_{L^{10}(\Omega)}\right].$$
Thus, by means of the H\"older inequality and Young's inequality with $\frac{5-\rho}{5}+\frac{\rho}{5}=1$, we have
\begin{eqnarray}
	\label{EST6b}
	\begin{aligned}
		& \|\mathcal{P}_Nf(v^N+w^N)\|_{L^1(0,t;L^2(\Omega))}\\
		&\;\le C\left[t+t^{(5-\rho)/5}\|v^N\|^{\rho}_{L^{5}(0,t;L^{10}(\Omega))}
		+t^{(5-\rho)/5}\|w^N\|^{\rho}_{L^{5}(0,t;L^{10}(\Omega))}\right]\\
		&\;\le C\left[t+t^{(5-\rho)/5}\|v^N\|^{\rho}_{L^{5}(0,t;L^{10}(\Omega))}\right]
		+C\|w^N\|^{5}_{L^{5}(0,t;L^{10}(\Omega))}.
	\end{aligned}
\end{eqnarray}
Now, for any admissible pair $(r, s)$ with $(r, s)\neq (5, 10)$,   define
\[
Z_N(t):=\|w^N\|_{L^{r}(0,t;L^{s}(\Omega))}
+\|w^N\|_{L^{5}(0,t;L^{10}(\Omega))}.
\]
Recalling that $\mathbf{w}^N\in C([0, T]; \mathscr{H})$, applying the Strichartz estimate \eqref{Strichartz-estimate} to system \eqref{EST02} with the admissible pairs $(r,s)\neq (5,10)$ and $(r,s)=(5,10)$ together with estimates \eqref{est-g0}, \eqref{EST55} and \eqref{EST6b}, we obtain that
\begin{eqnarray}
	\label{EST7b}
	\begin{aligned}
		Z_N(t)&\;\le C\int_0^tg\left(||\mathbf{u}^N(s)||^{2}_{\mathscr{H}}\right)\|w^N_t(s)\|ds+C\|\mathcal{P}_Nf(v^N+w^N)\|_{L^1(0,t;L^2(\Omega))}\\
		&\;\le \mathcal{Q}(||\mathbf{u}_0||_{\mathscr{H}})t+ C\left[t+t^{(5-\rho)/5}\|v^N\|^{\rho}_{L^{5}(0,t;L^{10}(\Omega))}\right]+C\|w^N\|^{5}_{L^{5}(0,t;L^{10}(\Omega))}.\\
		&\;\le \mathcal{Q}(||\mathbf{u}_0||_{\mathscr{H}})\left[t+ t^{(5-\rho)/5}\|v^N\|^{\rho}_{L^{5}(0,t;L^{10}(\Omega))}\right]+CZ_N(t)^5.
	\end{aligned}
\end{eqnarray}
We now choose $t>0$ sufficiently small so that the first term on the right-hand side is small. More precisely, taking $\varepsilon=\frac{1}{4}\left(\frac{1}{2C}\right)^{\frac{1}{4}}$, it follows from \eqref{est-vn} that there exists  $T = T\bigl(\|\mathbf{u}_0\|_{\mathscr{H}}\bigr)\leq 1$, which is nonincreasing with respect to $\|\mathbf{u}_0\|_{\mathscr{H}}$, such that
\[
\mathcal{Q}(||\mathbf{u}_0||_{\mathscr{H}})\left[t + t^{(5-\rho)/5}
\|v^N\|_{L^{5}(0,1;L^{10}(\Omega))}^{\rho}\right]
\le \varepsilon,
\qquad 0<t\le T.
\]
Substituting this estimate into \eqref{EST7b}, we obtain
\begin{align*}
	Z_N(t)\le \varepsilon+CZ_N(t)^5,\quad Z_N(0)=0,\quad\quad t\le T,
\end{align*}
where the constant $C$ is independent of $N$, $t$, and $\varepsilon$.
By Lemma~\ref{lemma-Strichartz},
\begin{align}\label{EST8}
	Z_N(t) \le \frac{1}{2}\left(\frac{1}{2C}\right)^{\frac{1}{4}},
	\qquad \forall\, t\in[0,T).
\end{align}
Therefore, noting that $\mathbf{w}^N\in C([0, T]; \mathscr{H})$, estimates \eqref{est-vn} and \eqref{EST8} yield the desired uniform estimate
\begin{equation}\label{EST9}
	\|u^{N}\|_{L^{r}(0,t;L^{s}(\Omega))}
	+ \|u^{N}\|_{L^{5}(0,t;L^{10}(\Omega))}
	\leq \mathcal{Q}(\|\mathbf{u}_0\|_{\mathscr{H}}),
	\qquad \forall, t\in[0,T].
\end{equation}
Here, $(r,s)$ is any admissible pair provided by the corresponding Strichartz estimate for the geometry of $\Omega$, and the positive, monotone function $\mathcal{Q}$ is independent of $N$, $t$, and the particular approximation $u^N$.

In what follows, we will make use of the estimates \eqref{EST9} in the cases $L^4(0,T;L^{12}(\Omega))$ and $L^5(0,T;L^{10}(\Omega))$ to prove the local existence of Shatah--Struwe solutions, as well as their uniqueness and continuous dependence on the initial data.
\medskip

\paragraph{\it The third estimate.}
Setting $u^{M,N}=u^M-u^N$, then $u^{M,N}$ satisfies
\begin{equation}\left\{\label{third-1}
	\begin{aligned}
		&u^{M,N}_{tt}-\Delta u^{M,N}+g\left(||\mathbf{u}^M||_{\mathscr{H}}^{2}\right)u^{M,N}_t +\left[\,g\left(||\mathbf{u}^M||_{\mathscr{H}}^{2}\right)-g\left(||\mathbf{u}^N||_{\mathscr{H}}^{2}\right)\,\right]u_t^{N}\\
		&\;+\left[\,\mathcal{P}_Mf(u^{M})-\mathcal{P}_Nf(u^{N})\,\right]
		=\mathcal{P}_M h-\mathcal{P}_N h.\\
		&\;{\mathbf{u}^{M,N}}(0)=\mathbf{u}^M(0)-\mathbf{u}^N(0).
	\end{aligned}\right.
\end{equation}
Multiplying \eqref{third-1} by $u^{M,N}_t$ and integrating over $\Omega$, we get
\begin{align}\label{eq3.18}
	\frac{1}{2}	\frac{d}{dt}||\mathbf{u}^M(t)-\mathbf{u}^N(t)||^2_{\mathscr{H}}
	+g\left(||\mathbf{u}^M(t)||_{\mathscr{H}}^{2}\right)\|u^{M,N}_t(t)\|^2
	=I_1(t)+I_2(t)+I_3(t),\end{align}
where
\begin{eqnarray*}	I_1&=&-\,\int_{\Omega}
	\left[\,\mathcal{P}_Mf(u^{M})-\mathcal{P}_Nf(u^{N})\,\right]u^{M,N}_tdx,\\ I_2&=&-\,\left[\,g\left(||\mathbf{u}^M(t)||_{\mathscr{H}}^{2}\right)
	-g\left(||\mathbf{u}^N(t)||_{\mathscr{H}}^{2}\right)\,\right]
	\int_{\Omega}u_t^{N}u^{M,N}_tdx,\\
	I_3&=&-\,\int_{\Omega}
	\left[\,\mathcal{P}_M h-\mathcal{P}_N h\,\right]u^{M,N}_tdx.
\end{eqnarray*}
Note that, $I_1$ can be rewritten as
\begin{align*}
	I_1&= \underbrace{-\int_{\Omega}\big(f(u^{M})-f(u^{N})\big)\, u^{M,N}_t\,dx}_{I_1'} - \underbrace{\int_{\Omega}(\mathcal{P}_M - I)f(u^{M})\, u^{M,N}_t\,dx}_{I_1''} \\
	&\quad + \underbrace{\int_{\Omega}(\mathcal{P}_N - I)f(u^{N})\, u^{M,N}_t\,dx}_{I_1'''}.
\end{align*}
From Assumption \eqref{hyp_f'}, we obtain
\begin{equation}
	\begin{split}
		\label{f-stable}
		|I_1'|
		&\leq\; C\left[\,1+\|u^{M}\|_{L^{12}(\Omega)}^{\rho-1}
		+\|u^{N}\|_{L^{12}(\Omega)}^{\rho-1}\,\right]
		\|\nabla u^{M,N}\|\|u^{M,N}_t\|\\
		&\leq C\left[\,1+\|u^{M}\|_{L^{12}(\Omega)}^{4}
		+\|u^{N}\|_{L^{12}(\Omega)}^{4}\,\right]
		||\mathbf{u}^M-\mathbf{u}^N||^2_{\mathscr{H}}.
	\end{split}
\end{equation}
The terms $I_1''$ and $I_1'''$ can be estimated as follows:
\begin{eqnarray*}
	| I_1''+I_1'''+I_3|
	\le r_{M,N}(t) \|u^{M,N}_t(t)\|
	+ \|(\mathcal{P}_N-I)f(u^N)\|\|u^{M,N}_t(t)\|
	\le r_{M,N}(t)\|\mathbf{u}^M(t)-\mathbf{u}^N(t)\|_{\mathscr{H}},
\end{eqnarray*}
where
\[
r_{M,N}(t) := \|(\mathcal{P}_M - I)f(u^M)\|
+ \|(\mathcal{P}_N - I)f(u^N)\|+\|(\mathcal{P}_M - \mathcal{P}_N)h\|.
\]
Note that, for \(q\ge \frac{1}{2}\), the function
\(r\mapsto r^{2q}\) belongs to \(C^{1}([0,\infty))\). Hence, by the
mean value theorem, we have
\begin{equation}\label{gen_sol_3}
	\begin{aligned}
		&\left|
		g\!\left(\|\mathbf u^M\|_{\mathscr{H}}^{2}\right)
		-
		g\!\left(\|\mathbf u^N\|_{\mathscr{H}}^{2}\right)
		\right|
		\\
		&\qquad
		=
		\gamma\left|
		\|\mathbf u^M\|_{\mathscr{H}}^{2q}
		-
		\|\mathbf u^N\|_{\mathscr{H}}^{2q}
		\right|
		\\
		&\qquad
		=
		2q\gamma
		\left|
		\int_{0}^{1}
		\left(
		\theta\|\mathbf u^M\|_{\mathscr{H}}
		+
		(1-\theta)\|\mathbf u^N\|_{\mathscr{H}}
		\right)^{2q-1}
		\,d\theta
		\right|
		\left|
		\|\mathbf u^M\|_{\mathscr{H}}
		-
		\|\mathbf u^N\|_{\mathscr{H}}
		\right|
		\\
		&\qquad
		\le
		2q\gamma
		\left(
		\|\mathbf u^M\|_{\mathscr{H}}
		+
		\|\mathbf u^N\|_{\mathscr{H}}
		\right)^{2q-1}
		\|\mathbf u^M-\mathbf u^N\|_{\mathscr{H}}.
	\end{aligned}
\end{equation}
Then, using again \eqref{Estimate_I0}, we have
\begin{align*}
	|I_2| &\le \mathcal{Q}(||\mathbf{u}_0||_{\mathscr{H}})
	||\mathbf{u}^M-\mathbf{u}^N||_{\mathscr{H}}
	\|u^N_t\|\|u^{M,N}_t\|\le \mathcal{Q}(||\mathbf{u}_0||_{\mathscr{H}})
	||\mathbf{u}^M-\mathbf{u}^N||^2_{\mathscr{H}}.
\end{align*}
Thus, substituting the estimates of the terms $I_1$ and $I_2$ into \eqref{eq3.18}, we obtain
\begin{align*}
	\frac{d}{dt}||\mathbf{u}^M-\mathbf{u}^N||^2_{\mathscr{H}}&\le r_{M,N}(t)\|\mathbf{u}^M-\mathbf{u}^N\|_{\mathscr{H}}\\ &+\mathcal{Q}(||\mathbf{u}_0||_{\mathscr{H}})
	\left[\,1+\|u^{M}\|_{L^{12}(\Omega)}^{4}+\|u^{N}\|_{L^{12}(\Omega)}^{4}\,\right]
	||\mathbf{u}^M-\mathbf{u}^N||^2_{\mathscr{H}}.
\end{align*}
Applying the Gronwall lemma to the above inequality and making use of  \eqref{EST9} with $(r, s)=(4, 12)$, we find
\begin{equation*}
	\begin{split}
		\max_{t\in[0, T]}\|\mathbf{u}^M(t)-\mathbf{u}^N(t)\|^2_{\mathscr{H}}
		& \leq \mathcal{Q}(||\mathbf{u}_0||_{\mathscr{H}})
		||\mathbf{u}^M(0)-\mathbf{u}^N(0)||^2_{\mathscr{H}}\\
		& + \mathcal{Q}(||\mathbf{u}_0||_{\mathscr{H}})
		\int_0^T r_{M,N}(t)\|\mathbf{u}^M(t)-\mathbf{u}^N(t)\|_{\mathscr{H}}\,dt.
	\end{split}
\end{equation*}
Since
\begin{equation*}
	\begin{split}
		\mathcal{Q}(||&\mathbf{u}_0||_{\mathscr{H}})
		\int_0^T r_{M,N}(t)\|\mathbf{u}^M(t)-\mathbf{u}^N(t)\|_{\mathscr{H}}\,dt\\
		\leq\, & \mathcal{Q}(||\mathbf{u}_0||_{\mathscr{H}}) \max_{t\in[0, T]}\|\mathbf{u}^M(t)-\mathbf{u}^N(t)\|_{\mathscr{H}}\int_0^T r_{M,N}(t)\,dt\\
		\leq\, & \frac12 \max_{t\in[0, T]}\|\mathbf{u}^M(t)-\mathbf{u}^N(t)\|^2_{\mathscr{H}} + \mathcal{Q}(||\mathbf{u}_0||_{\mathscr{H}}) \left(\int_0^T r_{M,N}(t)\,dt\right)^2,
	\end{split}
\end{equation*}
it follows that
\begin{eqnarray}
	\begin{aligned}\label{lipschtz-ineq}
		&\max_{t\in[0, T]}\|\mathbf{u}^M(t)-\mathbf{u}^N(t)\|^2_{\mathscr{H}}\\
		&\le \mathcal{Q}(||\mathbf{u}_0||_{\mathscr{H}})
		||\mathbf{u}^M(0)-\mathbf{u}^N(0)||^2_{\mathscr{H}}+ \mathcal{Q}(||\mathbf{u}_0||_{\mathscr{H}}) \left(\int_0^T r_{M,N}(t)\,dt\right)^2.
	\end{aligned}
\end{eqnarray}
On the other hand, by \eqref{hyp_f} and \eqref{EST9}, we deduce that
\begin{equation*}
	\|f(u^N)\|_{L^1(0, T; L^2(\Omega))}\leq C(1+ \|u^N\|_{L^5(0, T; L^{10}(\Omega))})\leq \mathcal{Q}(||\mathbf{u}_0||_{\mathscr{H}}), \qquad \forall N\in \mathbb{N},
\end{equation*}
that is, $\{f(u^N)\}$ is uniformly bounded in $L^1(0,T;L^2(\Omega))$. Since $\mathcal{P}_N \to I$ strongly in $L^2(\Omega)$, it follows that
\[
\int_0^T \|(\mathcal{P}_N-I)f(u^N(s))\| ds \longrightarrow 0
\quad \text{as } N\to +\infty.
\]
Thus, recalling that $h\in L^2(\Omega)$, we have
\begin{align*}
	\int_0^t r_{M,N}(s)\,ds
	& =\int_0^t \left[\|(\mathcal{P}_M-I)f(u^M(s))\|
	+ \|(\mathcal{P}_N-I)f(u^N(s))\|+\|(\mathcal{P}_M - \mathcal{P}_N)h\| \right]ds \\ & \longrightarrow_{M,N\to +\infty} 0.
\end{align*}
Therefore,  it follows from \eqref{convergence-1} and \eqref{lipschtz-ineq} that
\begin{align*}
	\max_{t\in [0, T]}||\mathbf{u}^M(t)-\mathbf{u}^N(t)||^2_{\mathscr{H}}\longrightarrow 0
	\quad \text{as } M, N\to +\infty.
\end{align*}
This proves that $\{\mathbf{u}^{N}\}_{N\in \mathbb{N}}$ is a Cauchy sequence in $C([0, T]; \mathscr{H})$. Therefore, there exists a function $\mathbf{u}\in C([0, T]; \mathscr{H})$ such that $\mathbf{u}^N\rightarrow \mathbf{u}$ in $ C([0, T]; \mathscr{H})$, that is
\begin{equation}\label{a1}
	\lim_{N\rightarrow +\infty}\max_{t\in [0,T]}||\mathbf{u}^N(t)-\mathbf{u}(t)||^2_{\mathscr{H}}=0.
\end{equation}
Moreover, the combination of \eqref{EST9} and \eqref{a1} implies that, up to a subsequence, $u^N\rightarrow u$ weakly in $ L^r(0, T; L^{s}(\Omega))$ and
\begin{equation}
	\label{4-12}
	\|u\|_{L^r(0, T; L^{s}(\Omega))}
	\leq \liminf_{N\rightarrow +\infty} \|u^N\|_{L^r(0, T; L^{s}(\Omega))}
	\leq \mathcal{Q}(||\mathbf{u}_0||_{\mathscr{H}}).
\end{equation}
In addition,  \eqref{a1} implies that $\mathbf{u}^N(0)\rightarrow \mathbf{u}(0)$ in $\mathscr{H}$. Combining this with \eqref{convergence-1} and $\eqref{approximate}_2$ implies that
\begin{equation*}
	\mathbf{u}(0)=(u(0),u_t(0))=(u_0, u_1)=\mathbf{u}_0.
\end{equation*}

\subsubsection{\bf Limiting process.} From \eqref{approximate}, we obtain
\begin{eqnarray}\label{a0}
	\begin{aligned}
		&\int_0^T\left[-\left(u^N_t,\psi_t\right)+(\nabla u^N,\nabla\psi)
		+g\left(||\mathbf{u}^N||^{2}_{\mathscr{H}}\right)(u^N_t,\psi)
		+(\mathcal{P}_N(f(u^N)),\psi)\right]dt\\ &\quad=-\left.\left(u^N_t,\psi\right)\right|_0^T
		+\int_0^T(\mathcal{P}_Nh,\psi)dt,
	\end{aligned}
\end{eqnarray}
holds for all $N\in \mathbb N$ and every $\psi\in C([0,T];H^1_0(\Omega))\cap C^1([0,T];L^2(\Omega))$. We pass to the limit in  \eqref{a0}. From estimate \eqref{a1}, it follows that
\begin{align*}
	\int_0^T(u^N_t,\psi_t)dt&\rightarrow  \int_0^T(u_t,\psi_t)dt\ \ \hbox{as}\ \ N\rightarrow +\infty,\\
	\int_0^T(\nabla u^N,\nabla \psi)dt&\rightarrow\int_0^T(\nabla u,\nabla \psi)dt\ \ \hbox{as}\ \ N\rightarrow +\infty,\\
	\left.\left(u^N_t,\psi\right)\right|_0^T&\rightarrow  \left.\left(u_t,\psi\right)\right|_0^T\ \ \hbox{as}\ \ N\rightarrow +\infty.
\end{align*}
Moreover, recalling that $T\leq 1$ and arguing as in the proof of  \eqref{gen_sol_3}, we obtain from \eqref{Estimate_I0}, \eqref{est-g0}, and \eqref{a1} that
\begin{eqnarray*}
	&&\left|\int_0^Tg\left(||\mathbf{u}^N(t)||^{2}_{\mathscr{H}}\right)
	(u^N_t,\psi)dt-\int_0^Tg\left(||\mathbf{u}(t)||^{2}_{\mathscr{H}}\right)
	(u_t,\psi)dt \right|\\
	&&\quad \le  \left|\int_0^Tg\left(||\mathbf{u}^N(t)||^{2}_{\mathscr{H}}\right)
	((u^N_t-u_t),\psi)dt \right|+ \left|\int_0^T\left[\,g\left(||\mathbf{u}^N(t)||^{2}_{\mathscr{H}}\right)
	-g\left(||\mathbf{u}(t)||^{2}_{\mathscr{H}}\right)\,\right]
	(u_t,\psi)dt \right|\\
	&&\quad \leq \mathcal{Q}(||\mathbf{u}_0||_{\mathscr{H}})\|\psi\|_{C([0, T]; L^2(\Omega))}\max_{t\in [0, T]}
	||\mathbf{u}^N(t)-\mathbf{u}(t)||_{\mathscr{H}}\rightarrow 0\quad \mbox{as}\quad N\rightarrow +\infty.
\end{eqnarray*}
Finally, arguing as in \eqref{f-stable}, and using \eqref{EST9}, \eqref{a1}, and \eqref{4-12}, we obtain
\begin{eqnarray*}
	&& \left|	\int_0^T\int_{\Omega}f(u^N)\psi dxdt-\int_0^T\int_{\Omega}f(u)\psi dxdt\right|\\
	&&\quad\le C\|\psi\|_{C([0, T]; L^2(\Omega))}\left[\,T+\|u^N\|^4_{L^4(0,T;L^{12}(\Omega))}
	+\|u\|^4_{L^4(0,T:L^{12}(\Omega))}\,\right]\max_{t\in [0,T]}\|\nabla u^N(t)-\nabla u(t)\|\\
	&&\quad \le \mathcal{Q}(||\mathbf{u}_0||_{\mathscr{H}})\|\psi\|_{C([0, T]; L^2(\Omega))}\max_{t\in [0, T]}
	||\mathbf{u}^N(t)-\mathbf{u}(t)||_{\mathscr{H}}\rightarrow 0\quad \mbox{as}\quad N\rightarrow +\infty.
\end{eqnarray*}
Therefore, letting $N\to +\infty$ in \eqref{a0} and using the above convergences, we conclude that
$\mathbf{u}=(u, u_t)$ is a Shatah--Struwe solution of problem \eqref{P} on the interval $[0,T]$.
\qed

\subsection{Energy identity}
\begin{theorem}\label{theo-energy-equality}
	Assume that the hypotheses of Theorem~\ref{theo-existence} are satisfied, and let
	$
	\mathbf{u}(t)=(u(t),u_t(t))
	$
	be a  Shatah--Struwe solution of problem~\eqref{P} on
	$[0,T]$. Then $\mathbf u$ satisfies the following energy identity:
	\begin{equation}\label{energy-equality}
		E(\mathbf{u}(t))
		+ \int_{\tau}^{t}
		g\left(\|\mathbf{u}(s)\|_{\mathscr{H}}^{2}\right)
		\|u_t(s)\|^{2}\,ds
		= E(\mathbf{u}(\tau)),
	\end{equation}
	for all $0\le \tau\le t\le T$, where $E(\cdot)$ is the energy functional defined in~\eqref{functional-energy}.
\end{theorem}
\begin{proof} Let $u^N(t)=\mathcal{P}_N u(t)$. Then $u^N$ satisfies 
	$$
	u^N_{tt}-\Delta u^N
	+g\left(\|\mathbf{u}\|_{\mathscr{H}}^{2}\right)u^N_t
	+ \mathcal{P}_N f(u)=\mathcal{P}_Nh.
	$$
	Multiplying the above equation by $u_t^N$ and integrating over
	$\Omega \times [\tau, t]$, we obtain
	\begin{eqnarray}\label{PP}
		\begin{aligned}
			& \frac12\|\mathbf{u}^N(t)\|_{\mathscr{H}}^2-\int_{\Omega} h\,u^N(t)\,dx+ \int_\tau^t g\left(
			\|\mathbf{u}(s)\|_{\mathscr{H}}^{2}\right)
			\|u^N_t(s)\|^2\,ds \\
			&\quad = -\int_\tau^t\!\!\int_{\Omega}
			\mathcal{P}_N f(u(s))\,u^N_t(s)\,dx\,ds
			+ \frac12\|\mathbf{u}^N(\tau)\|_{\mathscr{H}}^2-\int_{\Omega}\mathcal{P}_N h\,u^N(\tau)\,dx.
		\end{aligned}
	\end{eqnarray}
	Since $u_t\in L^2(0, T; L^2(\Omega))$, we claim that
	\begin{equation}\label{con-ut}
		u_t^N\to u_t \quad \text{in } L^2(\tau, t; L^2(\Omega)).
	\end{equation}
	Indeed, since $u_t^N(s)\to u_t(s)$ in $L^2(\Omega)$ for all $s\in [\tau, t]$ and
	$\|u_t^N(s)-u_t(s)\|^2 \le 4\|u_t(s)\|^2$,
	the dominated convergence theorem implies \eqref{con-ut}. On the other hand, since $\mathbf{u}\in C([0, T]; \mathscr{H})$ and $g\in C([0, \infty))$, we have
	\[
	\sup_{s\in [0, T]}g\left(\|\mathbf{u}(s)\|_{\mathscr{H}}^{2}\right)\leq C_0,
	\]
	where $C_0$ depends on $u$ and $T$. Hence, by \eqref{con-ut},
	\begin{eqnarray*}
		\begin{aligned}
			&\left|
			\int_\tau^t g\left(\|\mathbf{u}(s)\|_{\mathscr{H}}^{2}\right)
			\|u_t^N(s)\|^2\,ds
			- \int_\tau^t g\left(\|\mathbf{u}(s)\|_{\mathscr{H}}^{2}\right)
			\|u_t(s)\|^2\,ds
			\right|      \\
			&\quad \le C_0 \|u_t^N- u_t\|_{L^2(\tau, t; L^2(\Omega))}
			\to 0 \quad \text{as } N\to +\infty.
		\end{aligned}
	\end{eqnarray*}
	That is,
	\begin{align}\label{energy-eq-1}
		\int_\tau^t g\left(\|\mathbf{u}(s)\|_{\mathscr{H}}^{2}\right)
		\|u_t^N(s)\|^2\,ds
		\to
		\int_\tau^t g\left(\|\mathbf{u}(s)\|_{\mathscr{H}}^{2}\right)
		\|u_t(s)\|^2\,ds.
	\end{align}
	Moreover, using estimate~\eqref{hyp_f} and Remark~\ref{5-10}, we have
	\[
	\|\mathcal{P}_N f(u)\|_{L^1(\tau, t; L^2(\Omega))}
	\leq \|f(u)\|_{L^1(0, T; L^2(\Omega))}
	\leq C(1+ \|u\|_{L^5(0, T; L^{10}(\Omega))})
	\leq C_0.
	\]
	Noting that $\mathcal{P}_N f(u(s))\to f(u(s))$ in $L^2(\Omega)$ and
	\[
	\|\mathcal{P}_N f(u(s))-f(u(s))\|\leq 2\|f(u(s))\|,
	\]
	for all $s\in [\tau, t]$, we obtain
	\[
	\mathcal{P}_N f(u)\to f(u)\quad \text{in } L^1(\tau, t; L^2(\Omega)).
	\]
	This, together with $u_t\in C([0, T]; L^2(\Omega))$ and \eqref{con-ut}, yields
	\begin{align}\label{energy-eq2}
		\int_\tau^t \int_\Omega \mathcal{P}_N f(u)\, u_t^N\,dx\,ds
		\to
		\int_\tau^t \int_\Omega f(u)\,u_t\,dx\,ds.
	\end{align}
	Passing to the limit as $N\to\infty$ in \eqref{PP}, and using \eqref{energy-eq-1} and \eqref{energy-eq2}, we obtain
	\begin{eqnarray}\label{PPP}
		\begin{aligned}
			& \frac12\|\mathbf{u}(t)\|_{\mathscr{H}}^2
			- \int_{\Omega} h\,u(t)\,dx
			+ \int_\tau^t
			g\left(\|\mathbf{u}(s)\|_{\mathscr{H}}^{2}\right)
			\|u_t(s)\|^2\,ds \\
			&\quad=
			\frac12\|\mathbf{u}(\tau)\|_{\mathscr{H}}^2
			- \int_{\Omega} h\,u(\tau)\,dx
			- \int_\tau^t\int_{\Omega} f(u)\,u_t\,dx\,ds .
		\end{aligned}
	\end{eqnarray}
	Finally, since $f(u)u_t\in L^1([\tau, t]\times \Omega)$, we have
	\[
	\int_\tau^t \int_\Omega f(u)\,u_t\,dx\,ds
	= \int_\Omega F(u(t))\,dx
	- \int_\Omega F(u(\tau))\,dx.
	\]
	Substituting this into \eqref{PPP}, we obtain the desired energy identity \eqref{energy-equality}.
\end{proof}

\begin{remark}\label{unibound}
	Arguing as in \eqref{equiv_E} and \eqref{inienergy} for the Galerkin approximate solutions $u^N$, one can infer from the energy identity \eqref{energy-equality} that the Shatah-Struwe solution  ${\mathbf{u}}=(u,u_t)$  satisfies
	\begin{align}\label{Estimate_I0.1}
		\|\mathbf u(t)\|_{\mathscr{H}}^{2}
		+
		\int_0^t
		g\!\left(\|\mathbf u(s)\|_{\mathscr{H}}^{2}\right)
		\|u_t(s)\|^{2}\,ds
		\le
		\mathcal Q\!\left(\|\mathbf u(0)\|_{\mathscr{H}}\right),
		\qquad t\in[0,T].
	\end{align}
\end{remark}

\subsection{Continuous dependence and uniqueness}
\begin{theorem}[Continuous dependence]\label{theo-cont-depend}
	Assume that the hypotheses of Theorem~\ref{theo-existence} are satisfied.
	Let
	$
	\mathbf u^i=(u^i,u_t^i),$ $i=1,2,
	$
	be two Shatah--Struwe solutions of problem~\eqref{P} on
	$[0,T]$ such that
	$
	\|\mathbf u^i(0)\|_{\mathscr{H}}\le R,
	$ $i=1,2.
	$
	Then
	\begin{align}\label{relation-strong-solution}
		\|\mathbf u^1(t)-\mathbf u^2(t)\|_{\mathscr{H}}^{2}
		\le
		e^{\mathcal Q(R)
			\int_0^t
			\left(
			1+\|u^1(s)\|_{L^{12}(\Omega)}^{4}
			+\|u^2(s)\|_{L^{12}(\Omega)}^{4}
			\right)\,ds}
		\,
		\|\mathbf u^1(0)-\mathbf u^2(0)\|_{\mathscr{H}}^{2},
	\end{align}
	for all \(t\in[0,T]\). Consequently, problem~\eqref{P} admits 
	at most one Shatah--Struwe solution corresponding to a given initial
	datum.
\end{theorem}
\begin{proof}
	Let $w=u^1-u^2$. Then $w$ satisfies the following equation
	$$w_{tt}-\Delta w+g\left(||\mathbf{u}^1||_{\mathscr{H}}^{2}\right)w_t +\left[\,g\left(||\mathbf{u}^1||_{\mathscr{H}}^{2}\right)
	-g\left(||\mathbf{u}^2||_{\mathscr{H}}^{2}\right)\,\right]u_t^2
	+f(u^{1})-f(u^{2})=0.$$
	Multiplying this equation by $w_t$ and integrating over $\Omega$ (which is justified exactly as in Theorem~\ref{theo-energy-equality}),
	we proceed as in \eqref{f-stable} and \eqref{gen_sol_3}, and, in view of \eqref{r-s}, obtain
	\begin{eqnarray*}
		\frac{d}{dt}||\mathbf{u}^1-\mathbf{u}^2||^2_{\mathscr{H}}\le \mathcal{Q}(R)
		\left[\,1+\|u^{1}\|_{L^{12}(\Omega)}^{4}+\|u^{2}\|_{L^{12}(\Omega)}^{4}\,\right]
		||\mathbf{u}^1-\mathbf{u}^2||^2_{\mathscr{H}}.
	\end{eqnarray*}
	Applying the Gronwall inequality, we see that \eqref{relation-strong-solution} and uniqueness follows.
\end{proof}

\subsection{Global existence}
\begin{theorem}[Global existence]\label{theo-global-subcritical}
	Assume that the hypotheses of Theorem~\ref{theo-existence} are satisfied.
	Then, for every initial datum $\mathbf u_0\in\mathscr{H}$,
	problem~\eqref{P} admits a unique global Shatah--Struwe solution
	$\mathbf u(t)$. Moreover,
	\begin{align}\label{globalest}
		\|\mathbf u(t)\|_{\mathscr{H}}
		+
		\|u\|_{L^r(t,t+1;L^s(\Omega))}
		\le
		\mathcal Q\bigl(\|\mathbf u(0)\|_{\mathscr{H}}\bigr),
		\qquad \forall, t\ge0,
	\end{align}
	for every admissible pair $(r,s)$ covered by the corresponding Strichartz estimate for the geometry of $\Omega$, where $\mathcal Q$
	is a positive nondecreasing function depending only on
	$\|\mathbf u(0)\|_{\mathscr H}$.
\end{theorem}
\begin{proof}
	The uniqueness follows from Theorem~\ref{theo-cont-depend}, and the uniform  estimate~\eqref{Estimate_I0.1} for the energy norm $\|\mathbf{u}(t)\|_{\mathscr{H}}$ is obtained in Remark~\ref{unibound}.
	Moreover, estimate \eqref{Estimate_I0.1} ensures that $\|\mathbf{u}(t)\|_{\mathscr{H}}$ remains uniformly bounded on any finite time interval. Therefore, by Theorem~\ref{theo-existence}, the local Shatah--Struwe solution $\mathbf{u}(t)$ can be extended to a global solution and satisfies the uniform  estimate~\eqref{Estimate_I0.1} for all $t\geq 0$.
	
	To complete the proof, it remains to establish the uniform estimate for
	$\|u\|_{L^r(t,t+1;L^s(\Omega))}$. Indeed, Theorem~\ref{theo-existence} implies that, for each $t\geq 0$, there exists
	$\sigma=\sigma(\|\mathbf{u}(t)\|_{\mathscr{H}})>0$ such that
	\[
	\|u\|_{L^r(t, t+\sigma; L^s(\Omega))}
	\leq \mathcal{Q}(\|\mathbf{u}(t)\|_{\mathscr{H}}),
	\]
	for every admissible pair $(r,s)$. Recalling \eqref{Estimate_I0.1}, we have
	\[
	\|\mathbf{u}(t)\|_{\mathscr{H}}\leq \mathcal{Q}(\|\mathbf{u}(0)\|_{\mathscr{H}}),
	\quad \forall t\geq 0.
	\]
	Hence, $\sigma$ depends only on $\|\mathbf{u}(0)\|_{\mathscr{H}}$, and
	\begin{align}\label{ll}
		\|u\|_{L^r(t, t+\sigma; L^s(\Omega))}
		\leq \mathcal{Q}(\|\mathbf{u}(0)\|_{\mathscr{H}}),
	\end{align}
	for all $t\geq 0$. Let
	$N_0 := \left\lfloor \frac{1}{\sigma} \right\rfloor$.
	Then the interval $[t,t+1]$ can be covered by at most $N_0+1$ subintervals of length $\sigma$.
	Therefore, using \eqref{ll}, we obtain
	\begin{align*}
		\|u\|_{L^r(t,t+1;L^{s}(\Omega))}^r
		\leq
		\sum_{i=0}^{N_0}
		\|u\|_{L^r(t+i \sigma,\, t+(i+1)\sigma;L^{s}(\Omega))}^r \leq (N_0+1)\left[\mathcal{Q}(\|\mathbf{u}(0)\|_{\mathscr{H}}) \right]^r.
	\end{align*}
	Taking the $r$-th root, and using the fact $r>2$,  we deduce that
	\[
	\|u\|_{L^r(t,t+1;L^{s}(\Omega))}
	\le \mathcal{Q}(\|\mathbf{u}(0)\|_{\mathscr{H}})\, (N_0+1)^{1/2}.
	\]
	Since $N_0$ depends only on $\|\mathbf{u}(0)\|_{\mathscr{H}}$, the proof is complete.	
\end{proof}

\begin{remark}
	It follows directly from \eqref{Estimate_I0.1} that
	\begin{align}\label{est-u_t-a}
		\int_0^t \|u_t(s)\|^{2(q+1)} \, ds \le \mathcal{Q}(\|\mathbf{u}(0)\|_{\mathscr{H}}), \qquad t\geq 0.
	\end{align}
	In particular, if $\kappa>0$,
	\begin{align}\label{est-u-t-b}
		\int_0^t \|u_t(s)\|^{2} \, ds \le \mathcal{Q}(\|\mathbf{u}(0)\|_{\mathscr{H}}),
		\qquad  t\geq 0.
	\end{align}
\end{remark}
\section{Global Attractor}
\subsection{Statement of the main result}
Theorem \ref{theo-global-subcritical} guarantees that the  Shatah-Struwe solutions to problem \eqref{P} generate a dynamical system $(\mathscr{H},S_t)$. The evolution operator $S_t : \mathscr{H} \to \mathscr{H}$ is characterized by the relation
\begin{align*}
	S_t\mathbf{u}_0 = (u(t), u_t(t))=\mathbf{u}(t), \quad \mathbf{u}_0 \in \mathscr{H},
\end{align*}
where $\mathbf{u}$ is the unique Shatah-Struwe solution of problem \eqref{P} corresponding to the initial datum $\mathbf{u}_0 \in \mathscr{H}$. Moreover, it follows directly from Theorem \ref{theo-cont-depend} and estimate \eqref{globalest} with $(r, s)=(4, 12)$ that the dynamical system $(\mathscr{H}, S_t)$ possesses the following local \textbf{Lipschitz continuity property}
\begin{align*}
	\| S_t \mathbf{u}_0^1 - S_t \mathbf{u}_0^2 \|_{\mathscr{H}}^2 \leq e^{\mathcal{Q}(R)t} \| \mathbf{u}_0^1 - \mathbf{u}_0^2 \|_{\mathscr{H}}^2, \quad \forall t \ge 0,
\end{align*}
where $\mathbf{u}_0^i \in \mathscr{H}$ with $\|\mathbf{u}^{i}_0 \|_{\mathscr{H}} \leq R$, $i=1,2$, and the function $\mathcal{Q}(\cdot)$ is independent of $u^i$ and $t$.

Our main results on global attractors for dynamical systems $(\mathscr{H},S_t)$ are presented below.
\begin{theorem}\label{theo-main}
	Assume that the hypotheses of Theorem~\ref{theo-global-subcritical} are satisfied with $q\in[1/2,2]$. Then, the dynamical system $(\mathscr{H}, S_t)$ possesses a unique compact global attractor $\mathfrak{A} \subset \mathscr{H}$, which can be characterized as the unstable manifold
	$$
	\mathfrak{A} = \mathcal{M}^u(\mathcal{N}),
	$$
	emanating from the set $\mathcal{N}$ of stationary solutions. Moreover, $\mathfrak{A}$ consists of full trajectories
	$
	\{\mathbf{u}(t) = (u(t), u_t(t)) : t \in \mathbb{R}\}
	$
	such that
	$$
	\lim_{t \to -\infty} \mathrm{dist}_{\mathscr{H}}(\mathbf{u}(t), \mathcal{N})
	= \lim_{t \to +\infty} \mathrm{dist}_{\mathscr{H}}(\mathbf{u}(t), \mathcal{N}) = 0.
	$$
\end{theorem}
\begin{proof}
	The proof is a direct application of Theorems \ref{theo-chueshov} and
	\ref{theo-chueshov-caracterization}. More precisely, the verification that the
	dynamical system $(\mathscr{H}, S_t)$ associated with problem \eqref{P}
	satisfies the hypotheses of Theorems \ref{theo-chueshov} and
	\ref{theo-chueshov-caracterization} is carried out in Proposition
	\ref{Prop-gradient-system}, Proposition \ref{prop-absorbing-set}, and
	Corollary \ref{corollary-smootness2}. In Proposition
	\ref{Prop-gradient-system}, we show that the system $(\mathscr{H}, S_t)$
	is gradient; in Proposition \ref{prop-absorbing-set}, we prove that it is dissipative; and, in Corollary
	\ref{corollary-smootness2}, we establish that it is
	asymptotically smooth.
\end{proof}

In the remainder of this section, we focus on proving the results invoked in the proof above.
\subsection{Gradient system}
\begin{proposition}\label{Prop-gradient-system} Assume that the hypotheses of Theorem \ref{theo-global-subcritical} hold. Then the energy functional  $E:\mathscr{H}\to \mathbb{R}$ defined by \eqref{functional-energy} is a strict Lyapunov functional for the dynamical system $(\mathscr{H},S_t)$. Consequently, $(\mathscr{H},S_t)$ is a gradient dynamical system.
\end{proposition}
\begin{proof} By Theorem \ref{theo-energy-equality}, we obtain
	\begin{align}\label{inequality grad2}
		E(S_t \mathbf{u}_0)
		+ \int_0^t g\bigl(\|\mathbf{u}(\tau)\|^2_{\mathscr{H}}\bigr) \|u_t(\tau)\|^2 \, d\tau
		= E(\mathbf{u}_0)
	\end{align}
	for every $\mathbf{u}_0 \in \mathscr{H}$ and all $t\geq 0$.
	Since $f\in C^1(\mathbb{R})$, the functional $E(\cdot)$ is continuous on $\mathscr{H}$.  In view of \eqref{hyp-g}, it follows from \eqref{inequality grad2} that the mapping $t\mapsto E(S_t \mathbf{u}_0)$ is non-increasing for every $\mathbf{u}_0 \in \mathscr{H}$, and
	\[
	E(S_t \mathbf{u}_0) = E(\mathbf{u}_0) \quad \Rightarrow \quad
	\int_0^t \|u_t(\tau)\|^{2(q+1)} \, d\tau = 0,
	\]
	which implies
	\[
	\|u_t(t)\| = 0 \quad \text{for almost every } t>0.
	\]
	Moreover, since $u_t\in C([0, T];L^2(\Omega))$ for every $T>0$, we conclude  that $\|u_t(t)\|=0$ for all $t>0$. Consequently, $\mathbf{u} (t)=\mathbf{u}_0$ for all $t\geq 0$, which means that $\mathbf{u}_0\in \mathcal{N}$, where $\mathcal{N}$ denotes the  set of stationary points of the dynamical system $(\mathscr{H},S_t)$. Noting that
	$$
	\quad \mathbf{u}_0\in \mathcal{N} \quad  \Leftrightarrow \quad S_t(\mathbf{u}_0)=\mathbf{u}_0,   \quad  t>0,
	$$
	we infer that $E(\cdot)$ is a strict Lyapunov functional for the dynamical system $(\mathscr{H},S_t)$.
\end{proof}
	
	\subsection{Absorbing set}
	\begin{proposition}\label{prop-absorbing-set}
		Let $\mathbf{u}(t)=S_t\mathbf{u}_0$ be a Shatah--Struwe solution to problem~\eqref{P} with initial datum
		$\mathbf{u}_0=(u_0,u_1)\in\mathscr{H}$ satisfying $\|\mathbf{u}_0\|_{\mathscr{H}}\le R$.
		Under the assumptions of Theorem~\ref{theo-global-subcritical}, the dynamical system $(\mathscr{H},S_t)$ admits a bounded absorbing set $\mathcal{B}\subset\mathscr{H}$.
	\end{proposition}
	\begin{proof}
		To this end, we set
		\begin{align*}
			\widetilde{E}(\mathbf{u}(t)):=E(\mathbf{u}(t))+R_*,\quad \mbox{with}\quad R_* := \frac{1}{\lambda_1 \omega}\, \|h\|^{2} + c_f |\Omega|.
		\end{align*}
		Proceeding as in \eqref{existA}-\eqref{equiv_E}, we obtain
		\begin{eqnarray}
			\label{est-E-tilde}
			\widetilde{E}(\mathbf{u}(t))\ge \frac{1}{2}\|u_t(t)\|^2+\frac{\omega}{4}\|\nabla u(t)\|^2\ge \frac{\omega}{4}||\mathbf{u}(t)||^2_{\mathscr{H}},
		\end{eqnarray}
		where $\omega$ is defined in \eqref{def-w}.
		By Theorem~\ref{theo-energy-equality}, we have
		\begin{align}\label{dissip_A}
			\int_{t}^{t+1} g\left(\|\mathbf{u}(s)\|_{\mathscr{H}}^{2}\right)\|u_t(s)\|^2\,ds
			= E(\mathbf{u}(t)) - E(\mathbf{u}(t+1)) =: \mathcal{D}(t)^2.
		\end{align}
		Using \eqref{hyp-g}, it follows from \eqref{dissip_A} that
		\begin{align*}
			\int_{t}^{t+1}\|u_t(s)\|^{2(q+1)}\,ds \le \frac{1}{\gamma}\mathcal{D}(t)^2,
		\end{align*}
		which together with the H\"older inequality with $\frac{q}{q+1}+\frac{1}{q+1}=1$ gives
		\begin{align}\label{abs-01}
			\int_{t}^{t+1}\|u_t(s)\|^{2}\,ds
			&\le \left[\int_{t}^{t+1} 1\,ds\right]^{\frac{q}{q+1}}
			\left[\int_{t}^{t+1}\|u_t(s)\|^{2(q+1)}\,ds\right]^{\frac{1}{q+1}}
			\le \left(\frac{1}{\gamma}\right)^{\frac{1}{q+1}}\mathcal{D}(t)^{\frac{2}{q+1}}.
		\end{align}
		Hence, by the mean value theorem there exist $t_1\in[t,t+\frac{1}{4}]$ and $t_2\in[t+\frac{3}{4},t+1]$ such that
		\begin{align}\label{dissip_BB}
			\|u_t(t_i)\|^{2} \le 4\int_{t}^{t+1}\|u_t(s)\|^2\,ds
			\le 4\left(\frac{1}{\gamma}\right)^{\frac{1}{q+1}}\mathcal{D}(t)^{\frac{2}{q+1}}.
		\end{align}
		
		Next, multiplying the first equation in (\ref{P}) by $u$ and integrating over $\Omega\times[t_1,t_2]$, we get
		\begin{eqnarray}
			\label{abs_1}
			\begin{aligned}
				& \int_{t_1}^{t_2}\int_{\Omega}\left[\,|\nabla u|^2+f(u)u-hu\,\right]dxds\\
				& =\; \int_{t_1}^{t_2}\|u_t(s)\|^2ds -\int_{t_1}^{t_2}g\left(||\mathbf{u}(s)||^{2}_{\mathscr{H}}\right)\int_{\Omega}u_tudxds
				-\left[\,\int_{\Omega}u_t(s)u(s)dx\,\right]_{t_1}^{t_2}.
			\end{aligned}
		\end{eqnarray}
		Recalling the definition of $\widetilde{E}$, it follows from \eqref{abs_1} that
		\begin{eqnarray}\label{diss_C}
			\begin{aligned}	&\int_{t_1}^{t_2}\widetilde{E}(\mathbf{u}(s))ds+\frac{1}{2}\int_{t_1}^{t_2}\|\nabla u(s)\|^2\\
				& \leq\;R_*+c_f|\Omega|
				+\underbrace{\frac{3}{2}\int_{t_1}^{t_2}\|u_t(s)\|^2ds}_{I_1}
				+\underbrace{\int_{t_1}^{t_2}\int_{\Omega}\left[\,F(u)-f(u)u-c_f\right]dxds}_{I_2}\\
				&\quad-\underbrace{\int_{t_1}^{t_2}g\left(||\mathbf{u}(s)||^{2}_{\mathscr{H}}\right)\int_{\Omega}u_tudxds}_{I_3}-\underbrace{\left[\,\int_{\Omega}u_t(s)u(s)dx\,\right]_{t_1}^{t_2}}_{I_4}.
			\end{aligned}
		\end{eqnarray}
		In what follows we will estimate the terms $I_1,\cdots,I_4$. First, from \eqref{abs-01}, it is straightforward that
		\begin{align*}
			I_1\le \frac{3}{2}\left(\frac{1}{\gamma}\right)^{\frac{1}{q+1}}
			\mathcal{D}(t)^{\frac{2}{q+1}}.
		\end{align*}
		Using estimate \eqref{hyp-f-2} and the embedding $H^1_0(\Omega)\hookrightarrow L^2(\Omega)$, it follows that
		\begin{align*}
			I_2\le \frac{\theta}{2\lambda_1}\int_{t_1}^{t_2}\|\nabla u(s)\|^2ds.
		\end{align*}
		Next, using \eqref{globalest} and the continuity of $g$, we obtain		$$g(||\mathbf{u}(t)||^2_{\mathscr{H}})\le \mathcal{Q}(||\mathbf{u}_0||^2_{\mathscr{H}})\leq \mathcal{Q}(R),\qquad \forall t\geq 0.$$
		Hence, from  \eqref{est-E-tilde},  we have
		\begin{align*}
			I_3&\le \;\int_{t_1}^{t_2}g\left(||\mathbf{u}(s)||^2_{\mathscr{H}}\right)\|u_t\|\|u\|ds\\
			&\le \;\left[\int_{t}^{t+1}g\left(||\mathbf{u}(s)||^{2}_{\mathscr{H}}\right)
			\|u_t(s)\|^2ds\right]^{1/2}
			\left[\,\frac{1}{\lambda_1}\int_t^{t+1}
			\underbrace{g\left(||\mathbf{u}(s)||^{2}_{\mathscr{H}}\right)}_{\le \mathcal{Q}(R)}\|\nabla u(s)\|^2ds\right]^{1/2}\\
			&\le \;\mathcal{D}(t)\frac{2\mathcal{Q}(R)^{1/2}}{\omega^{1/2}\lambda_1^{1/2}}\sup_{t\le s\le t+1}\widetilde{E}(\mathbf{u}(s))^{1/2}\\
			&\le \frac{8\mathcal{Q}(R)}{\omega\lambda_1}\mathcal{D}(t)^2+\frac{1}{8}\sup_{t\le s\le t+1}\widetilde{E}(\mathbf{u}(s)).
		\end{align*}
		Finally, from \eqref{dissip_BB}, we have
		\begin{align*}
			I_4&\le\; \frac{1}{\lambda_1^{1/2}}\sum_{i=1}^2\|u_t(t_i)\|\lambda_1^{1/2}\|u(t_i)\|\\
			&\le\; \frac{8}{\omega^{1/2}\lambda_1^{1/2}\gamma^{1/2(q+1)}}
			\mathcal{D}(t)^{1/q+1} \sup_{t\le s\le t+1}\widetilde{E}(\mathbf{u}(s))^{1/2}\\
			&\le\; \frac{128}{\omega\lambda_1\gamma^{1/q+1}}
			\mathcal{D}(t)^{2/q+1}+\frac{1}{8}\sup_{t\le s\le t+1}\widetilde{E}(\mathbf{u}(s)).
		\end{align*}
		Substituting $I_1,\cdots,I_4$ into (\ref{diss_C}), we get
		\begin{equation}\label{dissip_D}
			\begin{split}
				&	\int_{t_1}^{t_2}\widetilde{E}(\mathbf{u}(s))ds
				+\frac{1}{2}\underbrace{\left(1-\frac{\theta}{\lambda_1}\right)}_{=\omega}
				\int_{t_1}^{t_2}\|\nabla u(s)\|^2ds\\
				& \le  R_*+c_f |\Omega|+\frac{1}{4} \sup_{t\le s\le t+1}\widetilde{E}(\mathbf{u}(s))\\
				&\,  +\frac{8\mathcal{Q}(R)}{\omega\lambda_1}\mathcal{D}(t)^2
				+\left(\frac{1}{\gamma}\right)^{1/(q+1)}\left(\frac{3}{2}
				+\frac{128}{\omega\lambda_1}\right)\mathcal{D}(t)^{2/q+1}.
			\end{split}
		\end{equation}
		Again from the mean value theorem, there exists $\tau\in [t_1,t_2]$ such that
		\begin{align}\label{dissip_E}
			\int_{t_1}^{t_2}\widetilde{E}(\mathbf{u}(s))ds=\widetilde{E}(\mathbf{u}(\tau))(t_2-t_1)\ge \frac{1}{2}\widetilde{E}(\mathbf{u}(t+1)),
		\end{align}
		where we use that $\widetilde{E}$ is nonincreasing.
		Moreover, from (\ref{dissip_A}), we have
		\begin{align}\label{dissip_F}
			\widetilde{E}(\mathbf{u}(t))=\widetilde{E}(\mathbf{u}(t+1))+\mathcal{D}(t)^2.
		\end{align}
		Combining (\ref{dissip_D}), (\ref{dissip_E}), (\ref{dissip_F}), and using that $\widetilde{E}(\mathbf{u}(t))=\sup_{t\le s\le t+1}\widetilde{E}(\mathbf{u}(s))$, we obtain
		\begin{eqnarray}\label{yx}
			\begin{aligned}
				\sup_{t\le s\le t+1}\widetilde{E}(\mathbf{u}(s))\le\;4R_*
				+\frac{1}{2}\sup_{t\le s\le t+1}\widetilde{E}(\mathbf{u}(s))\mathcal{Q}(R)\mathcal{D}(t)^{2/q+1},
			\end{aligned}
		\end{eqnarray}
		where we have used
		$$\mathcal{D}(t)\le \sqrt{\left|\widetilde{E}(\mathbf{u}(t))\right|+\left|\widetilde{E}(\mathbf{u}(t+1))\right|}\le \sqrt{2\left|\widetilde{E}(\mathbf{u}_0)\right|}
		,$$
		and 		$$\left[\,\left(1+\frac{16\mathcal{Q}(R)}{\omega\lambda_1}\right)
		\mathcal{D}(t)^{2q/q+1}+\left(\frac{1}{\gamma}\right)^{1/q+1}
		\left(3+\frac{256}{\omega\lambda_1}\right)\,\right]\le \mathcal{Q}(R).$$
		So,	it follows from \eqref{yx} that
		\begin{align*}
			\sup_{t\le s\le t+1}\widetilde{E}(\mathbf{u}(s))\le 8R_*+\mathcal{Q}(R)\mathcal{D}(t)^{\frac{2}{q+1}},
		\end{align*}
		which implies
		\begin{align*}
			\sup_{t\le s\le t+1}\widetilde{E}(\mathbf{u}(s))^{q+1}\le 2^q\left(8R_*\right)^{q+1}+2^q\mathcal{Q}(R)^{q+1}\mathcal{D}(t)^2.
		\end{align*}
		Letting $C_R=2^q\mathcal{Q}(R)^{q+1}$, the definition of $\mathcal{D}(t)$ gives
		\begin{align}\label{dissip_H}
			\sup_{t\le s\le t+1}\widetilde{E}(\mathbf{u}(s))^{q+1}\le (16R_*)^{q+1}+C_R\left[\,\widetilde{E}(\mathbf{u}(t))
			-\widetilde{E}(\mathbf{u}(t+1))\,\right].
		\end{align}
		Applying Nakao's Lemma (cf. [32, Lemma 2.1]) in \eqref{dissip_H}  yields
		\begin{align}\label{L66}
			\widetilde{E}(\mathbf{u}(t))	\le  \left[\,\frac{1}{C_R}(t-1)^{+}
			+\widetilde{E}(\mathbf{u}_0)^{-q}\,\right]^{-1/q}+16R_*,
		\end{align} where $(t-1)^{+}=\max\{t-1,0\}$.
		Thus, by using \eqref{est-E-tilde}, it follows from \eqref{L66} that
		\begin{align}\label{energy-estimateI}
			||\mathbf{u}(t)||^2_{\mathscr{H}}
			\;\le\;
			\frac{4}{\omega}\biggl[\, \frac{1}{C_R} (t-1)^{+} + \widetilde{E}(\mathbf{u}_0)^{-q} \biggr]^{-1/q}
			+ \frac{64R_*}{\omega}.
		\end{align}
		Then, using that
		$K_R=\sup_{\|\mathbf{u}_0\|_{\mathscr{H}}\leq R}\widetilde{E}(\mathbf{u}_0)<\infty,$
		it follows from \eqref{energy-estimateI} that
		\begin{align}\label{Diss-3}
			\limsup_{t\rightarrow +\infty}||\mathbf{u}(t)||_{\mathscr{H}}^2	\le  \frac{64R_*}{\omega}.
		\end{align}
		Therefore, choosing $\mathcal{R}^2>\frac{64R_*}{\omega}$, we conclude, from \eqref{Diss-3} that
		$$\mathcal{B}:=\{(u,v)\in \mathscr{H}:\,||(u,v)||_{\mathscr{H}}\le \mathcal{R}\}$$
		is an absorbing set in $\mathscr{H}$,
		which completes the proof of Proposition~\ref{prop-absorbing-set}.
	\end{proof}
	
	\subsection{Asymptotic smoothness property}
	We begin with an auxiliary result that will be used in the proof of Proposition~\ref{proposition-smoothness1}.
	\begin{lemma}\cite[Lemma 4.5]{Zhou-Yang-2026}\label{lemma-Zhou_Yang}
		For any $q>0$,
		\begin{align*}
			0<\frac{a^q-b^q}{a-b}\le \max\{1,q\}\left(\frac{a^q+b^q}{a+b}\right),
		\end{align*}
		for all $a, b\geq 0$ with $a+b>0$ and $a\neq b$.
	\end{lemma}
	
	\begin{proposition}\label{proposition-smoothness1}
		Let $\mathbf{u}^1(t)=S_t\mathbf{u}^1(0)$ and $\mathbf{u}^2(t)=S_t\mathbf{u}^2(0)$ be Shatah--Struwe solutions of problem \eqref{P} associated with initial data $\mathbf{u}^1(0), \mathbf{u}^2(0)$ belonging to a bounded set $B\subset\mathscr{H}$.
		Let $z(t)=u^1(t)-u^2(t), Z=(z,z_t)$ and fix $T>0$.
		Under the assumptions of Theorem~\ref{theo-main}, there exist positive constants $C_B$ and $C_{B,T}$ such that
		\begin{equation}\label{inequality-main}
			\|Z(T)\|_{\mathscr{H}}^2
			\le \frac{C_B}{ T^{1/q+1}}
			+ C_{B,T}\sup_{t\in[0,T]}\|z(t)\|_{L^{\rho+1}(\Omega)}.
		\end{equation}
	\end{proposition}
	
	\begin{proof} Obviously, \eqref{inequality-main} holds for $\mathbf{u}^1(t)=\mathbf{u}^2(t)$.  Let $\mathbf{u}^1(t)\neq \mathbf{u}^2(t)$ be two distinct  Shatah--Struwe solutions of problem \eqref{P}.
		Consider the difference $z=u^1-u^2$ and denote $Z=(z,z_t)$. Then,
		$z$ solves the problem:
		\begin{eqnarray}
			\left\{
			\begin{array}{l}\label{QS-1}
				\displaystyle
				z_{tt} - \Delta z
				+ \underbrace{\Bigl[\,g\left(\|\mathbf{u}^1(t)\|^{2}_{\mathscr{H}}\right)u^1_t
					- g\left(\|\mathbf{u}^2(t)\|^{2}_{\mathscr{H}}\right)u^2_t\,\Bigr]}_{J}
				+ \bigl[f(u^1) - f(u^2)\bigr] = 0, \\[2mm]
				\displaystyle
				Z(0) = (z_0, z_1)
				= (u^1_0 - u^2_0,\, u^1_1 - u^2_1).
			\end{array}
			\right.
		\end{eqnarray}
		Using the expression $g(s)=\gamma s^q+\kappa$ from Assumption \ref{hyp-g}, the term $J$ can be decomposed as
		\begin{align*}
			J
			= \kappa z_t+\frac{\gamma}{2}\Bigl[\,\|\mathbf{u}^1(t)\|^{2q}_{\mathscr{H}}
			+ \|\mathbf{u}^2(t)\|^{2q}_{\mathscr{H}}\,\Bigr] z_t
			+ \frac{\gamma}{2}\Bigl[\,\|\mathbf{u}^1(t)\|^{2q}_{\mathscr{H}}
			- \|\mathbf{u}^2(t)\|^{2q}_{\mathscr{H}}\,\Bigr](u^1_t + u^2_t).
		\end{align*}
		Let $T>0$.	Multiplying   \eqref{QS-1} by $z_t$ and integrating over $\Omega\times [\tau,T]$, with $\tau\ge 0$, we get
		\begin{eqnarray}
			\begin{aligned}\label{asymptotic00}
				&\frac{1}{2}||Z(T)||^2_{\mathscr{H}}
				+\kappa\int_{\tau}^T\|z_t(s)\|^2\,ds
				+\frac{\gamma}{2}\int_{\tau}^T\left[\,||\mathbf{u}^1(s)||^{2q}_{\mathscr{H}}+||\mathbf{u}^2(s)||^{2q}_{\mathscr{H}}\,\right]\| z_t(s)\|^2ds\\ &\quad+\,\underbrace{\frac{\gamma}{2}\int_{\tau}^T\left[\,||\mathbf{u}^1(s)||^{2q}_{\mathscr{H}}-||\mathbf{u}^2(s)||^{2q}_{\mathscr{H}}\,\right]\int_{\Omega}[\,u^1_t+u^2_t\,]z_tdxds}_{J_1}\\
				&=\frac{1}{2}||Z(\tau)||^2_{\mathscr{H}}
				-\underbrace{\int_{\tau}^T\int_{\Omega}[f(u^1) - f(u^2)\bigr]z_tdxds}_{J_2}.
			\end{aligned}
		\end{eqnarray}
		Observe that whenever
		$
		\|\mathbf{u}^1\|_{\mathscr{H}}
		=
		\|\mathbf{u}^2\|_{\mathscr{H}},
		$
		the coefficient defining $J_1$ vanishes, and consequently $J_1=0$. Hence, it is sufficient to consider only the case
		$
		\|\mathbf{u}^1\|_{\mathscr{H}}
		\neq
		\|\mathbf{u}^2\|_{\mathscr{H}}.
		$
		Under this condition, we may write
		\begin{eqnarray*}
			\|\mathbf{u}^1(t)\|_{\mathscr{H}}^{2q}
			-
			\|\mathbf{u}^2(t)\|_{\mathscr{H}}^{2q}
			&=&
			\Pi(t)
			\left[\,
			\|\mathbf{u}^1(t)\|_{\mathscr{H}}^{2}
			-
			\|\mathbf{u}^2(t)\|_{\mathscr{H}}^{2}
			\,\right]
			\\
			&=&
			\Pi(t)
			\left(
			\|u_t^1(t)\|^{2}
			-
			\|u_t^2(t)\|^{2}
			+
			\|\nabla u^1(t)\|^{2}
			-
			\|\nabla u^2(t)\|^{2}
			\right),
		\end{eqnarray*}
		where
		\begin{align*}
			\Pi(t):=
			\frac{
				\|\mathbf{u}^1(t)\|_{\mathscr{H}}^{2q}
				-
				\|\mathbf{u}^2(t)\|_{\mathscr{H}}^{2q}
			}{
				\|\mathbf{u}^1(t)\|_{\mathscr{H}}^{2}
				-
				\|\mathbf{u}^2(t)\|_{\mathscr{H}}^{2}
			}
			>0,
			\qquad
			\forall\,t\ge0.
		\end{align*}
		Then, we can rewrite the term $J_1$ as follows
		\begin{eqnarray*}
			J_1
			&=&\underbrace{\frac{\gamma}{2}\int_{\tau}^T\Pi(s)\left[\,\|u^1_t(s)\|^2-\|u^2_t(s)\|^2\,\right]^2ds}_{\ge 0}\nonumber\\
			&&+\underbrace{\frac{\gamma}{2}\int_{\tau}^T\Pi(s)\left[\,\|\nabla u^1(s)\|^2-\|\nabla u^2(s)\|^2\,\right]\left[\,\|u^1_t(s)\|^2-\|u^2_t(s)\|^2\,\right]ds}_{J'_1}.
		\end{eqnarray*}
		Then, returning to \eqref{asymptotic00}, we get
		\begin{eqnarray}
			\begin{aligned}
				\label{asymptotic001}
				&\frac{1}{2}||Z(T)||^2_{\mathscr{H}}+\underbrace{\kappa\int_{\tau}^T\|z_t(s)\|^2\,ds}_{\ge 0}+\frac{\gamma}{2}\int_{\tau}^T\left[\,||\mathbf{u}^1(s)||^{2q}_{\mathscr{H}}+||\mathbf{u}^2(s)||^{2q}_{\mathscr{H}}\,\right]\| z_t(s)\|^2ds\\
				&\quad+\frac{\gamma}{2}\int_{\tau}^T\Pi(s)\left[\,\|u^1_t(s)\|^2-\|u^2_t(s)\|^2\,\right]^2ds=\frac{1}{2}||Z(\tau)||^2_{\mathscr{H}}-J'_1-J_2.
			\end{aligned}
		\end{eqnarray}
		Next, we estimate the terms $-J'_1$ and $-J_2$ appearing on the right-hand side of \eqref{asymptotic001}. Using Young's inequality, the bound
		$
		\|\nabla u^{i}\|^2 \le \|\mathbf{u}^{i}\|_{\mathscr{H}}^2,
		$
		and Lemma~\ref{lemma-Zhou_Yang} with $a=||\mathbf{u}^1||^{2}_{\mathscr{H}}$ and $b=||\mathbf{u}^2||^{2}_{\mathscr{H}}$, we obtain
		\begin{eqnarray*}
			-J'_1
			&\le&\frac{\gamma}{2}\int_{\tau}^T\Pi(s)
			\left[\,\|u^1_t(s)\|^2-\|u^2_t(s)\|^2\,\right]^2ds\\
			&&+\frac{\gamma}{8}\int_{\tau}^T\Pi(s)\left[\,\|\nabla u^1(s)\|+\|\nabla u^2(s)\|\,\right]^2\|\nabla z(s)\|^2\\
			&\le&\frac{\gamma}{2}\int_{\tau}^T\Pi(s)
			\left[\,\|u^1_t(s)\|^2-\|u^2_t(s)\|^2\,\right]^2ds\\
			&&\qquad+\frac{\gamma}{4}\int_{\tau}^T\Pi(s)\left[\, ||\mathbf{u}^1(s)||^2_{\mathscr{H}}+||\mathbf{u}^2(s)||^2_{\mathscr{H}}\,\right]\|\nabla z(s)\|^2ds\\
			&\le&\frac{\gamma}{2}\int_{\tau}^T\Pi(s)
			\left[\,\|u^1_t(s)\|^2-\|u^2_t(s)\|^2\,\right]^2ds\\
			&&\qquad+\underbrace{ \frac{\gamma\max\{1,q\}}{4}\int_{\tau}^T\left[\,||\mathbf{u}^1(s)||_{\mathscr{H}}^{2q}+||\mathbf{u}^2(s)|_{\mathscr{H}}^{2q}\,\right]\|\nabla z(s)\|^2ds}_{J''_1}.
		\end{eqnarray*}
		From Assumption \eqref{hyp_f'}, \eqref{globalest} with $(r, s)=(4, 12)$,  H\"older's inequality with $\frac{\rho-1}{2(\rho+1)}+\frac{1}{\rho+1}+\frac{1}{2}=1$, embedding $L^{12}(\Omega)\hookrightarrow L^{2(\rho+1)}(\Omega)$, and Young's inequality with $\frac{\rho-1}{4}+\frac{5-\rho}{4}=1$, we have
		\begin{eqnarray*}-J_2&=&-\int_{\tau}^T\int_{\Omega} \int_0^1f'(\theta u^1+(1-\theta)u^2)d\theta zz_tdxds\\
			&\le&C\int_{\tau}^T\left(1+\|u^1(s)\|^{\rho-1}_{L^{2(\rho+1)}(\Omega)}
			+\|u^2(s)\|^{\rho-1}_{L^{2(\rho+1)}(\Omega)}\right)\|z(s)\|_{L^{\rho+1}
				(\Omega)}
			\|z_t(s)\|ds\\&\le&C\int_{\tau}^T\left(1
			+\|u^1(s)\|^{4}_{L^{12}(\Omega)}+\|u^2(s)\|^{4}_{L^{12}(\Omega)}\right)
			\|z_t(s)\|ds\sup_{t\in [0,T]}\|z(t)\|_{L^{\rho+1}
				(\Omega)}\\
			&\le&C_{B,T}\sup_{t\in [0,T]}\|z(t)\|_{L^{\rho+1}(\Omega)}.
		\end{eqnarray*}
		Thus, substituting the above estimates for $-J'_1$ and $-J_2$ into \eqref{asymptotic001}, we obtain
		\begin{eqnarray}
			\begin{aligned}
				\label{asymp_1}
				&	 \frac{1}{2}||Z(T)||^2_{\mathscr{H}}+\frac{\gamma}{2}\int_{\tau}^T\left[\,||\mathbf{u}^1(s)||^{2q}_{\mathscr{H}}+||\mathbf{u}^2(s)||^{2q}_{\mathscr{H}}\,\right]\| z_t(s)\|^2ds
				\\
				&\le\,\frac{1}{2}||Z(\tau)||^2_{\mathscr{H}}+C_{B, T}\sup_{t\in [0,T]}\|z(t)\|_{L^{\rho+1}(\Omega)}+J''_1.
			\end{aligned}
		\end{eqnarray}
		It remains to estimate the term $J''_1$. Multiplying \eqref{QS-1} by $\left[\,||\mathbf{u}^1(t)||^{2q}_{\mathscr{H}}+||{\mathbf{u}^2(t)}||^{2q}_{\mathscr{H}}\,\right]z$
		and integrating over $\Omega\times [\tau,T]$, we have
		\begin{eqnarray}
			\label{asymptotic2}
			\begin{aligned}
				&\int_{\tau}^T\left[\,||\mathbf{u}^1(s)||^{2q}_{\mathscr{H}}+||\mathbf{u}^2(s)||^{2q}_{\mathscr{H}}\,\right]\|\nabla z(s)\|^2ds\\
				&=\int_{\tau}^T\left[\,||\mathbf{u}^1(s)||^{2q}_{\mathscr{H}}+||\mathbf{u}^2(s)||^{2q}_{\mathscr{H}}\,\right]\|z_t(s)\|^2ds+\sum_{i=1}^5L_i,
			\end{aligned}
		\end{eqnarray}
		where
		\begin{eqnarray*}
			L_1&=& -\,\int_{\tau}^T \left[\,||\mathbf{u}^1(s)||^{2q}_{\mathscr{H}}+||\mathbf{u}^2(s)||^{2q}_{\mathscr{H}}\,\right]\int_{\Omega}[f(u^1) - f(u^2)\bigr]zdxds,\\
			L_2&=&-\,\frac{\gamma}{2}\int_{\tau}^T \left[\,||\mathbf{u}^1||^{2q}_{\mathscr{H}}+||\mathbf{u}^2||^{2q}_{\mathscr{H}}\,\right]\left[\,\frac{2\kappa}{\gamma}+||\mathbf{u}^1||^{2q}_{\mathscr{H}}+||\mathbf{u}^2||^{2q}_{\mathscr{H}}\,\right]\int_{\Omega}z_t z dxds,\\
			L_3&=&-\,\frac{\gamma}{2}\int_{\tau}^T \left[\,||\mathbf{u}^1||^{2q}_{\mathscr{H}}+||\mathbf{u}^2||^{2q}_{\mathscr{H}}\,\right]\left[\,||\mathbf{u}^1||^{2q}_{\mathscr{H}}-||\mathbf{u}^2||^{2q}_{\mathscr{H}}\,\right]\int_{\Omega}\left[\,
			u^1_t+ u^2_t\,\right]zdxds,\\
			L_4&=&\,\int_{\tau}^T\left\{\frac{d}{ds} \left[\,||\mathbf{u}^1(s)||^{2q}_{\mathscr{H}}+||\mathbf{u}^2(s)||^{2q}_{\mathscr{H}}\,\right]\cdot\int_{\Omega}z_tzdx\right\}ds,\\
			L_5&=&-\left\{\,\left[\,||\mathbf{u}^1(s)||^{2q}_{\mathscr{H}}+||\mathbf{u}^2(s)||^{2q}_{\mathscr{H}}\,\right]\int_{\Omega}z_tzdx\,\right\}_{\tau}^T.
		\end{eqnarray*}
		Combining (\ref{asymp_1}) and (\ref{asymptotic2}),
		we obtain
		\begin{eqnarray}
			\begin{aligned}
				\label{asymptotic3}
				& \frac{1}{2}||Z(T)||^2_{\mathscr{H}}+\frac{\gamma}{2}\left(1-\frac{\max\{1,q\}}{2}\right)\int_{\tau}^T\left[\,||\mathbf{u}^1(s)||^{2q}_{\mathscr{H}}+||\mathbf{u}^2(s)||^{2q}_{\mathscr{H}}\,\right]\| z_t(s)\|^2ds
				\\
				&\le\,\frac{1}{2}||Z(\tau)||^2_{\mathscr{H}}+C_{B, T}\sup_{t\in [0,T]}\|z(t)\|_{L^{\rho+1}(\Omega)}+\frac{\gamma\max\{1,q\}}{4}\sum_{i=1}^5L_i.
			\end{aligned}
		\end{eqnarray}
		Using that\begin{eqnarray*}||\mathbf{u}^1(t)||^{2q}_{\mathscr{H}}+||\mathbf{u}^2(t)||^{2q}_{\mathscr{H}}\ge \|u^1_t(t)\|^{2q}+\|u^2_t(t)\|^{2q}\ge \frac{1}{2^{2q-1}} \|z_t(t)\|^{2q}.
		\end{eqnarray*}
		Then, it follows from \eqref{asymptotic3} that
		\begin{eqnarray}\label{asymptotic4}
			\begin{aligned}
				||Z(T)||^2_{\mathscr{H}}+\mu_q\int_{\tau}^T\| z_t(s)\|^{2(q+1)}ds\le&\;||Z(\tau)||^2_{\mathscr{H}}+C_{B,T}\sup_{t\in [0,T]}\|z(t)\|_{L^{\rho+1}(\Omega)}\\
				&+\frac{\gamma\max\{1,q\}}{2}\sum_{i=1}^5L_i,
			\end{aligned}
		\end{eqnarray}
		where
		\begin{align}\label{const-mu}
			\mu_q:=\gamma
			\left(1-\frac{\max\{1,q\}}{2}\right)\frac{1}{2^{2q-1}}\ge0,
			\qquad q\in[1/2,2].
		\end{align}
		Moreover, by Assumption~\eqref{hyp_f'}, estimate~\eqref{globalest},
		H\"older's inequality with
		$
		\frac{\rho-1}{\rho+1}+\frac{2}{\rho+1}=1,
		$
		and the Sobolev embeddings
		$
		H_0^1(\Omega)\hookrightarrow L^6(\Omega)
		\hookrightarrow L^{\rho+1}(\Omega),
		$
		it follows that
		\begin{eqnarray*}
			L_1&\le& C_B\left|\int_{\tau}^T\int_{\Omega}[f(u^1) - f(u^2)\bigr]zdxds\right|\\
			&\le &C_B\int_{\tau}^{T}\int_{\Omega}\left[1+|u^1|^{\rho-1}+|u^2|^{\rho-1}\right]|z|^2dx\\
			&\le&C_B\int_{\tau}^T\left[\,1+\|u^1(s)\|^{\rho-1}_{L^{\rho+1}(\Omega)}
			+\|u^2(s)\|^{\rho-1}_{L^{\rho+1}(\Omega)}\,\right]\|z(s)\|^2_{L^{\rho+1}(\Omega)}ds\\
			&\le&\left(C_B\int_{\tau}^T\left[\,1+\|\nabla u^1(s)\|^4+\|\nabla u^2(s)\|^4\,\right]\|\nabla z(s)\|ds\right)\sup_{t\in[0,T]}\|z(s)\|_{L^{\rho+1}(\Omega)}\\
			&\le&C_{B,T}\sup_{t\in[0,T]}\|z(s)\|_{L^{\rho+1}(\Omega)}.
		\end{eqnarray*}
		The terms $L_2,L_3$ and $L_5$ can be estimated by
		\begin{align*}L_2,L_3,L_5\le C_{B,T}\sup_{t\in [0,T]}\|z(s)\|\le C_{B,T}\sup_{t\in [0,T]}\|z(s)\|_{L^{\rho+1}(\Omega)}.
		\end{align*}
		Next, we estimate the term $L_4$. By H\"older inequality, we have
		\begin{eqnarray*}
			|L_4|
			&=&
			\left|
			\int_\tau^T
			\frac{d}{ds}
			\Bigl(
			\|\mathbf u^1(s)\|_{\mathscr{H}}^{2q}
			+
			\|\mathbf u^2(s)\|_{\mathscr{H}}^{2q}
			\Bigr)
			\left(\int_\Omega z_tz\,dx\right)ds
			\right|\\
			&\le&
			\int_\tau^T
			\sum_{i=1}^{2}
			\left|
			\frac{d}{ds}
			\|\mathbf u^i(s)\|_{\mathscr{H}}^{2q}
			\right|
			\,
			\|z_t(s)\|\,\|z(s)\|
			\,ds.
		\end{eqnarray*}
		We claim that
		\begin{align}\label{a-0}
			\left|
			\frac{d}{dt}
			\|\mathbf u^i(t)\|_{\mathscr{H}}^{2q}
			\right|
			\le
			C_B
			\left[\,
			1+\|u^i(t)\|^4_{L^{12}(\Omega)}
			\,\right], \qquad i=1,2.
		\end{align}
		Since  \eqref{a-0} is immediate when
		$\|\mathbf u^i(t)\|_{\mathscr{H}}=0$, we may assume in the sequel that
		$\|\mathbf u^i(t)\|_{\mathscr{H}}>0$. Multiplying \eqref{P} by $u_t^i$ yields
		\begin{equation}\label{energy-derivative}
			\frac12
			\frac{d}{dt}
			\|\mathbf u^i(t)\|_{\mathscr{H}}^{2}
			=
			-\gamma
			\|\mathbf u^i(t)\|_{\mathscr{H}}^{2q}
			\|u_t^i(t)\|^{2}
			-\kappa \|u_t^i(t)\|^{2}
			-\int_\Omega f(u^i)u_t^i\,dx
			+\int_\Omega hu_t^i\,dx.
		\end{equation}
		From Assumption \eqref{hyp_f'}, H\"older's inequality with
		$
		\frac{\rho-1}{12}
		+
		\frac{7-\rho}{12}
		+
		\frac12
		=1,
		$
		and the embeddings
		$
		H_0^1(\Omega)
		\hookrightarrow
		L^6(\Omega)
		\hookrightarrow
		L^{\frac{12}{7-\rho}}(\Omega),
		$
		we obtain
		$$
		\left|
		\int_\Omega f(u^i)u_t^i\,dx
		\right|
		\le
		C
		\left[\,
		1+\|u^i(t)\|_{L^{12}(\Omega)}^{4}
		\,\right]
		\|\mathbf u^i(t)\|_{\mathscr{H}}^{2}.
		$$
		Moreover,
		$$
		\left|
		\int_\Omega hu_t^i\,dx
		\right|
		\le
		\|h\|\,\|u_t^i(t)\|
		\le
		\|h\|\,\|\mathbf u^i(t)\|_{\mathscr{H}}.
		$$
		Using the above estimates in \eqref{energy-derivative}, we obtain
		$$
		\|\mathbf u^i(t)\|_{\mathscr{H}}
		\left|\frac{d}{dt}\|\mathbf u^i(t)\|_{\mathscr{H}}\right|
		\le
		C_B
		\left[\,
		1+\|u^i(t)\|_{L^{12}(\Omega)}^{4}
		\,\right]
		\|\mathbf u^i(t)\|_{\mathscr{H}}.
		$$ 
		This implies that
		\begin{equation*}
			\left|\frac{d}{dt}\|\mathbf u^i(t)\|_{\mathscr{H}}\right|
			\le
			C_B
			\left[\,
			1+\|u^i(t)\|_{L^{12}(\Omega)}^{4}
			\,\right].
		\end{equation*}
		Consequently, recalling that $q\ge \frac12$, the function
		$r\mapsto r^{2q}$ belongs to $C^1([0,\infty))$. Therefore, by the chain rule,
		\begin{equation*}
			\left|
			\frac{d}{dt}
			\|\mathbf u^i(t)\|_{\mathscr{H}}^{2q}
			\right|
			=
			2q\,
			\|\mathbf u^i(t)\|_{\mathscr{H}}^{2q-1}
			\left|
			\frac{d}{dt}
			\|\mathbf u^i(t)\|_{\mathscr{H}}
			\right|
			\le
			C_B
			\left[\,
			1+\|u^i(t)\|_{L^{12}(\Omega)}^{4}
			\,\right].
		\end{equation*}
		This proves that inequality \eqref{a-0} holds. Therefore, the term $L_4$ can be estimated as
		\begin{eqnarray*}
			|L_4|
			&\le&
			C_B
			\int_\tau^T
			\left[\,
			1+\|u^1(s)\|_{L^{12}(\Omega)}^{4}
			+\|u^2(s)\|_{L^{12}(\Omega)}^{4}
			\,\right]
			\|z_t(s)\|\,\|z(s)\|
			\,ds\\
			&\le&
			C_B
			\left\{
			\int_\tau^T
			\left[\,
			1+\|u^1(s)\|_{L^{12}(\Omega)}^{4}
			+\|u^2(s)\|_{L^{12}(\Omega)}^{4}
			\,\right]
			\|z_t(s)\|\,ds
			\,\right\}
			\sup_{t\in[0,T]}\|z(t)\|\\
			&\le&
			C_{B,T}
			\sup_{t\in[0,T]}
			\|z(t)\|_{L^{\rho+1}(\Omega)}.
		\end{eqnarray*}
		Substituting $L_1,\cdots,L_5$ into \eqref{asymptotic4}, we get
		\begin{align}
			\label{asymptotic5}
			||Z(T)||^2_{\mathscr{H}}+\mu_q\int_{\tau}^T\| z_t(s)\|^{2(q+1)}ds\le\,||Z(\tau)||^2_{\mathscr{H}}
			+C_{B,T}\sup_{t\in[0,T]}\|z(s)\|_{L^{\rho+1}(\Omega)}.
		\end{align}
		Next, integrating \eqref{asymptotic5} from $0$ to $T$, we have
		\begin{align}
			\label{asymp-a}
			T||Z(T)||^2_{\mathscr{H}}+&\mu_q\int_0^T\int_{\tau}^T\| z_t(s)\|^{2(q+1)}dsd\tau\nonumber\\
			&\le\,\int_{0}^T||Z(\tau)||^2_{\mathscr{H}}d\tau
			+C_{B,T}\sup_{t\in[0,T]}\|z(s)\|_{L^{\rho+1}(\Omega)}.
		\end{align}
		Now, multiplying  \eqref{QS-1} by $z$ and integrating over $\Omega\times [0,T]$, we obtain
		\begin{align}\label{asymptotic7}
			\int_{0}^T||Z(s)||^2_{\mathscr{H}}ds=2\int_{0}^T\|z_t(s)\|^2ds+\sum_{i=1}^4\widetilde{L}_{i},
		\end{align}
		where
		\begin{eqnarray*}
			\widetilde{L}_1&=& -\,\int_0^T\int_{\Omega}\left[f(u^1)-f(u^2)\right]zdxds,\\
			\widetilde{L}_2&=&-\,\frac{\gamma}{2}\int_0^T
			\left[\,\frac{2\kappa}{\gamma}+||\mathbf{u}^1(s)||_{\mathscr{H}}^{2q}
			+||\mathbf{u}^2(s)||_{\mathscr{H}}^{2q}\,\right]\int_{\Omega}z_t z dxds,\\
			\widetilde{L}_3&=&-\,\frac{\gamma}{2}\int_0^T\left[\,||\mathbf{u}^1(s)||_{\mathscr{H}}^{2q}-||\mathbf{u}^2(s)||_{\mathscr{H}}^{2q}\,\right]\int_{\Omega}\left[\,u^1_t+u^2_t\,\right]zdxds,\\	 \widetilde{L}_4&=&-\,\ \left[\int_{\Omega}z_tzdx\right]^T_0.
		\end{eqnarray*}
		Since all the terms $\widetilde{L}_i$ are of lower order, we have
		\begin{align*}
			\sum_{i=1}^4\widetilde{L}_{i}\le C_{B,T}\sup_{t\in [0,T]}\|z(t)\|_{L^{\rho+1}(\Omega)}.
		\end{align*}
		Hence, returning to \eqref{asymptotic7}, we get
		\begin{align}\label{asymptotic8}
			\int_{0}^T||Z(s)||^2_{\mathscr{H}}ds\le2\int_{0}^T\|z_t(s)\|^2ds+C_{B,T}\sup_{t\in [0,T]}\|z(t)\|_{L^{\rho+1}(\Omega)}.
		\end{align}
		Combining \eqref{asymp-a} and \eqref{asymptotic8}, we obtain
		\begin{align*}
			T\|Z(T)\|_{\mathscr{H}}^{2}
			\le
			2\int_{0}^{T}\|z_t(s)\|^{2}\,ds
			+
			C_{B,T}\sup_{s\in[0,T]}
			\|z(s)\|_{L^{\rho+1}(\Omega)}.
		\end{align*}
		Applying H\"older's inequality yields
		\begin{align}
			\label{asymptotic10}
			T\|Z(T)\|_{\mathscr{H}}^{2}
			\le
			2T^{\frac{q}{q+1}}
			\left(
			\int_0^T \|z_t(s)\|^{2(q+1)}\,ds
			\right)^{\frac1{q+1}}
			+
			C_{B,T}\sup_{s\in[0,T]}
			\|z(s)\|_{L^{\rho+1}(\Omega)}.
		\end{align}
		By \eqref{est-u_t-a},
		\[
		\int_0^T \|u_t^i(s)\|^{2(q+1)}\,ds\le C_B,
		\qquad i=1,2,
		\]
		it follows that
		\[
		\int_0^T \|z_t(s)\|^{2(q+1)}\,ds\le C_B.
		\]
		Substituting this estimate into \eqref{asymptotic10}, we arrive at
		\begin{align}
			\label{asymptotic11}
			T\|Z(T)\|_{\mathscr{H}}^{2}
			\le
			C_B\,T^{\frac{q}{q+1}}
			+
			C_{B,T}\sup_{s\in[0,T]}
			\|z(s)\|_{L^{\rho+1}(\Omega)}.
		\end{align}
		Therefore, inequality \eqref{asymptotic11} yields \eqref{inequality-main}, which completes the proof.
	\end{proof}
	
	\begin{Corollary}
		\label{corollary-smootness2}
		Assume the hypotheses of Proposition \ref{proposition-smoothness1}. Then, the dynamical system $(\mathscr{H},S_t)$ is asymptotically smooth.
	\end{Corollary}
	\begin{proof}
		Using inequality \eqref{inequality-main} and noting that $$||Z(T)||_{\mathscr{H}}
		=||S_T\mathbf{u}^1(0)-S_T\mathbf{u}^2(0)||^2_{\mathscr{H}},$$
		we have
		\begin{align*}			||S_T\mathbf{u}^1(0)-S_T\mathbf{u}^2(0)||^2_{\mathscr{H}}\le&
			\;\frac{C_B}{ T^{1/q+1}}+C_{B,T}\sup_{t\in[0,T]}\|u^1(t)-u^2(t)\|_{L^{\rho+1}(\Omega)}.
		\end{align*}
		Note that, for any fixed $\varepsilon >0$, there exists  $T=T_{B}>0$ such that $\frac{C_B}{ T^{1/q+1}}<\varepsilon$ and
		\begin{align*}
			||S_T\mathbf{u}^1(0)-S_T\mathbf{u}^2(0)||^2_{\mathscr{H}}\le \varepsilon +\Psi_{B,T}(\mathbf{u}^1(0),\mathbf{u}^2(0)),\quad \mbox{for}\quad T\ge T_B,
		\end{align*}
		where the functional  $\Psi_{B,T}:B\times B\rightarrow \mathbb{R}$ is defined by
		\begin{align*}
			\Psi_{B,T}(\mathbf{u}^1(0),\mathbf{u}^2(0))
			:=C_{B,T}\sup_{t\in[0,T]}\|u^1(t)-u^2(t)\|_{L^{\rho+1}(\Omega)}.
		\end{align*}
		Let $\mathbf{u}^m=(u^m,u^m_t)$ be the Shatah--Struwe solution of \eqref{P} with $\mathbf{u}^m(0)=(u^m_0,u^m_1)\in B$, $m\in \mathbb{N}$.
		Thus, it follows immediately from the compact embedding
		$H^1_0(\Omega)\hookrightarrow\hookrightarrow L^{\rho+1}(\Omega) \hookrightarrow L^{2}(\Omega)$ [for $\rho<5$] and the Aubin-Lions-Simon compactness theorem  that $\{u^m\}_{m\in\mathbb{N}}$ is relatively compact in  $C([0,T];L^{\rho+1}(\Omega))$, which implies
		\begin{equation*}			\liminf_{m\rightarrow\infty}\liminf_{k\rightarrow\infty}
			\Psi_{B,T}(\mathbf{u}^M(0),\mathbf{u^k}(0))=\liminf_{m\rightarrow \infty}\liminf_{k\rightarrow \infty}C_{B,T}\sup_{t\in [0,T]}\| u^m(t)-u^k(t)\|_{L^{\rho+1}(\Omega)}=0.
		\end{equation*}
		Hence the dynamical system $(\mathscr{H},S_t)$ is  {\bf asymptotically smooth}.
	\end{proof}
	
	\subsection{Kolmogorov \texorpdfstring{$\varepsilon$}--entropy of the global attractor}
	In conclusion, this section will demonstrate that utilizing estimate \eqref{asymptotic11} allows us to derive an estimate for the Kolmogorov $\varepsilon$-entropy of the global attractor $\mathfrak{A}$ associated with the dynamical system $(\mathscr{H}, S_t)$. We  first recall the definition of the Kolmogorov $\varepsilon$-entropy for a compact set $\mathcal{M}$ in a Hilbert space $H$ (see \cite{chueshov-lasiecka-2005, Kolmogorov,Tikhomirov}).
	\begin{definition}\rm
		The  {\it Kolmogorov  $\varepsilon$-entropy} ${H}_{\varepsilon}(\mathcal{M})$ of a compact set $\mathcal{M}\subset H$   is given by
		\begin{align*}
			{H}_{\varepsilon}(\mathcal{M})=\ln N(\mathcal{M}, \varepsilon), \quad \varepsilon>0,
		\end{align*}	
		where $N(\mathcal{M}, \varepsilon)$ is the minimal number of closed sets of diameter not greater than $2 \varepsilon $ needed to cover $\mathcal{M}.$  The {\it fractal dimension} $\operatorname{dim}_{f} \mathcal{M}$ of $\mathcal{M}$ is defined by the formula
		$$
		\mbox{dim}^{f}_{H}\mathcal{M}=\limsup _{\varepsilon \rightarrow 0} \frac{{H}_{\varepsilon}(\mathcal{M})}{\ln (1 / \varepsilon)}.
		$$
	\end{definition}
	The following result from \cite{chueshov-lasiecka-2005} provides an estimate for the Kolmogorov \( \varepsilon \)-entropy \( H_{\varepsilon}(\mathcal{M}) \) of a compact set \( \mathcal{M} \subset H \).
	\begin{theorem}[{\cite[Theorem 4.2]{chueshov-lasiecka-2005}}]\label{theo-kolmog-entrop}
		Let $H$ be a separable Hilbert space and $\mathcal{M}$ be a bounded closed set in $H .$ Assume that there exists a mapping $V: \mathcal{M} \mapsto H$ such that:
		\begin{itemize}
			\item[$ 1. $] $\mathcal{M} \subseteq V \mathcal{M};$
			\item[$ 2. $]  $V$ is Lipschitz on $\mathcal{M},$ that is, there exists $L>0$ such that
			$$
			\big{\|}V z_{1}-V z_{2}\big{\|} \leq L\big{\|}z_{1}-z_{2}\big{\|}, \quad z_{1}, z_{2} \in \mathcal{M};
			$$
			\item[$ 3. $] There exist pseudometrics $\varrho_{1}$ and $\varrho_{2}$ on $H$ such that
			$$
			\big{\|}V z_{1}-V z_{2}\big{\|} \leq \varphi\big{(}\big{\|}z_{1}-z_{2}\big{\|}\big{)}+h\big{(}\big{[}\varrho_{1}\big{(}z_{1}, z_{2}\big{)}^{2}+\varrho_{2}\big{(}V z_{1}, V z_{2}\big{)}^{2}\big{]}^{1 / 2}\big{)}
			$$
			for all $z_{1}, z_{2} \in \mathcal{M},$ where $\varphi, h: \mathbb{R}^{+} \to \mathbb{R}^{+}$ are continuous non-decreasing functions such that
			$$
			\varphi(0)=0, \ \varphi(s)<s, \, s>0, \ s-\varphi(s) \text { is nondecreasing,}
			$$
			and the function $h(s)$ is strictly increasing in the interval $\big{[}0, s_{0}\big{]}$ for some $s_{0}>0$ with  $h(0)=0.$
			\item[$ 4. $]For any $q>0$ and for any closed bounded set $B \subset \mathcal{M}$ the maximal number $m(B, \psi)$ of elements $x_{j}^{B} \in B$ such
			$$
			\varrho_{1}\big{(}x_{j}^{B}, x_{i}^{B}\big{)}^{2}+\varrho_{2}\big{(}V x_{j}^{B}, V x_{i}^{B}\big{)}^{2}>\psi^{2},  \ \ i \neq j,  \, i, j=1, \ldots, m(B, \psi),
			$$
			is finite.
		\end{itemize}
		Then $\mathcal{M}$ is a compact set and there exists $0<\varepsilon_{0}<1$ such that for all $\varepsilon \leq \varepsilon_{0}<1,$ Kolmogorov's $\varepsilon$-entropy ${H}_{\varepsilon}(\mathcal{M})$ admits the following estimate
		$$
		{H}_{\varepsilon}(\mathcal{M}) \leq \int_{\varepsilon}^{\varepsilon_{0}} \frac{\ln m\big{(}g_{\delta}^{-1}(s), \psi(s)\big{)}}{s-g_{\delta}(s)} d s+H_{g_{0}\big{(}\varepsilon_{0}\big{)}}(\mathcal{M}),
		$$
		where $\varphi_{\delta}(s)=\frac{1-\delta}{2} \varphi(2 s)+\delta s$ with arbitrary $\delta \in(0,1),$ the function $\psi(s)$ is defined by the formula
		$$
		\psi(s)=\frac{1}{2} h^{-1}\{\delta[2 s-\varphi(2 s)]\}, 0<s<\varepsilon_{0},
		$$
		and
		$$
		m(r, \psi)=\sup \{m(B, \psi): B \subseteq \mathcal{M}, \operatorname{diam} B \leq 2 r\}.
		$$
	\end{theorem}
	We are now prepared to demonstrate the result that provides an estimate for the Kolmogorov \( \varepsilon \)-entropy of the global attractor \( \mathfrak{A} \).
	\begin{theorem}[{\bf Kolmogorov $\varepsilon$-Entropy}]
		Assumptions the hypotheses of Theorem~\ref{theo-main} with
		$0<q<2$. Then there exists $\varepsilon_0\in(0,1)$ such that, for every
		$\varepsilon\le\varepsilon_0$ and every $\delta\in(0,1)$, the
		Kolmogorov $\varepsilon$-entropy of the global attractor $\mathfrak A$
		satisfies
		\begin{equation}\label{eps-entrop1}
			H_\varepsilon(\mathfrak A)
			\le\int_\varepsilon^{\varepsilon_0}
			\frac{\ln m\bigl(\varphi_\delta^{-1}(s),\psi(s)\bigr)}{s-\varphi_{\delta}(s)}\,ds
			+
			H_{\varphi_\delta(\varepsilon_0)}(\mathfrak A),
		\end{equation}
		where
		\[
		\varphi_\delta(s)=\frac{1-\delta}{2}\varphi(2s)+\delta s,
		\qquad
		\psi(s)=\frac{1}{2}\bigl[\delta(2s-\varphi(2s))\bigr]^2,
		\qquad 0<s<\varepsilon_0,
		\]
		and
		\[
		\varphi(s)=\sqrt{G^{-1}(s^2)}, \qquad s\ge0.
		\]
		with the function $G$ defined in \eqref{def-G-proof}. Moreover,
		\[
		m(r,a)=\sup\Bigl\{
		m(B,a):
		B\subset\mathfrak A,\;
		\operatorname{diam}(B)\le 2r
		\Bigr\},
		\]
		where $m(B,a)$ denotes the maximal number of elements
		$\{u_j^B\}\subset B$ such that
		\[
		\varrho\bigl(S_{T_\star}u_j^B,\;S_{T_\star}u_i^B\bigr)>a,
		\qquad i\neq j,
		\]
		and $\varrho$ is a pseudometric on $\mathscr{H}$.
	\end{theorem}
	\begin{proof}
		Let
		$
		S_t\mathbf{u}^i(0)=(u^i(t),u_t^i(t)),$ $i=1,2,
		$
		be two trajectories of  problem \eqref{P} with initial data
		$\mathbf{u}^i(0)\in\mathfrak A$. Denote by
		$
		z=u^1-u^2,$ $Z=(z,z_t),
		$
		the difference of the corresponding solutions. Since $\mathfrak A$ is compact and invariant, it follows that
		$
		S_t\mathbf u^1(0),\, S_t\mathbf u^2(0)\in\mathfrak A,
		$ $t\ge 0.
		$
		Arguing as in Proposition~\ref{proposition-smoothness1} (see \eqref{asymptotic10}), there exists a time $T_\star(\mathfrak{A})>0$
		and a constant $C_{\mathfrak A,T_\star}>0$ such that
		\begin{align}\label{kom-1}
			\frac{T_{\star}^{\frac{1}{q+1}}}{2}\|Z(T_\star)\|_{\mathscr{H}}^{2}
			\le
			\left(
			\int_0^{T_\star} \|z_t(s)\|^{2(q+1)}\,ds
			\right)^{\frac1{q+1}}
			+
			C_{\mathfrak A,T_\star}
			\sup_{s\in[0,T_\star]}\|z(s)\|_{L^{\rho+1}}.
		\end{align}
		Using \eqref{kom-1} and the inequality
		$(a+b)^{q+1}\le 2^q(a^{q+1}+b^{q+1})$, we obtain
		\begin{align}\label{kom-2}
			\frac{T_\star}{2^{2q+1}}\|Z(T_\star)\|_{\mathscr{H}}^{2(q+1)}
			\le \int_0^{T_\star} \|z_t(s)\|^{2(q+1)}\,ds
			+ C_{\mathfrak A,T_\star}
			\sup_{s\in[0,T_\star]}\|z(s)\|^{q+1}_{L^{\rho+1}}.
		\end{align}
		By \eqref{asymptotic5} with $\tau=0$ and the positivity of
		$\mu_q$ for $q\in(0,2)$ (cf. \eqref{const-mu}), we infer that
		\[
		\mu_q\int_0^{T_\star} \|z_t(s)\|^{2(q+1)}ds
		\le
		\bigl(
		\|Z(0)\|_{\mathscr{H}}^2-\|Z(T_\star)\|_{\mathscr{H}}^2
		\bigr)
		+
		C_{\mathfrak A,T_\star}
		\sup_{s}\|z(s)\|_{L^{\rho+1}}.
		\]
		Substituting this estimate into \eqref{kom-2}, we obtain
		\begin{align*}
			\|Z(T_\star)\|_{\mathscr{H}}^2
			+
			\frac{\mu_qT_\star}{2^{2q+1}}
			\|Z(T_\star)\|_{\mathscr{H}}^{q+1}
			\le
			\|Z(0)\|_{\mathscr{H}}^2
			+
			C_{\mathfrak A,T_\star}
			\sup_t\|z(t)\|_{L^{\rho+1}}.
		\end{align*}
		Define
		\begin{align}\label{def-G-proof}
			G(s):=s+\frac{\mu_qT_\star}{2^{2q+1}}s^{q+1}.
		\end{align}
		Then
		\[
		G(\|Z(T_\star)\|_{\mathscr{H}}^2)
		\le
		\|Z(0)\|_{\mathscr{H}}^2
		+
		C_{\mathfrak A,T_\star}
		\sup_t\|z(t)\|_{L^{\rho+1}}.
		\]
		Applying $G^{-1}$,
		\[
		\|Z(T_\star)\|_{\mathscr{H}}^2
		\le
		G^{-1}(\|Z(0)\|_{\mathscr{H}}^2)
		+
		C_{\mathfrak A,T_\star}
		\sup_t\|z(t)\|_{L^{\rho+1}}.
		\]
		Taking square roots,
		\[
		\|S_{T_\star}\mathbf u^1(0)-S_{T_\star}\mathbf u^2(0)\|_{\mathscr{H}}
		\le
		\varphi(\|\mathbf u^1(0)-\mathbf u^2(0)\|_{\mathscr{H}})
		+
		C_{\mathfrak A,T_\star}
		\sup_t\|z(t)\|_{L^{\rho+1}}^{1/2}.
		\]
		By the continuous dependence estimate,
		\[
		\|S_{T_\star}\mathbf u^1(0)-S_{T_\star}\mathbf u^2(0)\|_{\mathscr{H}}
		\le
		K_{\mathfrak A,T_\star}
		\|\mathbf u^1(0)-\mathbf u^2(0)\|_{\mathscr{H}}.
		\]
		Hence the map $V=S_{T_\star}$ satisfies the assumptions of
		Theorem~\ref{theo-kolmog-entrop} with
		\[
		\varphi(s)=\sqrt{G^{-1}(s^2)},\qquad
		h(s)=s^{1/2},\qquad
		\varrho_2(S_{T_\star}u^1,S_{T_\star}u^2)
		=
		C_{\mathfrak A,T_\star}
		\sup_{t\in[0,T_\star]}\|u^1(t)-u^2(t)\|_{L^{\rho+1}}.
		\]
		The compact embedding $H_0^1(\Omega)\hookrightarrow L^{\rho+1}(\Omega)$
		ensures that $\varrho_2$ is a compact seminorm, so all assumptions of
		Theorem~\ref{theo-kolmog-entrop} hold. Applying that theorem yields estimate \eqref{eps-entrop1}.
	\end{proof}
	
	\section{Properties of the Global Attractor}
	To conclude the paper, we show that, under the hypotheses of Theorem~\ref{theo-main} together with the additional assumptions that $g(0)=\kappa>0$ and $f \in C^2(\mathbb{R})$, the dynamical system $(\mathscr{H},S_t)$ associated with problem~\eqref{P} is quasi-stable.
	As a consequence, the compact global attractor $\mathfrak{A}$ has finite fractal dimension and additional smoothness properties.
	Moreover, quasi-stability implies the existence of a generalized exponential attractor.
	\subsection{Finite dimensionality, regularity, and the existence of a generalized exponential attractor}	
	We now establish the second main result of this work.
	\begin{theorem}\label{theo-main2}
		In addition to the hypotheses of Theorem~\ref{theo-main}, assume that $\Omega$ has a smooth, strictly convex boundary and:
		\begin{itemize}
			\item $f \in C^2(\mathbb{R})$ and there exists a positive constant $L$ such that
			\begin{align}\label{hyp_f''}
				|f''(s)| \le L\bigl(1 + |s|^{\rho - 2}\bigr),\quad \rho\in[2,5);
			\end{align}
			\item $g$ is nondegenerate at the origin, that is, $g(0) = \kappa > 0$.
		\end{itemize}
		Then the following results hold:
		\begin{itemize}
			\item[{\bf I.}] {\bf (Finite dimensionality).} The attractor $\mathfrak{A}$ has finite fractal dimension $\dim^{\mathscr{H}}_f\mathfrak{A}$.
			\item[{\bf II.}] {\bf (Regularity).} Any full trajectory $\Upsilon=\{U(t):t\in \mathbb{R}\}$ belonging to the global attractor satisfies the following regularity properties:
			$$U_t=(u_t,u_{tt})\in L^{\infty}(\mathbb{R};\mathscr{H}).$$
			Moreover, there exists $R>0$ such that
			\begin{align}\label{regularity-attractor}\|\nabla u_t(t)\|^2+\|u_{tt}(t)\|^2+\|\Delta u(t)\|^2\le R^2.\end{align}
			\item[{\bf III.}] {\bf (Exponential attractor).} The dynamical system $(\mathscr{H},S_t)$ possesses a generalized exponential attractor $\mathfrak{A}_{\exp}$ with finite fractal dimension in the extended space
			\[
			\mathscr{H}_{-s}=H^{-s+1}\times H^{-s},\qquad \text{for}\quad 0<s\le 1,
			\]
			where $H^{2s}=D\big((-\Delta)^s\big)$, subject to the Dirichlet boundary condition.
		\end{itemize}
	\end{theorem}
	\begin{proof}
		Proposition~\ref{stabilizability-estimate} below guarantees that the system $(\mathscr{H},S_t)$ is quasi-stable. Consequently, by a direct application of Theorem~\ref{quasi-stable-result}, the attractor $\mathfrak{A}$ has finite fractal dimension. This establishes item~{\bf I}.
		
		Moreover,  Proposition~\ref{stabilizability-estimate} and Theorem~\ref{regularity-result} imply that there exists $R_1>0$ such that
		\begin{align}\label{reg-1}
			\|\nabla u_t(t)\|^2 + \|u_{tt}(t)\|^2 \le R_1^2,
			\qquad t\in\mathbb{R}.
		\end{align}
		Next, multiplying Eq.~\eqref{P} by $-\Delta u$ and integrating over $\Omega$, we obtain
		\begin{align*}
			\|\Delta u\|^2
			+ \int_{\Omega} f'(u)|\nabla u|^2 \, dx
			= \int_{\Omega} u_{tt}\Delta u \, dx
			+ g\!\left(\|\mathbf{u}\|_{\mathscr{H}}^2\right)\!\int_{\Omega} u_t \Delta u \, dx.
		\end{align*}
		Using Remark~\ref{lowerbounded},  the H\"older inequality, and the interpolation inequality
		$\|\nabla u\|^2 \le \|u\|\,\|\Delta u\|$, we deduce
		\begin{align*}
			\|\Delta u(t)\|^2
			\le {} \|u_{tt}(t)\|\,\|\Delta u(t)\|
			+ C_{\mathcal{A}} \|u_t(t)\|\,\|\Delta u(t)\| +
			+ K_f \|u(t)\|\,\|\Delta u(t)\|,
		\end{align*}
		where $K_f$ is defined in \eqref{lower1}, and $C_{\mathcal{A}}$ is a positive constant depending only on $\mathcal{A}$.
		Consequently, by \eqref{reg-1}, there exists a constant $C_{R_1}>0$ such that
		\begin{align}\label{reg-2}
			\|\Delta u(t)\| \le C_{R_1},
			\qquad t\in\mathbb{R}.
		\end{align}
		Therefore, combining \eqref{reg-1} and \eqref{reg-2}, we obtain \eqref{regularity-attractor}.
		This completes the proof of item {\bf II}.
		
		\medskip
		Finally, item {\bf III} follows directly from Theorem~\ref{theo_attrac_exponential}. Indeed,
		let $\mathcal{B}$ be the absorbing set provided by Proposition~\ref{prop-absorbing-set}.
		According to Proposition~\ref{stabilizability-estimate}, the dynamical system $(\mathscr{H},S_t)$ is quasi-stable on $\mathcal{B}$. For $0<s\le1$, we consider the family of Hilbert spaces
		$\mathscr{H}_{-s}=H^{-s+1}\times H^{-s}$ and the map
		$[0,T]\ni t \mapsto S_t z \in \mathscr{H}=H^{1}\times H^{0}$.
		We let $\mathbf{u}(t)=(u(t), u_t(t))=S_t z$ denote the Shatah-Struwe solution associated with the initial datum $z\in\mathcal{B}$. We next show that the mapping $t \mapsto S_t z$ is H\"older continuous in $\mathscr{H}_{-s}$ for every $z \in \mathcal{B}$ and $0 < s \le 1$. We begin with the case $s = 1$.
		From Theorem~\ref{theo-global-subcritical} it follows that
		$$
		(u_t, u_{tt}) \in L^{\infty}(0,T;\mathscr{H}_{-1}),
		\qquad \text{for all}\quad T > 0,
		$$
		and
		$$
		\|S_{t_2}z - S_{t_1}z\|_{\mathscr{H}_{-1}}
		\le \int_{t_1}^{t_2} \left\|\frac{d}{ds} \mathbf{u}(s)\right\|_{\mathscr{H}_{-1}} ds
		= \int_{t_1}^{t_2}\|(u_t(s),u_{tt}(s))\|_{\mathscr{H}_{-1}} ds
		\le C_{\mathcal{B}}|t_2 - t_1|,
		$$
		for all $t_1,t_2 \in [0,T]$, where $C_{\mathcal{B}}$ is a positive constant depending on $\mathcal{B}$. Hence, the mapping $t \mapsto S_t z$ is H\"older continuous in $\mathscr{H}_{-1}$.
		
		We next consider the case $0 < s < 1$.
		Applying the interpolation theorem componentwise to the factors of $\mathscr{H}_{-s}$ and using the H\"older continuity in $\mathscr{H}_{-1}$, we obtain
		$$
		\|S_{t_2}z - S_{t_1}z\|_{\mathscr{H}_{-s}}
		\le C_s
		\|S_{t_2}z - S_{t_1}z\|_{\mathscr{H}}^{1-s}
		\|S_{t_2}z - S_{t_1}z\|_{\mathscr{H}_{-1}}^{\,s}
		\le C_{s,\mathcal{B}}\, |t_2 - t_1|^{s},
		$$
		for all $t_1,t_2 \in [0,T]$. Hence, $t \mapsto S_t z$ is
		H\"older continuous in $\mathscr{H}_{-s}$. Therefore, it follows from Theorem~\ref{theo_attrac_exponential} that the system $(\mathscr{H},S_t)$ possesses a generalized exponential attractor in $\mathscr{H}_{-s}$.
	\end{proof}
	
	\subsection{Quasi-stability}
	Finally, we establish Proposition~\ref{stabilizability-estimate}, which guarantees that the system $(\mathscr{H}, S_t)$ is quasi-stable.
	\begin{proposition}[{\bf Quasi-stability estimate}]
		\label{stabilizability-estimate}
		Under the assumptions of Theorem~\ref{theo-main2}, let $B$ be a bounded subset of $\mathscr{H}$. Let $\mathbf{u}^i(t)=S_t\mathbf{u}^i_0$ ($i=1,2$) be Shatah-Struwe solutions of~\eqref{P} with initial data $\mathbf{u}^i_0=(u^i_0,u^i_1)\in B$.
		Then the following estimate holds:
		\begin{align*}
			\|S_t\mathbf{u}^1_0-S_t\mathbf{u}^2_0\|^2_{\mathscr{H}}
			\le b(t)\|\mathbf{u}^1_0-\mathbf{u}^2_0\|^2_{\mathscr{H}}
			+ c(t)\sup_{s\in[0,t]}\|u^1(s)-u^2(s)\|^2_{L^{\frac{36}{11-\rho}}(\Omega)},
		\end{align*}
		where $b(t)$ and $c(t)$ are nonnegative scalar functions such that
		$b\in L^1(\mathbb{R}^+)$,
		$\lim_{t\to+\infty} b(t)=0$,
		and $c(t)$ is locally bounded on $[0,\infty)$.
	\end{proposition}
	\begin{proof}
		Set $z = u^1 - u^2$ and $Z = (z, z_t)$. Then $Z$ satisfies \eqref{QS-1}. For any $\alpha, \theta\in [0, 1]$, let
		\[
		\zeta_{\alpha} = \alpha \|\mathbf{u}^1\|^2_{\mathscr{H}} + (1-\alpha)\|\mathbf{u}^2\|^2_{\mathscr{H}},\quad \chi_{\theta}=\theta u^1+(1-\theta)u^2.
		\]
		Multiplying \eqref{QS-1} by $z_t$ and integrating over $\Omega \times [\tau, T]$, we obtain the following energy identity:
		\begin{align}
			\label{QS-2}
			&E_z(T)-E_z(\tau)+\int_{\tau}^TD(z(s))ds \notag\\
			&=\underbrace{K_f\int_{\tau}^T\int_{\Omega}z_tz\,dxds}_{M_1}
			+\underbrace{\frac{1}{2}\int_{\tau}^T\int_{\Omega}\int_0^1 f''(\chi_{\theta})\bigl(\theta u^1_t+(1-\theta)u^2_t\bigr)\,d\theta\,|z|^2dxds}_{M_2}\\
			&\quad\underbrace{-\frac{1}{2}\int_{\tau}^T\!\!\int_0^1 g'\bigl(\zeta_{\alpha}(s)\bigr)
				\left[\|\mathbf{u}^1(s)\|^2_{\mathscr{H}} - \|\mathbf{u}^2(s)\|^2_{\mathscr{H}}\right]
				\bigl[\|u^1_t(s)\|^2 - \|u^2_t(s)\|^2\bigr] \, d\alpha \, ds}_{M_3}, \notag
		\end{align}
		where
		\begin{align*}
			E_z(t) &:= \frac{1}{2} \|Z(t)\|^2_{\mathscr{H}}
			+ \frac{1}{2} \int_\Omega \int_0^1 f'(\chi_\theta) \, d\theta \, |z|^2 \, dx
			+ \frac{K_f}{2} \|z(t)\|^2,\\
			D(z(t)) &:= \frac{1}{2} \left[ g(\|\mathbf{u}^1(t)\|^2_{\mathscr{H}}) + g(\|\mathbf{u}^2(t)\|^2_{\mathscr{H}}) \right] \|z_t(t)\|^2,
		\end{align*}
		and the constant $K_f$ is defined in Remark~\ref{lowerbounded}.
		\medskip
		
		To prove Proposition~\ref{stabilizability-estimate}, we first establish the following equivalence.
		\begin{lemma}
			The energy $E_z$ defined above satisfies
			\begin{equation}\label{equiva}
				\frac{1}{2} \|Z(t)\|^2_{\mathscr{H}} \le E_z(t) \le C_B \|Z(t)\|^2_{\mathscr{H}},
			\end{equation}
			for some positive constant $C_B>0$.
		\end{lemma}
		
		\begin{proof}
			It follows from \eqref{lower1} that
			\begin{equation}\label{ff}
				\frac{1}{2} \int_\Omega \int_0^1 f'(\chi_\theta(t)) \, d\theta \, |z(t)|^2 \, dx
				\ge -\frac{K_f}{2} \|z(t)\|^2.
			\end{equation}
			On the other hand, from \eqref{hyp_f'}, H\"older's inequality, and the Sobolev embedding $H^1_0(\Omega) \hookrightarrow L^6(\Omega) \hookrightarrow L^{3(\rho-1)/2}(\Omega)$, together with \eqref{Estimate_I0.1}, we have
			\begin{align}\label{est2-6}
				\frac{1}{2} \int_\Omega \int_0^1 f'(\chi_\theta) \, d\theta \, |z|^2 \, dx
				\le C \Big[ 1 + \|u^1\|^{\rho-1}_{L^6(\Omega)} + \|u^2\|^{\rho-1}_{L^6(\Omega)} \Big] \|\nabla z\|^2
				\le C_B \|Z\|^2_{\mathscr{H}}.
			\end{align}
			Combining \eqref{ff} and \eqref{est2-6} with the definition of $E_z$, we obtain the equivalence \eqref{equiva}.
		\end{proof}
		
		Returning to the proof of Proposition~\ref{stabilizability-estimate}, we use the energy identity \eqref{QS-2} to  estimate the terms $M_1$, $M_2$ and $M_3$ on the right-hand side.
		By Assumption~\ref{hyp-g} and the condition $g(0)=\kappa>0$, it follows that
		\begin{align}
			\label{est-g}
			D(z(t))
			\ge
			\kappa \|z_t(t)\|^2
			+\frac{\gamma}{2}
			\left[\,
			\|\mathbf{u}^1(t)\|_{\mathscr{H}}^{2q}
			+
			\|\mathbf{u}^2(t)\|_{\mathscr{H}}^{2q}
			\,\right]
			\|z_t(t)\|^2
			\ge
			\kappa \|z_t(t)\|^2.
		\end{align}
		It follows directly from \eqref{est-g} and Young's inequality that
		\begin{eqnarray*}
			|M_1| &\le& K_f \left( \int_{\tau}^T \|z_t(s)\|^2 ds \right)^{1/2} \left( \int_{\tau}^T \|z(s)\|^2 ds \right)^{1/2} \\
			&\le& K_f \left( \frac{1}{\kappa} \int_{\tau}^T D(z(s)) ds \right)^{1/2} \left( \int_{\tau}^T \|z(s)\|^2 ds \right)^{1/2} \\
			&\le& \frac{1}{4} \int_{\tau}^T D(z(s)) ds + \frac{K_f^2 T}{\kappa} \sup_{t \in [0,T]} \|z(t)\|^2.
		\end{eqnarray*}
		Hereafter, $C_B$ denotes a generic constant that may change from line to line, depending on $B$. To estimate $M_2$, we apply \eqref{hyp_f''}, H\"older's inequality with
		$
		\frac{\rho - 2}{18} + \frac{1}{2} + \frac{11 - \rho}{18} = 1,
		$
		and \eqref{globalest} with $(r, s)=(3, 18)$, to obtain
		\begin{eqnarray*}
			|M_2| &\le&
			C \int_{\tau}^T \left[\, 1 + \|u^1(s)\|^{\rho-2}_{L^{18}(\Omega)} + \|u^2(s)\|^{\rho-2}_{L^{18}(\Omega)} \,\right]
			\left[ \|u^1_t(s)\| + \|u^2_t(s)\| \right] \|z(s)\|^2_{L^{\frac{36}{11-\rho}}(\Omega)} ds \\
			&\le& C_B \left[ T + \|u^1\|^{3}_{L^3(0,T; L^{18}(\Omega))} + \|u^2\|^{3}_{L^3(0,T; L^{18}(\Omega))} \right]\, \sup_{t \in [0,T]} \|z(t)\|^2_{L^{\frac{36}{11-\rho}}(\Omega)} \\
			&\le& C_B(T+1) \sup_{t \in [0,T]} \|z(t)\|^2_{L^{\frac{36}{11-\rho}}(\Omega)}.
		\end{eqnarray*}
		Finally, using assumption~\eqref{hyp-g}, Young's inequality, \eqref{equiva} and \eqref{est-g}, we deduce
		\begin{align*}
			|M_3|
			\le&\; \frac{\gamma q}{2} \int_{\tau}^T \int_0^1
			\bigl[\zeta_{\alpha}(s)\bigr]^{q-1} \, d\alpha \,
			(\|\mathbf{u}^1(s)\|_{\mathscr{H}} + \|\mathbf{u}^2(s)\|_{\mathscr{H}}) \, \|Z(s)\|_{\mathscr{H}} \, (\|u^1_t(s)\| + \|u^2_t(s)\|) \, \|z_t(s)\| ds \\
			\le&\; \frac{1}{4} \int_{\tau}^T D(z(s)) ds+C_B \int_{\tau}^T \varphi(s) E_z(s)ds,
		\end{align*}
		where $\varphi(t):= \|u^1_t(t)\|^2 + \|u^2_t(t)\|^2$.
		
		Thus, using the estimates $M_1$, $M_2$ and $M_3$, together with the embedding $L^{\frac{36}{11-\rho}}(\Omega) \hookrightarrow L^{2}(\Omega)$, we obtain from \eqref{QS-2} that
		\begin{eqnarray}\label{QS-3}
			\begin{aligned}
				&E_z(T)-E_z(\tau)+\frac{1}{2}\int_{\tau}^TD(z(s))ds\\
				&\le  C_{B}(T+1)\sup_{t\in [0,T]}\|z(t)\|^2_{L^{\frac{36}{11-\rho}}(\Omega)}+
				C_{B}\,\int_{0}^T\varphi(s)E_z(s)ds.
		\end{aligned}\end{eqnarray}
		Integrating \eqref{QS-3} from $0$ to $T$, we obtain
		\begin{eqnarray}
			\label{QS-33}
			\begin{aligned}
				&	TE_z(T)+\frac{1}{2}\int_0^T\int_{\tau}^TD(z(s))dsd\tau\\
				& \le \int_0^TE_z(\tau)d\tau+ C_{B}(T+1)\sup_{t\in [0,T]}\|z(t)\|^2_{L^{\frac{36}{11-\rho}}(\Omega)}+
				C_{B}T\int_{0}^T\varphi(s)E_z(s)ds.
			\end{aligned}
		\end{eqnarray}	
		Next, multiplying  Eq. \eqref{QS-1} by $z$ and integrating over $\Omega\times [0,T]$, we have
		\begin{eqnarray}\label{QS-333}
			\begin{aligned}
				2\int_{0}^TE_z(s)ds=
				&\;\underbrace{2\int_{0}^T\|z_t(s)\|^2ds+K_f\int_0^T\|z(s)\|^2ds}_{N_1}-\underbrace{\left[\int_{\Omega}z_t(s)z(s)dx\right]^{s=T}_{s=0}}_{N_2}\\
				&\;-\underbrace{\frac{1}{2}\int_0^T\left[\,g\left(||\mathbf{u}^1(t)||_{\mathscr{H}}^{2}\right)+g\left(||\mathbf{u}^2(t)||_{\mathscr{H}}^{2}\right)\,\right]\int_{\Omega}z_t z dxds}_{N_3}\\
				&\; -\underbrace{\frac{1}{2}\int_0^T\left[\,g\left(||\mathbf{u}^1(s)||_{\mathscr{H}}^{2}\right)-g\left(||\mathbf{u}^2(s)||_{\mathscr{H}}^{2}\right)\,\right]\int_{\Omega}\left(\,u^1_t+u^2_t\,\right)zdxds}_{N_4}.\end{aligned}
		\end{eqnarray}
		From \eqref{est-g} and the embedding $L^{\frac{36}{11-\rho}}(\Omega)\hookrightarrow L^2(\Omega)$, we get
		\begin{align*}
			|N_1|\le \frac{2}{\kappa} \int_0^TD(z(s))ds+CT\sup_{t\in[0,T]}\|z(t)\|^2_{L^{\frac{36}{11-\rho}}(\Omega)}.
		\end{align*}
		From  \eqref{equiva}, we obtain
		\begin{align}
			\label{Asymptotic-2}
			|N_2|\le \frac{1}{\lambda_1^{1/2}}\left[\,E_z(T)+E_z(0)\,\right].
		\end{align}
		Moreover, using the energy identity \eqref{QS-2} and the estimates for $M_1,\,M_2$, and $M_3$ (with $\tau=0$), together with the embedding $L^{\frac{36}{11-\rho}}(\Omega) \hookrightarrow L^{2}(\Omega)$, we obtain
		\begin{align*}
			E_z(0)\le E_z(T)+\frac{3}{2}\int_0^TD(z(s))ds+C_{B}(T+1)\sup_{t\in [0,T]}\|z(t)\|^2_{L^{\frac{36}{11-\rho}}(\Omega)}+C_{B}\int_{0}^T\varphi(s)E_z(s)ds.
		\end{align*}
		Hence, substituting the above inequality into \eqref{Asymptotic-2}, we obtain
		\begin{eqnarray*}
			|N_2|&\le& \frac{2}{\lambda_1^{1/2}}E_z(T)+\frac{3}{2\lambda_1^{1/2}}\int_0^TD(z(s))ds+C_{B}(T+1)\sup_{t \in [0,T]} \|z(t)\|^2_{L^{\frac{36}{11-\rho}}(\Omega)}\\
			&&+\,C_{B}\int_{0}^T\varphi(s)E_z(s)ds.
		\end{eqnarray*}	
		Using the continuity of $g$, \eqref{est-g} and the embedding $L^{\frac{36}{11-\rho}}(\Omega)\hookrightarrow L^2(\Omega)$, we get
		\begin{eqnarray*}
			|N_3|
			\le  C_B\int_0^T\|z_t(s)\|\|z(s)\|ds
			\le \int_0^TD(z(s))ds + C_{B}T\sup_{t \in [0,T]} \|z(t)\|^2_{L^{\frac{36}{11-\rho}}(\Omega)},
		\end{eqnarray*}
		and
		\begin{eqnarray*}
			|N_4|
			&\le& C_B\int_0^T
			\|Z(s)\|_{\mathscr{H}}\left[\|u^1_t(s)\|+\|u^2_t(s)\|\right]\|z\|ds\\
			&\le &C_B\int_0^T\varphi(s)E_z(s)ds+C_BT\sup_{t\in[0,T]}\|z(t)\|^2_{L^{\frac{36}{11-\rho}}(\Omega)}.
		\end{eqnarray*}
		Returning to~\eqref{QS-333}, and using~\eqref{QS-3} together with the estimates for the terms \(N_1, N_2, N_3\), and $N_4$, we deduce
		\begin{eqnarray}\label{QS-4}
			\begin{aligned}
				\int_{0}^TE_z(s)ds\le&\; \frac{1}{\lambda_1^{1/2}}E_z(T)+C\int_0^TD(z(s))ds\\
				&+C_{B}(T+1)\sup_{t\in[0,T]}\|z(t)\|^2_{L^{\frac{36}{11-\rho}}(\Omega)}+C_B\int_0^T\varphi(s)E_z(s)ds,\end{aligned}
		\end{eqnarray}	
		where $C=\frac{2}{\kappa}+\frac{3}{2\lambda_1^{1/2}}+1$.
		
		Combining \eqref{QS-33} and \eqref{QS-4} we obtain that
		\begin{eqnarray}
			\label{QS-33'}
			\begin{aligned}
				\left(T-\frac{1}{\lambda_1^{1/2}}\right)E_z(T)&+\frac{1}{2}\int_0^T\int_{\tau}^TD(z(s))dsd\tau\\
				& \le  C\int_0^TD(z(s))ds+C_{B}(T+1)\sup_{t\in [0,T]}\|z(t)\|^2_{L^{\frac{36}{11-\rho}}(\Omega)}\\
				&+C_{B}(T+1)\int_{0}^T\varphi(s)E_z(s)ds.
			\end{aligned}
		\end{eqnarray}	
		From \eqref{QS-3} with $\tau=0$, we get
		\begin{eqnarray}
			\label{QS-03}
			\begin{aligned}
				\int_{0}^TD(z(s))ds\le &\;2\left[\,E_z(0)-E_z(T)\,\right]+  C_{B}(T+1)\sup_{t \in [0,T]} \|z(t)\|^2_{L^{\frac{36}{11-\rho}}(\Omega)}\\
				&+
				C_{B}\,\int_{0}^T\varphi(s)E_z(s)ds.
			\end{aligned}
		\end{eqnarray}	
		Substituting \eqref{QS-03} in \eqref{QS-33'}, we have
		\begin{eqnarray*}
			\left(T-\frac{1}{\lambda_1^{1/2}}\right)E_z(T) &\le &2C\left[\,E_z(0)-E_z(T)\,\right]+C_{B}(T+1)\sup_{t\in [0,T]}\|z(t)\|^2_{L^{\frac{36}{11-\rho}}(\Omega)}\\
			&&\;+
			C_{B}(T+1)\int_{0}^T\varphi(s)E_z(s)ds.
		\end{eqnarray*}	
		Without loss of generality, we may choose $T>0$ such that $\frac{T}{2}-\frac{1}{\lambda_1^{1/2}}>0$. Then
		\begin{align*}
			E_z(T) \le \eta E_z(0)+C_{B}\left[\,\sup_{t\in [0,T]}\|z(s)\|^2_{L^{\frac{36}{11-\rho}}(\Omega)}+\int_{0}^T\varphi(s)E_z(s)ds\right],
		\end{align*}
		for $C_B\ge\frac{C_B(T+1)}{T+4C}$, where $\eta=\frac{4C}{T+4C}$.
		
		Reiterating the estimate on the intervals $[mT,(m+1)T]$ yields
		\begin{align*}
			E_z((m+1)T)\le\eta E_z(mT)+C_{B}b_{m},\quad m=0,1,2,\cdots,
		\end{align*}
		where
		$$0<\eta<1\quad\mbox{and}\quad b_m\equiv \sup_{t\in [mT,(m+1)T]}\|z(t)\|^2_{L^{\frac{36}{11-\rho}}(\Omega)}+\int_{mT}^{(m+1)T}\varphi(s)E_z(s)ds.$$
		This yields
		$$E_z(mT)\le\eta^{m}E_z(0)+C_B\sum_{l=1}^m\gamma^{m-l}b_{l-1}.$$
		Since $\eta<1$,  arguing as in \cite[Remark 3.30]{chueshov-lasiecka-2008}, along with the definition of $b_m$ we obtain that there exists $\varrho>0$ such that
		\begin{align*}
			E_z(t)\le C_1e^{-\varrho t}E_z(0)+C_2\left[\,\sup_{s\in [0,t]}\|z(s)\|^2_{L^{\frac{36}{11-\rho}}(\Omega)}+\int_0^te^{-\varrho(t-s)}\varphi(s)E_z(s)ds\,\right],\quad \forall t\ge 0.
		\end{align*}
		Therefore, applying Gronwall's lemma, we find
		\begin{align*}
			E_z(t)\le \left[\,C_1e^{-\varrho t}E_z(0)+C_2\sup_{s\in [0,t]}\|z(s)\|^2_{L^{\frac{36}{11-\rho}}(\Omega)}\,\right]e^{C_2\int_0^t\varphi(s)ds}.
		\end{align*}
		Using the equivalence \eqref{equiva},
		we obtain
		\begin{align*}
			\|\mathbf{u}^1(t) - \mathbf{u}^2(t)\|^2_{\mathscr{H}}
			\le b(t) \|\mathbf{u}^1_0 - \mathbf{u}^2_0\|^2_{\mathscr{H}}
			+ c(t) \sup_{s \in [0,t]} \|z(s)\|^2_{L^{\frac{36}{11-\rho}}(\Omega)},
		\end{align*}
		where the embedding $H^1_0(\Omega) \hookrightarrow\hookrightarrow L^{\frac{36}{11-\rho}}(\Omega)$ is compact, and
		\begin{align*}
			b(t) := 2 C_B C_1 e^{-\varrho t} e^{C_2 \int_0^t \varphi(s) \, ds},
			\quad
			c(t) := 2 C_2 e^{C_2 \int_0^t \varphi(s) \, ds}.
		\end{align*}
		By \eqref{est-u-t-b}, we have
		\[
		\varphi(s) = \|u^1_t(s)\|^2 + \|u^2_t(s)\|^2 \in L^1(0,t) 
		\]
		uniformly with respect to  $t>0$. It follows that $b(t) \in L^1(\mathbb{R}^+)$ with $\lim_{t \to +\infty} b(t) = 0$, while $c(t)$ is bounded on $\mathbb{R}^+$ due to the $L^1$-integrability of $\varphi(t)$. Therefore, the proof of Proposition~\ref{stabilizability-estimate} is complete.
	\end{proof}
	
	\appendix
	\section{Abstract Results on Attractors for Dynamical Systems}
	
	This appendix collects fundamental results from the theory of continuous dynamical systems concerning the existence, characterization, and qualitative properties of global attractors. The relevant theoretical framework is developed in the monographs by Chueshov and Lasiecka \cite{chueshov-lasiecka-2008,chueshov-lasiecka-2010}.
	
	\subsection{Existence and characterization of global attractors}
	
	\begin{definition}
		A \textbf{dynamical system} is a pair $(X,S_t)$, where $X$ is a complete metric space and $S_t:X \to X$ is a family of continuous maps such that $y\mapsto S_ty$ is continuous in $X$ for every $y\in X$, and the semigroup property holds:
		\[
		S_0=I, \quad S_{t+s}=S_t\circ S_s, \quad \forall\, s,t\ge 0.
		\]
		$S_t$ is called the \textbf{evolution semigroup} or \textbf{evolution operator}.
	\end{definition}
	
	\begin{definition}
		Let $(X,S_t)$ be a dynamical system. A \textbf{global attractor} $\mathfrak{A} \subset X$ is a compact set that is fully invariant and uniformly attracting:
		\[
		S_t\mathfrak{A} = \mathfrak{A} \quad \forall t\ge 0,
		\]
		and for every bounded set $B \subset X$,
		\[
		\mbox{dist}_X(S_tB,\mathfrak{A}) := \sup_{x\in S_tB} \inf_{y\in \mathfrak{A}} d_X(x,y) \to 0 \quad \text{as } t\to\infty.
		\]
	\end{definition}
	
	\begin{definition}
		A bounded set $\mathcal{B} \subset X$ is called an \textbf{absorbing set} if, for every bounded set $B\subset X$, there exists $t_B \ge 0$ such that
		\[
		S_tB \subset \mathcal{B}, \quad \forall t \ge t_B.
		\]
		A dynamical system admitting a bounded absorbing set is called \textbf{dissipative}.
	\end{definition}
	
	\begin{definition}
		A dynamical system $(X,S_t)$ is \textbf{asymptotically smooth} in $X$ if, for any bounded positively invariant set $B \subset X$, there exists a compact set $K \subset \overline{B}$ such that
		\[
		\mbox{dist}_X(S_tB,K) \to 0 \quad \text{as } t\to \infty.
		\]
	\end{definition}
	
	\begin{definition}
		Let $\mathcal{N}$ denote the set of stationary points of $(X,S_t)$:
		\[
		\mathcal{N} = \{v \in X : S_tv = v \text{ for all } t\ge 0\}.
		\]
		The \textbf{unstable manifold} $\mathrm{M}^u(\mathcal{N})$ emanating from $\mathcal{N}$ is defined as the set of all $y\in X$ for which there exists a full trajectory $\Upsilon = \{u(t): t\in \mathbb{R}\}$ satisfying
		\[
		u(0)=y, \quad \lim_{t\to -\infty} \mbox{dist}_X(u(t),\mathcal{N}) = 0.
		\]
	\end{definition}
	
	\begin{theorem}[Theorem 2.3, \cite{chueshov-lasiecka-2008}] \label{theo-chueshov}
		Let $(X,S_t)$ be a dissipative dynamical system on a complete metric space $X$. Then $(X,S_t)$ possesses a compact global attractor $\mathfrak{A}$ if and only if $(X,S_t)$ is asymptotically smooth.
	\end{theorem}
	
	\begin{theorem}[Theorem 7.5.6, \cite{chueshov-lasiecka-2010}] \label{theo-chueshov-caracterization}
		Let $(X,S_t)$ possess a compact global attractor $\mathfrak{A}$. Assume that there exists a strict Lyapunov function on $\mathfrak{A}$. Then
		\[
		\mathfrak{A} = \mathcal{M}^u(\mathcal{N}),
		\]
		and $\mathfrak{A}$ consists of full trajectories $\Upsilon = \{z(t): t\in \mathbb{R}\}$ such that
		\[
		\lim_{t\to -\infty} \mbox{dist}_X(z(t),\mathcal{N}) = 0, \quad
		\lim_{t\to +\infty} \mbox{dist}_X(z(t),\mathcal{N}) = 0.
		\]
	\end{theorem}
	
	\subsection{Finite dimension, regularity, and exponential attractors}
	
	\begin{definition}
		The \textbf{fractal dimension} of a compact set $K$ in a metric space $X$ is
		\[
		\mbox{dim}_f^X K = \limsup_{\varepsilon \to 0} \frac{\ln n(K,\varepsilon)}{\ln(1/\varepsilon)},
		\]
		where $n(K,\varepsilon)$ is the minimal number of closed balls of radius $\varepsilon$ covering $K$.
	\end{definition}
	
	\begin{definition}
		A compact set $\mathfrak{A}_{\exp} \subset X$ is a \textbf{fractal exponential attractor} if it is positively invariant, has finite fractal dimension, and for every bounded set $B\subset X$, there exist $t_B, C_B, \sigma_B > 0$ such that
		\[
		\sup_{z\in B} \mbox{dist}_X(S(t)z, \mathfrak{A}_{\exp}) \le C_B e^{-\sigma_B (t-t_B)}, \quad t \ge t_B.
		\]
		If such a set is defined in an extended space $\widetilde{X} \supseteq X$, it is called a \textbf{generalized fractal exponential attractor}.
	\end{definition}
	
	\begin{definition} \label{Def-quasi-stable}
		Let $Y, Z, W$ be reflexive Banach spaces with $Y \hookrightarrow\hookrightarrow Z$ and set $X = Y \times Z \times W$, allowing $W = \{0\}$. Consider $(X,S_t)$ with evolution operator
		\[
		S_t z = (u(t), u_t(t), \theta(t)), \quad z=(u_0,u_1,\theta_0)\in X,
		\]
		where
		\[
		u \in C(\mathbb{R}^+; Y) \cap C^1(\mathbb{R}^+; Z), \quad
		\theta \in C(\mathbb{R}^+; W).
		\]
		The system is said to be  \textbf{quasi-stable} on $B\subset X$ if there exist a compact seminorm $n_Y$ on $Y$ and nonnegative locally bounded functions $a(t)$, $c(t)$  and $b(t)\in L^1(\mathbb{R}^+)$ with $b(t)\to 0$ as $t\to\infty$, such that for any $z^1,z^2\in B$,
		\begin{align*}
			\Vert S_t z^1 - S_t z^2 \Vert_X^2 &\le a(t) \Vert z^1 - z^2 \Vert_X^2, \\
			\Vert S_t z^1 - S_t z^2 \Vert_X^2 &\le b(t) \Vert z^1 - z^2 \Vert_X^2 + c(t) \sup_{s\in[0,t]} [n_Y(u^1(s)-u^2(s))]^2.
		\end{align*}
		The second inequality is called the \textbf{stabilizability inequality}.
	\end{definition}
	
	\begin{theorem}[Theorem 7.9.6, \cite{chueshov-lasiecka-2010}] \label{quasi-stable-result}
		If $(X,S_t)$ possesses a compact global attractor $\mathfrak{A}$ and is quasi-stable on $\mathfrak{A}$, then $\mathfrak{A}$ has finite fractal dimension $\dim_f^X \mathfrak{A}$.
	\end{theorem}
	
	\begin{theorem}[Theorem 7.9.8, \cite{chueshov-lasiecka-2010}] \label{regularity-result}
		Under the assumptions of Theorem \ref{quasi-stable-result}, if $c_\infty := \sup_{t\ge 0} c(t) < \infty$, then any full trajectory $(u(t), u_t(t), \theta(t))$ in $\mathfrak{A}$ satisfies
		\[
		u_t \in L^\infty(\mathbb{R};Y) \cap C(\mathbb{R};Z), \quad
		u_{tt} \in L^\infty(\mathbb{R};W), \quad
		\theta \in L^\infty(\mathbb{R};W),
		\]
		and there exists $R>0$ (depending on $c_\infty$, the semigroup, and the embedding $Y \hookrightarrow Z$) such that
		\[
		\|u_t(t)\|_Y^2 + \|u_{tt}(t)\|_Z^2 + \|\theta_t(t)\|_W^2 \le R^2, \quad t\in\mathbb{R}.
		\]
	\end{theorem}
	
	\begin{theorem}[Theorem 7.9.9, \cite{chueshov-lasiecka-2010}] \label{theo_attrac_exponential}
		Let $(X,S_t)$ satisfy the assumptions of Definition \ref{Def-quasi-stable}, be dissipative, and quasi-stable on some bounded absorbing set $\mathcal{B}$. Assume there exists $\widetilde{X} \supseteq X$ such that $t \mapsto S(t)z$ is H\"older continuous in $\widetilde{X}$ for each $z\in \mathcal{B}$:
		\[
		\|S_{t_2} z - S_{t_1} z\|_{\widetilde{X}} \le C_{\mathcal{B},T} |t_2 - t_1|^\alpha, \quad t_1,t_2 \in [0,T], \ z\in \mathcal{B},
		\]
		for some $0<\alpha\le 1$. Then $(X,S_t)$ possesses a (generalized) fractal exponential attractor with finite dimension in $\widetilde{X}$.
	\end{theorem}

	\section*{Acknowledgements}
	The authors are indebted  to Nicolas Burq for  his valuable comments and  for  bringing to our attention the work of Ivanovici, Lascar, Lebeau, and Planchon \cite{Ivanovici-Lascar-Lebeau-Planchon-2023}  which is critical  in the present proof of quasi-stability.

	\section*{Declarations}
	
	\subsubsection*{\bf Ethical Approval}
	Not applicable to this article.
	
	\subsubsection*{\bf Availability of data and materials}
	Data availability does not apply to this article, as no datasets were generated or analyzed and no materials were used during the current study.
	
	\subsubsection*{\bf Conflict of interest}
	On behalf of all authors, the corresponding author declares that there is no conflict of interest.

\end{document}